\documentclass[7pt]{article}
\usepackage{graphicx} 
\usepackage[colorlinks=True, linkcolor=blue, citecolor=red]{hyperref}
\usepackage{subfigure}
\usepackage{fancyvrb}
\usepackage{amsmath, amsfonts, amsthm, amssymb, bbm}
\usepackage{stmaryrd}
\usepackage{changepage}
\usepackage[title,toc,titletoc]{appendix}
\usepackage{verbatim}
\usepackage{graphicx}
\usepackage{nicematrix}
\usepackage{setspace}
\usepackage{multirow}
\usepackage{float}
\usepackage{comment}
\usepackage{enumitem}
\usepackage{xcolor}
\usepackage[all]{xy}
\usepackage{lipsum} 
\usepackage{mathtools}
\usepackage{mathrsfs}
\usepackage{shuffle}
\usepackage{enumitem}
\usepackage{csquotes}
\usepackage{authblk}
\usepackage[pagewise]{lineno}

\allowdisplaybreaks

\usepackage[a4paper,
            bindingoffset=0.2in,
            left=0.5in,
            right=0.75in,
            top=1in,
            bottom=1in,
            footskip=.25in]{geometry}

\newcommand{\floor}[1]{\left \lfloor {#1}\right \rfloor}

\newtheorem{theorem}{Theorem}[section]
\newtheorem{corollary}{Corollary}[theorem]
\newtheorem{lemma}[theorem]{Lemma}
\newtheorem{definition}[theorem]{Definition}
\newtheorem{proposition}[theorem]{Proposition}
\newtheorem{remark}[theorem]{Remark}
\theoremstyle{remark}
\newtheorem{example}[theorem]{Example}

\author{Thomas Cass\thanks{Thomas Cass has been supported by the EPSRC Programme Grant EP/S026347/1 and acknowledges the support of the Erik Ellentuck Fellowship at the Institute for Advanced Study.} , Dan Crisan, and Andrea Iannucci \thanks{Corresponding author: \texttt{andrea.iannucci22@imperial.ac.uk}}  }

\affil{Department of Mathematics, Imperial College London}

\title{Pathwise Optimal Control and Rough Fractional Hamilton-Jacobi-Bellman Equations for Rough-Fractional Dynamics}
\begin{document}

\setlength{\abovedisplayshortskip}{0pt}
\setlength{\belowdisplayshortskip}{0pt}

\maketitle

\begin{abstract}
In this work, we investigate the degeneracy problem in pathwise control,
extending the framework developed in \cite{allan2020pathwise} to a more
general class of driving signals and a broader set of admissible controls.
Our approach consists in choosing admissible controls from a suitable class of
derivatives and transforms the original control equation into a fractional
dynamics system. Within this setting, we derive sufficient conditions ensuring
that the control problem remains non-degenerate. We then build on the analysis developed in
\cite{gomoyunov2020dynamic,gomoyunov2020theory,gomoyunov2021viscosity}
to study the resulting value function and the associated
Hamilton--Jacobi--Bellman equation.
\end{abstract}

\section{Introduction}
Control problems driven by irregular signals arise naturally when one seeks a
pathwise formulation of stochastic control problems. Rough path theory provides
a robust framework for treating such dynamics, since it allows one to define
controlled differential equations for fixed driving signals of low regularity.
In this paper we work in this pathwise setting and consider controlled state dynamics
of the form,
\begin{equation}\label{controlled_equation_introduction}
    dX^{x,\gamma}_t
    =
    b(X^{x,\gamma}_t,\gamma_t)\,dt
    +
    \lambda(X^{x,\gamma},\gamma)_t\,d\boldsymbol{\zeta}_t,
    \qquad
    X^{x,\gamma}_0=x,
    \qquad
    \gamma\in\mathcal A,
\end{equation}
where \(\mathbf{\zeta}\) is a given rough path and \(\gamma\) ranges over a
class of admissible controls $\mathcal{A}$. Given a running cost \(f\) and a terminal cost \(g\), the associated cost
functional is
\[
    J^{\boldsymbol{\zeta}}(t,x,\gamma)
    =
    \int_t^T f(X^{x,\gamma}_r,\gamma_r)\,dr
    +
    g(X^{x,\gamma}_T),
\]
and the pathwise value function is defined by
\begin{equation}\label{pathwise_control_problem}
    v^{\boldsymbol{\zeta}}(t,x)
    =
    \inf_{\gamma\in\mathcal A}
    J^{\boldsymbol{\zeta}}(t,x,\gamma).
\end{equation}
The optimisation is therefore performed after the driving rough path has been
fixed. In particular, once \(\boldsymbol{\zeta}\) is prescribed, the problem is
deterministic.\newline

Control problems driven by rough signals have been studied from several
perspectives. In the case where the diffusion coefficient is not controlled,
rough-path methods provide a robust framework for formulating pathwise
Hamilton--Jacobi--Bellman equations and dynamic programming principles; see,
for instance, \cite{caruana2011rough}. In this setting the
rough path lift of the driving signal allows one to give a pathwise meaning to
the dynamics and to the associated value function.

The situation changes substantially when the diffusion coefficient
\(\lambda\) is allowed to depend on the control. In that case, the irregularity
of the driving signal can be exploited by the control, and the resulting
pathwise problem may become degenerate. This phenomenon was analysed by Allan
and Cohen \cite{allan2020pathwise}, who showed that the degeneracy can be
removed by imposing suitable regularity on the controls and by adding an
appropriate penalisation to the running cost. More precisely, their approach
restricts the admissible controls to a Sobolev-type class and penalises the
weak derivative of the control.\newline

As mentioned above, the interest in such pathwise formulations is also motivated by stochastic
control. If the rough signal \(\mathbf{\zeta}\) is obtained as the lift of a
stochastic process, then one may compare the usual stochastic value function
\[
    \overline v(t,x)
    =
    \inf_{\gamma\in\mathcal A}
    \mathbb E\left[J(t,x,\gamma)\right]
\]
with the pathwise-optimised quantity
\[
    v(t,x)
    =
    \mathbb E\left[
        \inf_{\gamma\in\mathcal A}
        J(t,x,\gamma)
    \right].
\]
The relation between these two formulations has a long history. Wets
\cite{wets1975relation} obtained an early result linking stochastic and
pathwise optimisation problems, up to the nonanticipativity constraint on the
controls. Haussmann \cite{ha1992deterministic}, building on the flow
decomposition of anticipating stochastic differential equations due to Ocone
and Pardoux \cite{ocone1989generalized}, showed how stochastic control problems
can be decomposed into deterministic problems indexed by the realisation of the
noise. Related ideas were discussed by Lions and Souganidis
\cite{lions1998fully}, who connected these decompositions with
Hamilton--Jacobi--Bellman equations, and this connection was later developed by
Buckdahn and Ma \cite{buckdahn2007pathwise}. Pathwise dual formulations also
appear in the work of Rogers on optimal stopping and discrete-time Markov
control problems \cite{rogers2002monte,rogers2007pathwise}.\newline

These results provide part of the motivation for studying pathwise optimisation
problems in their own right. In particular, they show that such problems arise
naturally when stochastic control problems are decomposed according to the
realisation of the driving noise. This leads to the pathwise viewpoint adopted
here, where the driving signal is fixed throughout the optimisation and the
resulting problem is analysed as a deterministic control problem driven by a
rough path.\newline

The main difficulty, when the diffusion coefficient depends on the control, is
the degeneracy phenomenon described above.
In this work, we further explore the degeneracy problem by building on the framework developed in \cite{allan2020pathwise}, extending it to encompass a broader class of driving paths and a wider set of admissible controls. The admissible controls are selected from the space
\(AC^\alpha([0,T],\mathbb R^\ell)\) of fractional absolutely continuous
paths, consisting of all paths of the form
\(
    \gamma_t
    =
    \gamma_0 + I_{0+}^\alpha u(t),
    u\in L^\infty([0,T],\mathbb R^\ell)
\) and $I^\alpha_{0^{+}}$ is the left-sided fractional integral of order $\alpha$ with base point $0$.
These paths are H\"older continuous, and the function
\(
    u
    =
    D_{0+}^\alpha(\gamma-\gamma_0)
\)
is referred to as the pseudo-control. As a result, the controlled process is transformed from \eqref{controlled_equation_introduction} into:
    
    \begin{equation} 
    \label{fractional_dynamics_introduction}
    \begin{aligned}
        &dX^{x, a, u}_t = b(X^{x, a, u}_t, \gamma^{a, u}_t) \, dt + \lambda(X^{x, a, u}, \gamma^{a, u})_t \, d\boldsymbol{\zeta}_t,  \quad X^{x, a, u}_0 = x,  \\
        &D^{\alpha}_{0^+}(\gamma^{a, u} - a)(t) = u_t,  \quad \gamma^{a,u}_0 = a.
    \end{aligned}
    \end{equation}

    The rest of the paper is organised as follows. Section \ref{section2} recalls basic results from the theory of weakly geometric rough paths, which provide the analytic foundation for our work. Section \ref{section_rough} introduces the rough control problem underlying the paper. 
    We formulate the controlled rough differential equation and derive estimates for the corresponding state process in terms of the control $\gamma$. These estimates identify the mechanism by which the value functional can become degenerate, and motivate the need for an additional penalisation of the control. Following the philosophy of \cite{allan2020pathwise}, this leads us to introduce a pseudo-control and to modify the cost functional by an additional penalisation term.
    
    In contrast with the Sobolev formulation of \cite{allan2020pathwise}, the pseudo-control used here is the fractional derivative of the control. This leads to an extended fractional dynamics: the control path is reconstructed from its fractional derivative via the Riemann--Liouville integral. This enlargement of the dynamics is a central step in the paper. It permits a broader class of admissible controls while retaining the  structure needed for dynamic programming and the subsequent Hamilton--Jacobi--Bellman analysis.
    
    Section~\ref{section_recovering} develops this fractional formulation.
Motivated by the estimates obtained in the previous section, we introduce a
functional $\tilde J$ depending on the pseudo-control $u$, which is added to the
original cost functional $J$. We then address the nonlocal nature of fractional
integration by defining a suitable extension of past controls by future
pseudo-controls. This yields the fractional cost and value functionals used
throughout the rest of the paper. The section concludes by proving that the
resulting control problem is well posed and that the corresponding value
functional is bounded. Section \ref{section_properties} studies the structural properties of the fractional value functional. In particular, we prove the Dynamic Programming Principle in the extended path-dependent formulation. We then establish continuity of the value functional with respect to its variables, where the nonlocality of the fractional derivative requires careful treatment of the concatenation of past and future controls. These results provide the regularity and stability properties needed to pass from the control problem to its Hamilton--Jacobi--Bellman characterisation.
    
    Finally, Section \ref{section_well} carries out this passage. The equation obtained from the DPP is both fractional, in the sense of Gomoyunov's path-dependent HJB theory \cite{gomoyunov2020dynamic, gomoyunov2020theory, gomoyunov2021viscosity}, and rough, through its dependence on the driving rough path. We therefore combine Gomoyunov's fractional viscosity framework with the
stability approach to rough viscosity solutions of
\cite{caruana2011rough}. For a fixed smooth driving path, viscosity solutions
of the corresponding fractional HJB equation are defined through the
truncated Hamiltonians \(H_n\): the limiting functional is a viscosity
solution if it is the locally uniform limit of viscosity solutions of the
\(n\)-truncated problems. For a rough driving path
\(\boldsymbol{\zeta}\), we then define the corresponding rough viscosity
solution through smooth approximation of the driver. More precisely, for
every sequence of smooth paths whose canonical lifts converge to
\(\boldsymbol{\zeta}\) in the \(p\)-variation topology, the associated
smooth-driver value functionals converge locally uniformly to
\(v^{\boldsymbol{\zeta}}\) on
\([0,T]\times \mathbb{R}^e\times AC^\alpha\).

With these notions, we prove that the value functional
\(v^{\boldsymbol{\zeta}}\) is a rough viscosity solution of the rough
fractional HJB equation. Moreover, for every smooth approximation of the
driving rough path, the corresponding value functional is the unique
viscosity solution of the approximating fractional HJB equation within the
class of viscosity solutions whose truncated approximating sequences
satisfy property~{\rm (L)}.

The following theorem collects, under a common set of standing
assumptions, the main results proved throughout the paper.
\begin{theorem}[Summary of the main results]
\label{thm:main_result}
In the setting described above, let \(p\in[2,\infty)\setminus\mathbb N\), and let
\(\boldsymbol{\zeta}\) be a \(p\)-rough path over
\(\mathbb R^d\). Suppose that
\[
    \frac{\lfloor p\rfloor}{p}<\alpha<1,
    \quad
    1<\kappa<
    \frac{1}{
        1-\alpha+\lfloor p\rfloor/p
    },\quad      q>
    \lfloor p\rfloor p
    \vee
    \frac{\kappa}{\kappa-1}.
\]
The admissible controls are
\(\gamma\in AC^\alpha([0,T],\mathbb R^\ell)\), with pseudo-control
\(
    u:=D_{0+}^\alpha(\gamma-\gamma_0).
\)
Assume that
\[
    b\in
    \operatorname{Lip}_b
    (\mathbb R^e\times\mathbb R^\ell,\mathbb R^e),
    \qquad
    \lambda\in
    C_b^{\floor p +1}
    \bigl(
        \mathbb R^e\times\mathbb R^\ell,
        \mathcal L(\mathbb R^d,\mathbb R^e)
    \bigr).
\]
Here, \(\mathcal L(\mathbb R^d,\mathbb R^e)\) denotes the space of
linear maps. Assume moreover that
\[
    f
    \in
    \operatorname{Lip}_b
    \bigl(
        \mathbb R^e\times\mathbb R^\ell,
        \mathbb R
    \bigr),
    \qquad
    g
    \in
    \operatorname{Lip}
    \bigl(
        \mathbb R^e,
        \mathbb R
    \bigr).
\]
Let the additional pseudo-control penalisation be given by
\[
    \widetilde J(t,u)
    :=
    \int_t^T \widetilde f(u_r)\,dr,
\]
where
\(\widetilde f:\mathbb R^\ell\to\mathbb R\)
is continuous and satisfies
\(
    \widetilde f(u)
    \geq
    f_0|u|^q-C_0
\)
for some constants \(f_0,C_0>0\).

Then the associated fractional control problem with total cost
\(J+\widetilde J\) is well posed. More precisely, for every admissible
initial condition and pseudo-control, the extended fractional controlled
rough differential equation admits a unique solution, and the
corresponding cost functional is well defined.

The associated value functional
\[
    v:
    [0,T]\times\mathbb R^e
    \times AC^\alpha([0,T],\mathbb R^\ell)
    \to \mathbb R
\]
is locally bounded and continuous with respect to the product topology
induced by
\[
    \varrho
    \bigl(
        (t,x,\gamma),(s,y,\nu)
    \bigr)
    :=
    \max
    \left\{
        |t-s|,
        |x-y|, |\gamma_0 - \nu_0| + 
        \|\gamma-\nu\|_{p / \floor{p}-\text{var};[0,T]}
    \right\}.
\]
Moreover, \(v\) satisfies the Dynamic Programming Principle in the
extended path-dependent state space.

Finally, \(v\) is a rough viscosity solution of an
associated rough fractional Hamilton--Jacobi--Bellman equation.

The precise formulations of these conclusions are given in
Proposition~\ref{boundedness_val_fun}, Lemma~\ref{dpp_principle},
Propositions~\ref{regularity_value_functional} and~\ref{continuity_in_noise},
Theorem~\ref{comparison_thm}, Corollary~\ref{uniqueness_thm_viscosity},
and Theorem~\ref{final_theorem_viscosity}.
\end{theorem}

\subsection{A toy model of insider trading with exogenous market impact} \label{example_subsection}
We consider a motivating example, inspired by insider trading under exogenous
market perturbations, which illustrates the degeneracy problem discussed above arising in a
simple control problem driven by a path of unbounded variation. Let the time horizon be \([0,T]\) for some $T>0$. The \enquote{unperturbed price process} \(S_t\) is modelled as a simple stochastic process driven by Brownian motion:
\begin{align*}
dS_t &= \sigma \, dB_t, \\
S_0 &= s_0,
\end{align*}
where \(B\) is a standard Brownian motion under the probability measure \(\mathbb{P}\), and \(\sigma > 0\) is a volatility parameter. Here, \enquote{unperturbed} means that \(S_t\) evolves without any external interventions or additional market effects.

We introduce an \enquote{exogenous market impact} modeled by a linear propagator:
\[
I_t = \int_0^t k(t-r) \, dW_r,
\]
where \(k(t) = \sigma C  t^{-\rho}\) for some constants \(C>0\) and \(\rho \in (\frac{1}{6}, \frac{1}{2})\), and \(W\) is a standard Brownian motion independent of \(B\). This process \(I\) represents an external perturbation to the price,
e.g., from other market participants, and is independent of the
trader's actions. It admits a modification whose paths are
\(\beta\)-H\"older continuous for every
\(\beta<\frac12-\rho\), and therefore have finite \(p\)-variation for
every
\(
    p>1 / ({\frac12-\rho})
\). Such a propagator kernel is used in \cite{gatheral2010no, bucci2018slow, webster2023handbook,bouchaud2003fluctuations} and \cite{abi2022optimal}.

We now consider a trader who observes the full path of the combined process \(P_t := B_t + I_t\). Under the assumption of no additional frictions, the trader's wealth process \(X_t\) is influenced by a control \(\gamma_t\), representing the trader's chosen trading rate. Formally, the trader seeks to maximise
\[
v(t, x) = \sup_{\gamma \in B_M} X^{x, \gamma}_T,
\]
where \(B_M\) denotes the set of bounded measurable functions with supremum norm bounded by \(M\), and the wealth process \(X_t^{x, \gamma}\) is heuristically written as
\begin{align*}
dX^{x, \gamma}_t &= \sigma \, \gamma_t \, dP_t, \\
X_0 &= x.
\end{align*}

The expression above is not well-defined pathwise for arbitrary bounded
controls \(\gamma\). A natural way to make the problem rigorous is to first
consider discrete-time approximations, where the wealth process is defined
along a fixed partition. However, this also reveals the basic obstruction:
as the mesh of the partition tends to zero, the corresponding control problem
may become degenerate.

Indeed, let \((\mathcal P_n)_{n \geq 1}\) be a sequence of partitions of
\([0,T]\) with mesh tending to zero. For a fixed path \(P\), consider the
pathwise choice
\[
    {\gamma}_{t_i}
    =
    M \operatorname{sign}(P_{t_{i+1}}-P_{t_i}).
\]
Then
\[
\begin{aligned}
X^{x,{\gamma}}_T
&=
\sigma
\sum_{t_i \in \mathcal P_n}
{\gamma}_{t_i}
\bigl(P_{t_{i+1}}-P_{t_i}\bigr) 
=
\sigma M
\sum_{t_i \in \mathcal P_n}
\bigl|P_{t_{i+1}}-P_{t_i}\bigr|.
\end{aligned}
\]
Hence, whenever the path \(P\) has infinite variation, the right-hand side
diverges as \(|\mathcal P_n| \to 0\). In particular, in the absence of
additional regularity constraints or penalisation, the limiting continuous-time
problem cannot be expected to give a finite value.

This example identifies the obstruction addressed in the rest of the paper.
One needs a class of controls which is rich enough to contain genuinely
rough, path-dependent strategies, but sufficiently regular to prevent the
wealth process and the associated value function from degenerating.

\section{Weakly geometric rough paths and controlled paths}\label{section2}

The preceding discussion shows that the control problem cannot be formulated
directly at the level of arbitrary bounded controls. However, the admissibility
of the controls is only one of the issues. Even after imposing suitable
regularity or penalisation on the control, the state dynamics themselves may be
driven by a signal which is too irregular for classical Riemann--Stieltjes or
Young integration. The next step is therefore to specify the pathwise analytic
framework in which these dynamics are defined.

We briefly recall the elements of rough path theory needed in the sequel. In
their simplest form, the dynamics underlying the control problem are described
by equations of the form
\[
     dX_t = \lambda(X_t)\,d\zeta_t,\quad X_0 = x_0,
\]
where the driving signal \(\zeta\) is allowed to be too irregular to have
bounded variation. In this setting, \(\zeta\) must be enhanced with additional
algebraic information satisfying suitable analytic and algebraic conditions,
such as Chen's relation and the shuffle identity. These conditions ensure that
the enhanced signal behaves consistently with the formal calculus of smooth
paths, while still allowing for genuinely rough driving signals. They are the
basis for the pathwise interpretation of the state equation, and hence for the
formulation of the control problem considered below.

We recall these conditions first, before introducing rough differential
equations. For a more detailed treatment, we refer the reader to
\cite{cass2022combinatorial,hairer2015geometric}.
\begin{definition}
Let $n_1,\ldots,n_m\in\mathbb N$. The set of
$(n_1,\ldots,n_m)$-shuffles, denoted by
$\operatorname{Sh}(n_1,\ldots,n_m)$, is the subset of the permutation group
$\mathfrak S_{n_1+\cdots+n_m}$ consisting of all permutations $\sigma$ such
that, for each $i=1,\ldots,m$,
\[
    \sigma(n_1+\cdots+n_{i-1}+1)
    <
    \sigma(n_1+\cdots+n_{i-1}+2)
    <
    \cdots
    <
    \sigma(n_1+\cdots+n_i),
\]
with the convention that $n_1+\cdots+n_0=0$. The set of ordered shuffles,
denoted by $\overline{\operatorname{Sh}}(n_1,\ldots,n_m)$, is the subset of
$\operatorname{Sh}(n_1,\ldots,n_m)$ consisting of those permutations satisfying
\[
    \sigma(n_1)
    \leq
    \sigma(n_1+n_2)
    \leq
    \cdots
    \leq
    \sigma(n_1+\cdots+n_m).
\]
\end{definition}

Let $\beta=(\beta_1,\ldots,\beta_n)$ and
$\delta=(\delta_1,\ldots,\delta_m)$ be two multi-indices, and denote by
\[
    \epsilon := (\beta_1,\ldots,\beta_n,\delta_1,\ldots,\delta_m)
\]
their concatenation. We define the shuffle of $\beta$ and $\delta$ by
\[
    \operatorname{Sh}(\beta,\delta)
    :=
    \left\{
    (\epsilon_{\sigma(1)},\ldots,\epsilon_{\sigma(n+m)})
    \;:\;
    \sigma\in\operatorname{Sh}(n,m)
    \right\}.
\]
Similarly, we define the ordered shuffle of $\beta$ and $\delta$ by
\[
    \overline{\operatorname{Sh}}(\beta,\delta)
    :=
    \left\{
    (\epsilon_{\sigma(1)},\ldots,\epsilon_{\sigma(n+m)})
    \;:\;
    \sigma\in\overline{\operatorname{Sh}}(n,m)
    \right\}.
\]

Let $\beta$ be a multi-index of length $|\beta|$ over an alphabet $1, \dots, d$. We define the inverse ordered
shuffle $\overline{\operatorname{Sh}}^{-1}(\beta)$ as the set of all finite
block decompositions
\[
    (\beta^1,\ldots,\beta^m)
\]
such that
\[
    \beta = (\beta^1,\ldots,\beta^m),
\]
and such that the block decomposition is compatible with an ordered shuffle;
that is, there exists
\[
    \sigma\in
    \overline{\operatorname{Sh}}(|\beta^1|,\ldots,|\beta^m|)
\]
corresponding to this decomposition. We also define the restricted inverse
ordered shuffle by
\[
    \overline{\operatorname{Sh}}^{-1}_1(\beta)
    :=
    \left\{
    (\beta^1,\ldots,\beta^m)\in
    \overline{\operatorname{Sh}}^{-1}(\beta)
    \;:\;
    |\beta^i|\geq 1
    \text{ for every } i=1,\ldots,m
    \right\}.
\]
Throughout the paper, sums over words are always understood to be taken over
words satisfying the condition displayed below the summation
symbol. Thus, for instance,
\[
    \sum_{|\beta|=n}
    \quad\text{means}\quad
    \sum_{\substack{\beta\in \{1, \dots, d\} \\ |\beta|=n}} .
\]

We also fix real finite-dimensional Banach spaces
$U$, $V$ and $E$. Throughout this work, $V^{\otimes k }$ and $T((V))$ denote respectively the $k$-fold tensor product of $V$ and the tensor algebra over $V$ with unit $\mathbf{1}$. We use the convention that $V^{\otimes 0} = \mathbb{R}$. For each $N \geq 0$, we denote by $T^N(V)$ its truncation at level $N$.

Let \(\{e_1,\ldots,e_d\}\) denote the basis of \(V\) dual to
\(\{e^1,\ldots,e^d\}\). For a word
\(\beta=(\beta_1,\ldots,\beta_n)\), we write
\[
    e_\beta := e_{\beta_1}\otimes\cdots\otimes e_{\beta_n}
    \in V^{\otimes n},
\]
and for the empty word we set
\[
    e_\emptyset := 1 \in V^{\otimes 0}=\mathbb R.
\]
Let \(\mathcal{L}(V,U)\) denote the space of bounded linear maps from \(V\) to \(U\) equipped with the operator norm. Then, if \(Y\in \mathcal L(T^N(V),U)\), we define its \(\beta\)-component by
\[
    Y_\beta := Y(e_\beta),
    \qquad |\beta|\leq N.
\]
Thus \(Y_\beta\in U\). If, in addition, \(\{u^1,\ldots,u^e\}\) is a fixed
basis of \(U^*\), we define its scalar coordinates by
\[
    Y^h_\beta := u^h(Y_\beta)
    =
    u^h\bigl(Y(e_\beta)\bigr),
    \qquad h=1,\ldots,e .
\]

Let $T > 0$ and for every $0 \leq s \leq t \leq T$, we denote by $\Delta_{[s,t]}$ the simplex
\[
\Delta_{[s,t]} := \{ (r,u) \in [s,t]^2 : r \leq u \}.
\]
\begin{definition}
Let $p\geq 1$. We denote by
$\mathcal C^p(\Delta_{[0,T]},V)$ the space of continuous
two-parameter maps $\gamma:\Delta_{[0,T]}\to V$ with finite $p$-variation,
that is, such that
\[
    \|\gamma\|_{p\text{-var};[0,T]}
    :=
    \left(
    \sup_{\mathcal P}
    \sum_{[s,t]\in\mathcal P}
    \|\gamma_{st}\|_V^p
    \right)^{1/p}
    <\infty,
\]
where the supremum is taken over all finite partitions $\mathcal P$ of $[0,T]$.
\end{definition}

The same notation extends naturally to ordinary paths. If
$\gamma:[0,T]\to V$ is continuous, we write
$\gamma\in \mathcal C^p([0,T],V)$ whenever its increment map
\(
    (s,t)\mapsto \gamma_t-\gamma_s
\)
belongs to $\mathcal C^p(\Delta_{[0,T]},V)$. From now on, we use
the notation
\(
    \gamma_{st}:=\gamma_t-\gamma_s
\)
for the increments of a path. More generally, for any subinterval
$[s,t]\subseteq [0,T]$, we define
\[
    \|\gamma\|_{p\text{-var};[s,t]}
    :=
    \left(
    \sup_{\mathcal P}
    \sum_{[r,u]\in\mathcal P}
    \|\gamma_{ru}\|_V^p
    \right)^{1/p},
\]
where the supremum is taken over all finite partitions $\mathcal P$ of $[s,t]$. When the meaning is clear from the context, we abbreviate $\|\gamma\|_{p;[s,t]}$ in place of $\|\gamma\|_{p\text{-var};[s,t]}$.
\begin{definition}
Let  $p\geq 1$. We denote by
$\mathscr C^p([0,T],V)$ the space of $p$-weakly geometric rough
paths over $V$. An element of this space is a map
\[
    \boldsymbol{\zeta}:\Delta_T\to T^{\lfloor p\rfloor}(V),
\]
such that, for all $0\leq s\leq u\leq t\leq T$, the following conditions hold:
\begin{enumerate}
    \item 
\(
    \boldsymbol{\zeta}^{\emptyset}_{st}=1,
    \qquad
    \boldsymbol{\zeta}_{tt}=\mathbf 1.
\)
    \item For every multi-index $\beta$ with
    $1\leq |\beta|\leq \lfloor p\rfloor$,
    \[
        \|\boldsymbol{\zeta}^{\beta}\|_{\frac{p}{|\beta|}\text{-var};[0,T]}<\infty.
    \]

    \item Chen's relation holds:
    \[
        \boldsymbol{\zeta}^{\beta}_{st}
        =
        \sum_{(\epsilon,\delta)=\beta}
        \boldsymbol{\zeta}^{\epsilon}_{su}
        \boldsymbol{\zeta}^{\delta}_{ut}.
    \]

    \item The shuffle identity holds:
    \[
        \boldsymbol{\zeta}^{\epsilon}_{st}
        \boldsymbol{\zeta}^{\delta}_{st}
        =
        \sum_{\beta\in \operatorname{Sh}(\epsilon,\delta)}
        \boldsymbol{\zeta}^{\beta}_{st}.
    \]
\end{enumerate}
\end{definition}

 For a given $p$-weakly geometric path $\boldsymbol{\zeta}$ we introduce the $p$-variation seminorm for rough paths: 
\begin{equation*}
    \left\|\boldsymbol{\zeta}\right\|_{p; [s,t]} := \sum_{|\beta| = 1}^{\floor{p}} \left\|\boldsymbol{\zeta}^{\beta}\right\|_{\frac{p}{|\beta|}; [s, t]},  
\end{equation*}
to which we associate a norm defined via the map $\boldsymbol{\zeta} \rightarrow |\boldsymbol{\zeta}_s| + \left\|\boldsymbol{\zeta}\right\|_{p; [s,t]}$.\newline
We denote by 
\[
C^k(V, U) \quad \text{the space of $k$-times Fr\'echet differentiable functions } f: V \to U,
\]
and by 
\[
C^k_b(V, U) \subset C^k(V, U)
\]
the subspace of functions with bounded derivatives up to order $k$, with the usual extension to \(k = \infty\).\newline
Notice that any path $\zeta \in C^1([0,T], V)$ can be made into a $p$-weakly geometric rough path via the map 
\[ \Delta_{[0, T]} \ni (s,t) \rightarrow \left(1,  \zeta_{st}, \int_{s< t_1 <t_2 < t}  d\zeta_{t_1} \otimes d\zeta_{t_2} ,..., \int_{s< t_1 < ... < t_{\floor{p}} < t}  d\zeta_{t_1} \otimes ... \otimes d\zeta_{t_{\floor{p}}} \right).
\]
The image of this map is called \enquote{canonical lift} of the path $\zeta$ to a $p$-weakly geometric rough path.

\begin{definition}
Let \(p \geq 1\). The space of
\(p\)-geometric rough paths over \(V\) is defined as the closure, with respect
to the \(p\)-variation rough path norm, of the set of canonical
lifts of smooth paths \(\zeta \in C^{\infty}([0,T],V)\) to $p$-geometric rough paths.
\end{definition}

We next introduce controlled rough paths, which provide the natural class of
integrands for integration against a fixed rough path. This notion encodes the
way an integrand depends on the increments of the driving signal and is the
basic framework in which rough integrals are defined.
\begin{definition}
For a given path \(\boldsymbol{\zeta} \in \mathscr{C}^{p}([0,T], V)\), the
class of \(\boldsymbol{\zeta}\)-controlled paths
\(\mathcal{D}_{\boldsymbol{\zeta}}(U)\) is defined as the set of paths
\[
    \overline{X} \in
    \mathcal{C}^p([0,T], \mathcal{L}(T^{\lfloor p \rfloor - 1}(V), U))
\]
such that the following holds. For every \(h=1,\ldots, \mathrm{dim}(U)\), every word
\(\beta\) with \(0 \leq |\beta| \leq \lfloor p \rfloor-2\), and every
\(0\leq s\leq t\leq T\),
\[
    \overline{X}^h_{\beta; t}
    =
    \sum_{|\epsilon|=0}^{\lfloor p \rfloor - 1 - |\beta|}
    \overline{X}^h_{(\epsilon, \beta); s}
    \boldsymbol{\zeta}^{\epsilon}_{st}
    +
    R^{\beta,h}_{st}.
\]
Here, for notational convenience, the component labels of the remainder are
written as superscripts, \(R^{\beta,h}_{st}\), so that the subscript is reserved
for the time increment \((s,t)\).
The remainder satisfies
\[
    R^{\beta}
    \in
    \mathcal{C}^{\frac{p}{\lfloor p\rfloor-|\beta|}}
    \bigl(\Delta_{[0,T]},\mathcal{L}(V^{\otimes |\beta|},U)\bigr).
\]
\end{definition}
Observe that the definition of $\overline X$ ensures that the
regularity condition on the remainder is automatically satisfied when
$|\beta|=\lfloor p\rfloor-1$. In this case, we may therefore set
\[
    R^\beta_{s,t}:=\overline X_{\beta;s,t}, \quad 0 \leq s \leq t \leq T.
\]
When there is possible ambiguity about the underlying controlled rough path, we
write the corresponding remainder as $R^X$.

To simplify notation, whenever
\[
    \overline X
    \in
    \mathcal{D}_{\boldsymbol{\zeta}}\bigl(\mathcal L(T^{\lfloor p\rfloor-1}(V),
    \mathcal L(V,U))\bigr),
\]
we use the following convention.

For any multi-index
$\beta=(\beta_1,\ldots,\beta_n)$ with
$1\leq |\beta|\leq \lfloor p\rfloor$, we write for every $t \in [0, T]$
\[
    \overline X^h_{\beta,t}
    :=
    \overline X^{(\beta^\cdot,h)}_{\beta^-,t},
\]
where
\[
    \beta^- := (\beta_1,\ldots,\beta_{n-1}),
    \qquad
    \beta^\cdot := \beta_n .
\]
For maps with values in \(\mathcal L(V,U)\), we shall use the canonical
levelwise identification
\[
    \mathcal L\!\left(V^{\otimes(k-1)},\mathcal L(V,U)\right)
    \cong
    \mathcal L\!\left(V^{\otimes k},U\right),
\]
given by
\[
    f \mapsto
    \left[
        (v_1 \otimes \cdots \otimes v_k)
        \mapsto
        f(v_1 \otimes \cdots \otimes v_{k-1})(v_k)
    \right].
\]
For a controlled rough path $\overline X$, we define its trace to be the process
\(
    X_t := \overline X_{\emptyset ,t}
\). The higher levels of $\overline X$ are usually referred to as Gubinelli derivatives of $\overline X$.

Finally, the space $\mathcal D_{\boldsymbol\zeta}(U)$ can be endowed with a
Banach-space structure by introducing the norm
\begin{equation}\label{controlled_path_norm}
    \|\overline X\|_{\mathcal{D}_{\boldsymbol{\zeta}};[s,t]}
    :=
    |\overline X_s|
    +
    \sum_{|\beta|=0}^{\lfloor p\rfloor-1}
    \|R^\beta\|_{\frac{p}{\lfloor p\rfloor-|\beta|};[s,t]}, \quad 0 \leq s \leq t \leq T.
\end{equation}
Additionally, if we consider the rough path $\boldsymbol{\eta} \in \mathscr{C}^p([0, T], V)$ and $\overline Y \in \mathcal{D}_{\boldsymbol{\eta}}(U)$, we can define the distance between this controlled path and the path $\overline{X}$ defined above as 
\[
d(\overline X, \overline Y)_{[s,t]} 
    :=
    |\overline X_s - \overline Y_s|
    +
    \sum_{|\beta|=0}^{\lfloor p\rfloor-1}
    \|R^{X, \beta} - R^{Y, \beta}\|_{\frac{p}{\lfloor p\rfloor-|\beta|};[s,t]}, \quad 0 \leq s \leq t \leq T. 
\]
Having introduced weakly geometric rough paths as driving signals and controlled rough paths as admissible integrands, we can now define the rough integral between these two objects.
\begin{proposition}\label{definition_rough_integral}
Let \(\boldsymbol{\zeta}\in \mathscr C^p([0,T],V)\) and let
\(
    \overline{X} \in
    \mathcal{D}_{\boldsymbol{\zeta}}(\mathcal{L}(V,U)).
\)
Then, for every \(0 \leq s < t \leq T\), the limit
\[
     \lim_{|\mathcal{P}| \rightarrow 0}
     \sum_{[r,u] \in \mathcal{P}}
     \sum_{|\beta|=1}^{\floor{p}}
     \overline{X}_{\beta,r}\boldsymbol{\zeta}^{\beta}_{ru}
\]
exists and is independent of the choice of the sequence of finite partitions
\(\mathcal P\) of \([s,t]\). We define the rough integral of
\(\overline X\) against \(\boldsymbol{\zeta}\) by
\[
     \int_s^t \overline{X}_r\, d\boldsymbol{\zeta}_r
     :=
     \lim_{|\mathcal{P}| \rightarrow 0}
     \sum_{[r,u] \in \mathcal{P}}
     \sum_{|\beta|=1}^{\floor{p}}
     \overline{X}_{\beta,r}\boldsymbol{\zeta}^{\beta}_{ru}.
\]
Moreover, the rough integral satisfies the estimate
\begin{equation}\label{continuity_rough_integral}
    \bigg|
    \int_s^t \overline{X}_r\, d\boldsymbol{\zeta}_r
    -
    \sum_{|\beta|=1}^{\floor{p}}
    \overline{X}_{\beta,s}\boldsymbol{\zeta}^{\beta}_{st}
    \bigg|
    \leq
    C_p
    \sum_{|\beta|=1}^{\floor{p}}
    \|\boldsymbol{\zeta}^\beta\|_{\frac{p}{|\beta|};[s,t]}
    \|R^{\beta}\|_{\frac{p}{\floor{p}-|\beta|+1};[s,t]},
\end{equation}
where \(C_p>0\) depends only on \(p\).
\end{proposition}

The next proposition shows that given a sufficiently regular function $\lambda$ and a controlled rough path $X$, it is possible to define a controlled rough path $\lambda(\overline{X})$ with trace $\lambda(X)$.
\begin{proposition}[Composition of controlled paths and functions]
Let $\boldsymbol{\zeta}$ and $\overline{X}$ be as above. For a function $\lambda \in C_b^{\floor{p} }(U, E)$, it is possible to lift the composition $\lambda(X)$ to a controlled rough path $\lambda(\overline{X})$ by defining    
\begin{alignat*}{2}
\lambda(\overline X)_\beta
&:= \lambda(X)^\beta,
&\qquad& |\beta|=0,
\\[0.5em]
\lambda(\overline X)_\beta
&:= \sum_{(\beta^1,\ldots,\beta^m) \in \overline{Sh}_1^{-1}(\beta^-)}
\frac{\partial^k \lambda(X)^\beta}{\partial x^k}\,
\prod_{j=1}^m \overline X^{k_j}_{\beta^j},
&\qquad& 1\le |\beta|\le \lfloor p\rfloor-1 .
\end{alignat*}
\end{proposition}

With these results in place it is finally possible to give meaning to the solution to differential equations driven by a rough path. 
\begin{definition}[Solution to RDE] \label{definition_solution_RDE}
    Let $\boldsymbol{\zeta}$ be as above and consider the equation 
    \begin{align*}
        & X_t - x_0 = \int_{0}^{t} \lambda(\overline{X})_r d\boldsymbol{\zeta}_r, \nonumber   
    \end{align*}
    where $\lambda \in C_b^{\floor{p}}(U, \mathcal{L}(V, U))$ and $ t \in [0, T]$.\newline
    We say that $X$ solves the previous equation if there exists a controlled rough path $\overline{X}  \in \mathcal{D}_{\boldsymbol{\zeta}}(U)$ such that 
    \begin{align*}
        & X_t - x_0 = \int_{0}^{t} \lambda(\overline{X})_r d\boldsymbol{\zeta}_r,  \\
        & \overline{X}_{\beta, t} = \lambda(\overline{X})_{\beta, t}, \quad 1\leq |\beta| \leq \floor{p}-1   .      
    \end{align*}

\end{definition}

\section{Optimal control of rough differential equations}\label{section_rough}
This section provides the setup for the control problem studied in the paper. We first introduce dynamics of the state process as a rough differential equation driven by a geometric rough path. We then show that, without a sufficiently regularising cost functional, the control problem is ill posed. We subsequently derive estimates for the associated state process in terms of the control that will be instrumental in obtaining a non-degeneracy condition for the value functional later in the paper.
\subsection{Setup}

Consider a geometric rough path $\boldsymbol{\zeta} \in \mathscr{C}^{p}([0,T], \mathbb{R}^d)$, with $p \in [2, \infty) \setminus \mathbb N$, for convenience we will assume that for some $L>0$, $\|\boldsymbol{\zeta}\|_{p; [0, T]} < L$.

The state process $\overline{X} \in \mathcal{D}_{\boldsymbol{\zeta}}(\mathbb{R}^e)$, controlled (in the sense of optimal control) by $\gamma \in \mathcal{C}^{p/ \lfloor p \rfloor}([0,T], \mathbb{R}^\ell)$, satisfies the rough differential equation
\begin{align}
dX_t &= b(X_t, \gamma_t)\, dt + \lambda(\overline{X}_t, \gamma_t)\, d\boldsymbol{\zeta}_t, \qquad t \in [0, T], \label{control_process}\\
X_0 &= x_0. \nonumber
\end{align}
Before stating the well-posedness result for equation \eqref{control_process}, we introduce the function spaces appearing in the assumptions. 
The space $\mathrm{Lip}_b(U \times V, E)$ denotes the set of bounded Lipschitz functions. 
\begin{proposition}\label{existence_uniqueness_stability_control_process}
Let
\(
    b \in \mathrm{Lip}_b(\mathbb{R}^e\times\mathbb{R}^\ell,\mathbb{R}^e),
    \) and \(
    \lambda \in
    C_b^{\lfloor p\rfloor+1}
    \bigl(
        \mathbb{R}^e\times\mathbb{R}^\ell,
        \mathcal{L}(\mathbb{R}^d,\mathbb{R}^e)
    \bigr).
\)
Let \(L>0\) be fixed, and let
\(\boldsymbol{\zeta},\boldsymbol{\eta}
\in \mathscr{C}^{p}([0,T],\mathbb{R}^d)\) be such that
\[
    \|\boldsymbol{\zeta}\|_{p;[0,T]}
    ,\,
    \|\boldsymbol{\eta}\|_{p;[0,T]}
    \leq L .
\]
For any \(x_0,y_0\in\mathbb{R}^e\) and any
\(\gamma,\nu\in \mathcal{C}^{p/\lfloor p\rfloor}([0,T],\mathbb{R}^\ell)\)
satisfying
\[
    \|\gamma\|_{\frac{p}{\lfloor p\rfloor};[0,T]}
    ,\,
    \|\nu\|_{\frac{p}{\lfloor p\rfloor};[0,T]}
    \leq M ,
\]
there exists a unique solution
\(\overline X\in\mathcal{D}_{\boldsymbol{\zeta}}(\mathbb{R}^e)\) to the RDE
\[
    X_t
    =
    x_0
    +
    \int_0^t b(X_r,\gamma_r)\,dr
    +
    \int_0^t \lambda(\overline X,\gamma)_r\,d\boldsymbol{\zeta}_r,
    \qquad t\in[0,T],
\]
where, for \(1\leq |\beta|\leq \lfloor p\rfloor-1\),
\[
    \overline X_\beta
    =
    \sum_{(\beta^1,\dots,\beta^m)\in
    \overline{Sh}^{-1}_1(\beta^-)}
    \frac{\partial^k \lambda(\overline X,\gamma)^{\beta^.}}{\partial x^k}
    \prod_{j=1}^m \overline X^{k_j}_{\beta^j}.
\]
This solution also satisfies the bound
\[
\|\overline X\|_{\mathcal{D}_{\boldsymbol{\zeta}}; [0, T]} \leq  C_{p,L,\lambda,b,M,x_0} \left( |x_0| + |\gamma_0|
        +
        \|\gamma\|_{\frac{p}{\lfloor p\rfloor};[0,T]}
        +
        \|\boldsymbol{\zeta}\|_{p;[0,T]}
    \right).
\]
Moreover, if
\(\overline Y\in\mathcal{D}_{\boldsymbol{\eta}}(\mathbb{R}^e)\) denotes the
corresponding solution driven by \(\boldsymbol{\eta}\), with initial condition
\(y_0\) and control \(\nu\), then
\[
    d(\overline X, \overline Y)_{[0,T]}
    \leq
    C_{p,L,\lambda,b,M,x_0,y_0}
    \left(
        |x_0-y_0|
        +
        |\gamma_0-\nu_0|
        +
        \|\gamma-\nu\|_{\frac{p}{\lfloor p\rfloor};[0,T]}
        +
        \|\boldsymbol{\zeta}-\boldsymbol{\eta}\|_{p;[0,T]}
    \right).
\]
\end{proposition}

Although this result is classical, we provide a proof in order to carefully highlight the role played by the control~$\gamma$ in the estimates and in the construction of the solution.
\begin{remark}
We regard the path $\gamma$ as an element of $\mathcal{D}_{\boldsymbol{\zeta}}(\mathbb{R}^\ell)$: the regularity assumption on $\gamma$ ensures that it is controlled by $\boldsymbol{\zeta}$ with Gubinelli derivatives that can be chosen to be $0$. An alternative route to well-posedness would be to lift the pair $(\boldsymbol{\zeta}, \gamma)$ and work with \enquote{weighted} rough paths (see \cite{bel23}). This approach, however, would lead to more cumbersome notation without offering any real advantage for the estimates developed later.
\end{remark}

We define the cost functional associated with a control $\gamma$ by 
\[
J(t,x,\gamma) \;=\; \int_t^T f\bigl(X^{t,x,\gamma}_r, \gamma_r\bigr)\,dr \;+\; g\!\left(X^{t,x,\gamma}_T\right).
\]
The corresponding value function is then given by 
\begin{equation}\label{value_function_first_definition}
v(t,x)
=
\inf_{\gamma \in \mathcal{C}^{p/\lfloor p \rfloor}([t,T],\mathbb{R}^\ell)}
J(t,x,\gamma).
\end{equation}
Here we assume that $f \in \mathrm{Lip}_b(\mathbb{R}^e \times \mathbb{R}^\ell, \mathbb{R})$, and that $g: \mathbb{R}^e \to  \mathbb{R}$ satisfies a linear growth condition, namely: there exists $C_g > 0$ such that 
\[
|g(x)| \leq C_g \,(1 + |x|), \quad \forall x \in \mathbb{R}^e.
\]
\begin{remark}
Notice that the restriction to controls in the class 
$\mathcal{C}^{p / \lfloor p \rfloor}$ is essential. Keeping in mind example  \ref{example_subsection}, if one 
were to allow $\gamma \in \mathcal{C}^q$ with 
\[
\frac{1}{q} + \frac{1}{p} \leq 1,
\]  
then in the general case it is possible to construct a sequence of smooth 
paths $\{\gamma^n\}_{n \in \mathbb{N}}$ converging to $\gamma$ in the 
$q$-variation norm, but such that 
\[
\int \gamma^n\, d\zeta \;\to\; \infty 
\quad \text{as } n \to \infty.
\]  
This shows that extending admissible controls beyond the conjugate exponent 
$p/(p-1)$ is not feasible in general.  

On the other hand, we do not consider controls with regularity 
$q \in \bigl(p/\lfloor p \rfloor, p/(p-1)\bigr)$. The reason is that it is 
more natural to restrict attention to controls that are 
immediately controlled by $\boldsymbol{\zeta}$ in the sense of rough 
path theory, i.e.\ they admit a Gubinelli derivative that can be chosen to be 
zero. This choice yields a simpler and more transparent dynamics for the 
controlled system.
\end{remark}

The previous remark also shows that our restriction is not sufficient to 
ensure the well-posedness of the control problem \eqref{value_function_first_definition}. Indeed, by applying the 
same approximation procedure one can in general construct a sequence of 
smooth controls $\{\gamma^n\} \subseteq C^\infty$ such that 
\[
v(t,x) \;\leq\; \inf_{\gamma \in \{\gamma^n\}_{n \in \mathbb{N}}} 
J(t,x,\gamma) \;=\; -\infty.
\]
This is the motivating reason to investigate the dependence of the solution to \eqref{control_process} with respect to $\gamma$ in order to find suitable conditions to impose on this class of controls and our cost functions to make the problem non-degenerate.

\subsection{Dependence on the control and admissibility conditions}

We begin this section by studying the remainder terms associated with the
solution of the controlled state equation \eqref{control_process}. In
particular, we estimate their dependence on the
$p/\floor{p}$-variation norm of the control path $\gamma$. The following preliminary result provides a quantitative estimate, which will serve as a foundation for subsequent bounds on the state process and for controlling the behaviour of the solution with respect to variations in $\gamma$. These bounds also offer an intuitive justification for the restrictions imposed later in the paper on the control and cost functional, which are necessary to guarantee the well-posedness of the problem.

\begin{lemma}\label{lemma_remainders_X_final}
Let \(\boldsymbol{\zeta}\) be as above, and let
\(\overline X\in\mathcal D_{\boldsymbol{\zeta}}(\mathbb R^e)\) denote the
controlled rough path lift of the solution to \eqref{control_process}. Then,
for every \(0\leq s<t\leq T\), the following estimates hold:
\begin{equation}\label{estimate_trace_remainder_control_process}
\begin{aligned}
&\|R^{X}\|_{\frac{p}{\floor{p}}; [s,t]}\leq C_{\lambda, b, p} (t-s + \|\boldsymbol{\zeta}\|_{p; [s,t]})(1+  \|\gamma\|_{\frac{p}{\floor{p}}; [s,t]})\left( 1 + \sum\limits_{j=1}^{\floor{p} - 1 }\|R^{X}\|^j_{\frac{p}{\floor{p}}; [s, t]} \right), \\ 
&\left\|R^{X, \beta}\right\|_{\frac{p}{\floor{p} - |\beta| +1} ; [s,t]} \leq C_{\lambda,p, L}(1+\left\|\gamma\right\|_{\frac{p}{\floor{p}}; [s,t]})\left( 1 + \sum_{j=1}^{\floor{p} - 1 }\left\|R^{X}\right\|^j_{\frac{p}{\floor{p}};[s,t]} \right), \quad |\beta| \geq 1.
\end{aligned}
\end{equation}

\end{lemma}
The proof follows from standard arguments in the theory of RDEs and is
included in the Appendix for the reader's convenience.

The recursive structure of the remainder estimates established in the previous
lemma allows us to control \(R^X\) on sufficiently small subintervals of
\([0,T]\). By combining these local bounds, one can obtain an estimate on the global $p/\floor{p}$-variation norm of $R^X$. This idea is formalized in the following lemma, which is a variant of Lemma 2.3 in \cite{allan2020pathwise}. We postpone the proof of this result to Appendix \ref{controlling_remainders} in order to maintain the flow of the exposition.
\begin{lemma}\label{bound_allen_revisited}
    Let \(\mathcal P=\{0=t_0<t_1<\cdots<t_n=T\}\) be a partition of \([0,T]\).
Then the remainder \(R^X\) of the trace of \(\overline X\) satisfies 
    \begin{equation} 
        \|R^X\|_{\frac{p}{\floor{p}};[0,T ]} \leq C_{p, L, \lambda}  n^{1 -\frac{\floor{p}}{p}} \left( \sum_{i=1}^{n} \left\| R^X\right\|^{\frac{p}{\floor{p}}}_{{\frac{p}{\floor{p}}}; [{t_{i-1},t_{i}}]} +  \sum_{i=1}^{n}  \sum_{|\beta|=1}^{\floor{p}}   \left\| R^{\lambda(X, \gamma), \beta}\right\|^{\frac{p}{\floor{p} - |\beta| +1}}_{{\frac{p}{\floor{p} - |\beta| +1}}; [t_{i-1}, t_i]} +  \left\| \boldsymbol{\zeta}\right\|_{p; [0, T]} \right)^{\frac{\floor{p}}{p}} .\label{partition_p_var_inequality}
    \end{equation}
\end{lemma}

A combination of the results of the two lemmas above
 leads to  explicit global estimates for the remainders of $\overline X$ in terms of the $p / \lfloor p \rfloor$-variation norm of the control $\gamma$, which we formalise below.

\begin{lemma}\label{lemma_3_6}
Let \(\boldsymbol{\zeta}\) be as above, and let
\(\overline X\in\mathcal D_{\boldsymbol{\zeta}}(\mathbb R^e)\) be the
controlled rough path lift of the solution to \eqref{control_process}.
Then, for every \(0\leq s<t\leq T\), the following estimates hold:
\begin{align*}
&\|R^{X}\|_{\frac{p}{\floor{p}}; [s,t]} \leq C_{\lambda, b, p, L, T}  \bigg( 1 + \|\gamma\|^p_{\frac{p}{\floor{p}}; [s,t]}\bigg),  \\
&\|R^{X, \beta}\|_{\frac{p}{\floor{p} - |\beta| +1}; [s,t]} \leq C_{\lambda, b, p, L, T}  \bigg( 1 + \|\gamma\|^{\floor{p}p}_{\frac{p}{\floor{p}}; [s,t]}\bigg),  \quad |\beta| \geq 1 .
\end{align*}
\end{lemma}
\begin{proof}
Denote by $\mathcal{P}$  the partition of $[s,t]$  defined as
\begin{equation*}
s_0 := s \quad \quad s_i := \sup\Big\{ z>s_{i-1} : \|R^{X}\|_{\frac{p}{\floor{p}}; [ s_{i-1}, z]} \leq 1 \Big\} \wedge t. 
\end{equation*}
Using this partition and inequality \eqref{estimate_trace_remainder_control_process} yields
\[
1 \leq C \left(s_i - s_{i-1} + \|\boldsymbol{\zeta}\|_{p; [s_{i-1}, s_i]}\right)\left(1+\|\gamma\|_{\frac{p}{\floor{p}}; [s_{i-1}, s_i]}\right). 
\]
Implying that the number of intervals $n$ in $\mathcal{P}$ satisfies
\begin{align*}
    n &= \sum_{(s_{i-1}, s_i) \in \mathcal{P}} 1 = \sum_{(s_{i-1}, s_i) \in \mathcal{P}} 1^p \leq C^p \sum_{(s_{i-1}, s_i) \in \mathcal{P}} \left(s_i - s_{i-1} + \|\boldsymbol{\zeta}\|_{p; [s_{i-1}, s_i]}\right)^p\left(1+\|\gamma\|_{\frac{p}{\floor{p}}; [s_{i-1}, s_i]}\right)^p \\
    &\lesssim_{\lambda, b, p, L, T} \left(1 + \|\gamma\|^p_{\frac{p}{\floor{p}}; [s, t]} \right).
\end{align*}
This last inequality in conjunction with  \eqref{partition_p_var_inequality}  allows us to obtain the bound on the remainder of the trace of $X$. More precisely recalling that
\begin{align*}
    &\|R^{X}\|_{\frac{p}{\floor{p}}; [s, t]}\\
    &\leq C_{p, L, \lambda}  n^{1 -\frac{\floor{p}}{p}} \left( \sum_{i=1}^{n} \left\| R^X\right\|^{\frac{p}{\floor{p}}}_{{\frac{p}{\floor{p}}}; [{s_{i-1},s_{i}}]} +  \sum_{i=1}^{n}  \sum_{|\beta|=1}^{\floor{p}}   \left\| R^{\lambda, \beta}\right\|^{\frac{p}{\floor{p} - |\beta| +1}}_{{\frac{p}{\floor{p} - |\beta| +1}}; [s_{i-1}, s_{i}]} +  \left\| \boldsymbol{\zeta}\right\|_{\frac{p}{\floor{p}}; [0, T]} \right)^{\frac{\floor{p}}{p}},
    \end{align*}
     and that by using Proposition \ref{preliminary_bound_controlled_process} and Lemma \ref{lemma_remainders_X_final} we can write for every $|\beta| \ge 1$
    \[
    \left\|R^{\lambda(X), \beta}\right\|_{\frac{p}{\floor{p} - |\beta| +1} ; [s,t]} \leq C_{\lambda,p, L}(1+\left\|\gamma\right\|_{\frac{p}{\floor{p}}; [s,t]})\left( 1 + \sum_{j=1}^{\floor{p} - 1 }\left\|R^{X}\right\|^j_{\frac{p}{\floor{p}};[s,t]} \right).
    \]
    These two results and  the definition of the partition $\mathcal{P}$ yield the estimate
    \begin{align*}
    &\|R^{X}\|_{\frac{p}{\floor{p}}; [s, t]}\\
    &\leq C_{p, L, \lambda}  n^{1 -\frac{\floor{p}}{p}} \left( \sum_{i=1}^{n} 1 +  \sum_{i=1}^{n} \left(1 + \|\gamma\|_{\frac{p}{\floor{p}}; [s_{i-1}, s_i]}   ^{p}\right) +  \left\| \boldsymbol{\zeta}\right\|_{\frac{p}{\floor{p}}; [0, T]} \right)^{\frac{\floor{p}}{p}}\\
    &\leq C_{p, L, \lambda}  n^{1 -\frac{\floor{p}}{p}} \left( n  + \|\gamma\|_{\frac{p}{\floor{p}}; [s,t]}   ^{p} + \left\| \boldsymbol{\zeta}\right\|_{\frac{p}{\floor{p}}; [0, T]}  \right)^{\frac{\floor{p}}{p}}\\
    &\leq C_{p, L, \lambda} \left( 1 + \|\gamma\|^{p}_{\frac{p}{\floor{p}}; [s,t]} \right)^{1 -\frac{\floor{p}}{p}} \left( 1 + \|\gamma\|^{p}_{\frac{p}{\floor{p}}; [s,t]} \right)^{\frac{\floor{p}}{p}}\\
    &\leq C_{p, L, \lambda} \left( 1 + \|\gamma\|_{\frac{p}{\floor{p}}; [s,t]}^p \right) 
    ,
\end{align*} 
where, in the penultimate line we used the estimate of $n$ in terms of $\gamma$ that we found above.\newline
The bound on the remainder of the Gubinelli derivatives of $\overline X$ follows immediately from the inequality we have just recovered and Lemma \ref{lemma_remainders_X_final}.
\end{proof}

From the previous result and Proposition \ref{definition_rough_integral} we notice that, since $g$ has linear growth, there exist $C_{\lambda, b, p, L, T, g}, C_{\lambda, b, p, L, T,g,x} >0$ such that
\begin{equation} \label{bound_rough_integral_psi}
\left|g(X^{t,x,\gamma}_T) \right| \leq C_{\lambda, b, p, L, T, g}\left
(1 + |x| + |X^{t,x,\gamma}_T - x|\right) \leq C_{\lambda, b, p, L, T,g,x}  \bigg( 1 + \|\gamma\|^{\floor{p}p }_{\frac{p}{\floor{p}}; [t,T]}\bigg).
\end{equation}

In light of estimate \eqref{bound_rough_integral_psi}, there are at least two ways to ensure well-posedness of the control problem. The first is to restrict the class of admissible controls to those with uniformly bounded \(p / {\lfloor p \rfloor}\)-variation, say by a constant \(C_\gamma > 0\). However, as already noted in \cite{diehl2017stochastic} and later in \cite{allan2020pathwise}, this restriction is unnatural: it breaks the Dynamic Programming Principle (DPP) and prevents the use of HJB-type pointwise optimizations.  

The second method, which lies at the core of \cite{allan2020pathwise}, is to augment the cost functional \(J\) by an auxiliary penalty term \(\tilde{J} : [0, T] \times \mathcal{C}^{\frac{p}{\lfloor p \rfloor}} \to \mathbb{R}\), chosen so that the growth of the control norm is sufficiently penalized. Concretely, one requires
\begin{equation}\label{auxiliary_definition}
\frac{J(t, x ,\gamma) + \tilde{J}(t, \gamma)}{\|\gamma\|^{\lfloor p \rfloor p}_{\frac{p}{\lfloor p \rfloor}; [t, T]}} \to \infty 
\quad \text{as } 
\|\gamma\|^{\lfloor p \rfloor p}_{\frac{p}{\lfloor p \rfloor}; [t, T]} \to \infty.
\end{equation}

To preserve the additivity structure required by the DPP, \cite{allan2020pathwise} restricts the controls further to the Sobolev space \(W^{1,q}\) for a suitably chosen \(q\). This allows one to introduce a “pseudo-control” \(u: [0, T] \to \mathbb{R}^\ell\) and reformulate the dynamics as
\begin{align*}
dX_t &= b(X_t, \gamma_t)\, dt + \lambda(\overline X, \gamma)_t\, d\boldsymbol{\zeta}_t, \quad X_0 = x_0, \\
d\gamma_t &= u_t \, dt, \quad \gamma_0 = a.
\end{align*}

It is then immediate that, for an appropriately chosen \(q\), the additive functional 
\[
\tilde{J}(t, u) = c \int_t^T |u_r|^q \, dr, \quad c > 0,
\]
satisfies condition \eqref{auxiliary_definition}. This ensures that the modified control problem
\[
v(t,x,a) = \inf_{u \in L^q} \big( J(t,x,\gamma) + \tilde{J}(t,u) \big)
\]
is well posed.\newline

The question we aim to answer in the remainder of this paper is whether we can consider a larger class of controls than \(W^{1,q}\). As we will see, the answer is positive. Assuming, without loss of generality, that \(p\) is not an integer
(which can always be achieved by viewing the rough path as having slightly
lower regularity), we choose the class  
\[
AC^{\alpha} := \Big\{ \gamma \in C^{\text{Höl}-\alpha}([0, T], \mathbb{R}^\ell) : \gamma_t = \gamma_0 + I_{0^+}^\alpha u(t),\ u \in L^{\infty}([0, T]) \Big\}, 
\]  
where \(\alpha \in [\floor{p} / p , 1)\),  \(C^{\text{Höl}-\alpha}\) denotes the class of $\alpha$-H\"older continuous functions (see Definition \ref{definition_alpha_holder}) and \(I^\alpha_{0^+} u\) denotes the left-sided Riemann--Liouville integral of order \(\alpha\) with base point \(0\) defined in appendix \ref{fundamentals_fractional_integrations}.\newline

With this choice, the left-sided Riemann--Liouville derivative of order \(\alpha\) defines the pseudo-control $u$ through $\gamma$ via the identity  
\[
D_{0^+}^\alpha (\gamma - \gamma_0)(t) = u_t \quad \text{for all } t \in [0,T].
\]  
This is a genuine extension of the framework in
\cite{allan2020pathwise}, since in their framework \(q\) may be chosen
arbitrarily large subject to \(q\geq \floor{p}(p+1)\). In particular,
we may choose \(\alpha<1-\frac{1}{q}\), so that the following chain
of inclusions holds:
\[
W^{1,q}([0,T],\mathbb{R}^\ell)
\subset
AC^\alpha([0,T],\mathbb{R}^\ell)
\subset
C^{\mathrm{H\ddot{o}l}-\alpha}([0,T],\mathbb{R}^\ell)
\subset
\mathcal{C}^{p/\floor{p}}([0,T],\mathbb{R}^\ell).
\]
Indeed, if \(\gamma\in W^{1,q}([0,T],\mathbb{R}^\ell)\), then, by the
semigroup property of the Riemann--Liouville integral,
\[
\gamma-\gamma_0
=
I_{0^+}^{1}\dot{\gamma}
=
I_{0^+}^{\alpha}
\bigl(I_{0^+}^{1-\alpha}\dot{\gamma}\bigr).
\]
Since \(q>1/(1-\alpha)\), Proposition 2.15 of
\cite{carbotti2021note} implies that
\(I_{0^+}^{1-\alpha}\dot{\gamma}\in L^\infty\), and hence
\(\gamma\in AC^\alpha([0,T],\mathbb{R}^\ell)\).

The last inclusion is classical; see, for instance, Section 5.1.1 of
\cite{caruana2011rough}. The second inclusion is proved in Proposition
\ref{properties_AC}.
For a self-contained and brief introduction to Riemann--Liouville
integrals and derivatives, as well as the space \(AC^\alpha\), we
refer the reader to Appendix
\ref{fundamentals_fractional_integrations}.

In light of the preceding discussion, we reformulate the state dynamics of the
control problem as follows:
\begin{equation}\label{smooth_fractional_control_dynamics}
\begin{aligned}
&dX^{0, x, \gamma, u}_s = b(X^{0, x, \gamma, u}_s, \gamma^{a, u}_s) \, ds + \lambda(\overline X^{0, x, \gamma, u}, \gamma^{a, u})_s \, d\boldsymbol{\zeta}_s,  &&X^{0, x, \gamma, u}_0 = x,  \\
&D^{\alpha}_{0^+}(\gamma^{a, u} - a)(s) = u_s,  && \gamma^{a,u}_0 = a 
\end{aligned}
\end{equation}
and we will construct an appropriate \(\tilde{J}(t, u)\) so that condition \eqref{auxiliary_definition} is satisfied.

\begin{remark}
Notice that the framework in \cite{allan2020pathwise} considers a cost functional of the form
\[
J(t,x, \gamma) = \int_t^T f(X_r^{t, x, \gamma}, \gamma_r)\, dr + \int_t^T \psi(\overline X^{t, x, \gamma}, \gamma)_r\, d\boldsymbol{\zeta}_r + g(X_T^{t, x, \gamma}),
\]
where \(g\) is assumed bounded below, and \(\psi\) is sufficiently regular to ensure the rough integral is well-defined. In this setting, the degeneracy in the value function arises solely from the rough integral term. Importantly, however, whenever \(g\) is allowed to be unbounded below, the regularity of the terminal cost and the rough integral with respect to the control \(\gamma\) are of the same order. Consequently, our framework, which omits the term \(\int \psi\, d\boldsymbol{\zeta}\) in favour of allowing an unbounded \(g\), does not reduce the mathematical complexity of the problem. On the contrary, it accommodates cost functionals that are more representative of those commonly considered in the literature, while preserving the full analytical challenges associated with rough control problems. 
\end{remark}

\section{Recovering the non-degeneracy of the control problem}\label{section_recovering}

Similarly to the previous section, in order to identify a suitable form of the additional cost functional 
$\tilde{J}$ that guarantees the well-posedness of the control problem, we begin 
by deriving bounds on the norm of $\gamma$ in terms of the $L^q$-norm of the 
pseudo-control $u$. These estimates will naturally suggest an appropriate choice 
of $\tilde{J}$.

This is achieved through a simple argument that exploits the representation of a 
path $\gamma \in AC^{\alpha}([0, T], \mathbb{R}^\ell)$ as the Riemann--Liouville integral of the control $u$.

\begin{proposition} \label{cost_function_fractional}
    Let $\gamma \in AC^{\alpha}([0,T], \mathbb{R}^\ell)$. Then for any $1 < \kappa \leq \frac{1}{1 - \alpha  + \frac{\floor{p}}{p}}$ we have
    \begin{align*}
        \|\gamma\|^{\frac{p}{\floor{p}}}_{\frac{p}{\floor{p}}; [t,T]} \leq C_{\alpha, p, T} \Bigg( \Big(\int_t^T |u_r|^\frac{\kappa}{\kappa -1} dr\Big)^\frac{p (\kappa-1)}{\floor{p}\kappa} + \Big(\int_0^t |u_r|^\frac{\kappa}{\kappa -1} dr\Big)^\frac{p (\kappa-1)}{\floor{p}\kappa}  \Bigg), 
    \end{align*}
    where $u := D_{0^+}^\alpha(\gamma - \gamma_0)$.
\end{proposition}

\begin{proof}
From the definition of \(AC^\alpha([0,T],\mathbb R^\ell)\), for any
\(0\leq s<v\leq T\),
\begin{align*}
    |\gamma_v-\gamma_s|
    &=
    |I_{0^+}^\alpha u(v)-I_{0^+}^\alpha u(s)|
    \\
    &=
    \bigg|
        \frac{1}{\Gamma(\alpha)}
        \int_0^v
        \frac{u_r}{(v-r)^{1-\alpha}}\,dr
        -
        \frac{1}{\Gamma(\alpha)}
        \int_0^s
        \frac{u_r}{(s-r)^{1-\alpha}}\,dr
    \bigg|
    \\
    &\leq
    \frac{1}{\Gamma(\alpha)}
    \bigg|
        \int_s^v
        \frac{u_r}{(v-r)^{1-\alpha}}\,dr
    \bigg|
    \\
    &\quad
    +
    \frac{1}{\Gamma(\alpha)}
    \bigg|
        \int_0^s
        u_r
        \bigg(
            \frac{1}{(v-r)^{1-\alpha}}
            -
            \frac{1}{(s-r)^{1-\alpha}}
        \bigg)\,dr
    \bigg|.
\end{align*}
By H\"older's inequality, we have
\begin{align*}
    |\gamma_v-\gamma_s|
    &\leq
    C_\alpha
    \left(
        \int_s^v
        |u_r|^{\frac{\kappa}{\kappa-1}}\,dr
    \right)^{\frac{\kappa-1}{\kappa}}
    \left(
        \int_s^v
        (v-r)^{\kappa(\alpha-1)}\,dr
    \right)^{\frac1\kappa}
    \\
    &\quad
    +
    C_\alpha
    \left(
        \int_0^s
        |u_r|^{\frac{\kappa}{\kappa-1}}\,dr
    \right)^{\frac{\kappa-1}{\kappa}}
    \left(
        \int_0^s
        \bigg(
            \frac{1}{(s-r)^{1-\alpha}}
            -
            \frac{1}{(v-r)^{1-\alpha}}
        \bigg)^\kappa
        \,dr
    \right)^{\frac1\kappa}.
\end{align*}
The restriction on \(\kappa\) implies that
\(\kappa(1-\alpha)<1\), and therefore
\[
    \int_s^v
    (v-r)^{\kappa(\alpha-1)}\,dr
    =
    \frac{
        (v-s)^{1-\kappa(1-\alpha)}
    }{
        1-\kappa(1-\alpha)
    }.
\]
Moreover, using the elementary inequality
\[
    (x-y)^\kappa\leq x^\kappa-y^\kappa,
    \qquad x\geq y\geq0,
\]
we obtain
\begin{align*}
    &\int_0^s
    \bigg(
        \frac{1}{(s-r)^{1-\alpha}}
        -
        \frac{1}{(v-r)^{1-\alpha}}
    \bigg)^\kappa
    \,dr
    \\
    &\qquad\leq
    \int_0^s
    \left(
        \frac{1}{(s-r)^{\kappa(1-\alpha)}}
        -
        \frac{1}{(v-r)^{\kappa(1-\alpha)}}
    \right)\,dr
    \\
    &\qquad=
    \frac{
        s^{1-\kappa(1-\alpha)}
        +(v-s)^{1-\kappa(1-\alpha)}
        -v^{1-\kappa(1-\alpha)}
    }{
        1-\kappa(1-\alpha)
    }
    \\
    &\qquad\leq
    \frac{
        (v-s)^{1-\kappa(1-\alpha)}
    }{
        1-\kappa(1-\alpha)
    }.
\end{align*}
Consequently,
\begin{align*}
    |\gamma_v-\gamma_s|
    &\leq
    C_{\alpha,p}
    \Bigg[
        \left(
            \int_s^v
            |u_r|^{\frac{\kappa}{\kappa-1}}\,dr
        \right)^{\frac{\kappa-1}{\kappa}}
        +
        \left(
            \int_0^s
            |u_r|^{\frac{\kappa}{\kappa-1}}\,dr
        \right)^{\frac{\kappa-1}{\kappa}}
    \Bigg]
    (v-s)^{\alpha-1+\frac1\kappa}.
\end{align*}
Taking the \(p/\floor{p}\)-power on both sides yields
\begin{align*}
    |\gamma_v-\gamma_s|^{\frac{p}{\floor{p}}}
    &\leq
    C_{\alpha,p}
    \Bigg[
        \left(
            \int_s^v
            |u_r|^{\frac{\kappa}{\kappa-1}}\,dr
        \right)^{
            \frac{p(\kappa-1)}{\floor{p}\kappa}
        }
        +
        \left(
            \int_0^s
            |u_r|^{\frac{\kappa}{\kappa-1}}\,dr
        \right)^{
            \frac{p(\kappa-1)}{\floor{p}\kappa}
        }
    \Bigg]
    (v-s)^{
        \frac{p}{\floor{p}}
        \left(\alpha-1+\frac1\kappa\right)
    }.
\end{align*}

Now let
\[
    t=t_0<t_1<\cdots<t_n=T
\]
be an arbitrary partition of \([t,T]\). Applying the preceding
estimate to each interval \([t_i,t_{i+1}]\) gives
\begin{align*}
    \sum_{i=0}^{n-1}
    |\gamma_{t_{i+1}}-\gamma_{t_i}|^{\frac{p}{\floor{p}}}
    &\leq
    C_{\alpha,p}
    \left(
        \int_0^T
        |u_r|^{\frac{\kappa}{\kappa-1}}\,dr
    \right)^{
        \frac{p(\kappa-1)}{\floor{p}\kappa}
    }
    \\
    &\quad\times
    \sum_{i=0}^{n-1}
    (t_{i+1}-t_i)^{
        \frac{p}{\floor{p}}
        \left(\alpha-1+\frac1\kappa\right)
    }.
\end{align*}
By the definition of \(\kappa\),
\[
    \frac{p}{\floor{p}}
    \left(\alpha-1+\frac1\kappa\right)
    \geq1.
\]
It follows that
\begin{align*}
    \sum_{i=0}^{n-1}
    (t_{i+1}-t_i)^{
        \frac{p}{\floor{p}}
        \left(\alpha-1+\frac1\kappa\right)
    }
    &\leq
    (T-t)^{
        \frac{p}{\floor{p}}
        \left(\alpha-1+\frac1\kappa\right)
    }.
\end{align*}
Taking the supremum over all partitions of \([t,T]\), and absorbing
the last factor into the constant, we obtain
\[
    \|\gamma\|_{\frac{p}{\floor{p}};[t,T]}^{\frac{p}{\floor{p}}}
    \leq
    C_{\alpha,p,T}
    \left(
        \int_0^T
        |u_r|^{\frac{\kappa}{\kappa-1}}\,dr
    \right)^{
        \frac{p(\kappa-1)}{\floor{p}\kappa}
    }.
\]
Finally, splitting the integral at \(t\), we have
\begin{align*}
    \left(
        \int_0^T
        |u_r|^{\frac{\kappa}{\kappa-1}}\,dr
    \right)^{
        \frac{p(\kappa-1)}{\floor{p}\kappa}
    }
    &\leq
    C_p
    \Bigg[
        \left(
            \int_t^T
            |u_r|^{\frac{\kappa}{\kappa-1}}\,dr
        \right)^{
            \frac{p(\kappa-1)}{\floor{p}\kappa}
        } +
        \left(
            \int_0^t
            |u_r|^{\frac{\kappa}{\kappa-1}}\,dr
        \right)^{
            \frac{p(\kappa-1)}{\floor{p}\kappa}
        }
    \Bigg],
\end{align*}
which concludes the proof.
\end{proof}

 From the bound we have just recovered and equation \eqref{bound_rough_integral_psi} we obtain
\begin{align}
&\left| g(X^{t,x,\gamma, u}_T) \right| \leq C_{\lambda, b, p, L, T,g, x}  \bigg( 1 + \|\gamma\|^{\floor{p} p }_{\frac{p}{\floor{p}}; [t,T]}\bigg)\notag\\
& \leq C_{\lambda, b, p, L, T,g, x, \alpha}  \bigg( 1 +  \Big(\int_t^T |u_r|^\frac{\kappa}{\kappa -1} dr\Big)^\frac{\floor{p}p (\kappa-1)}{\kappa } + \Big(\int_0^t |u_r|^\frac{\kappa}{\kappa -1} dr\Big)^\frac{\floor{p}p  (\kappa-1)}{\kappa}  \bigg)\notag \\
&\leq  C_{\lambda, b, p, L, T,g, x, \alpha}  \bigg( 1 +  \Big(\int_t^T |u_r|^{ \floor{p}p \vee \frac{\kappa}{\kappa -1}} dr\Big) + \Big(\int_0^t |u_r|^\frac{\kappa}{\kappa -1} dr\Big)^\frac{\floor{p}p (\kappa-1)}{\kappa}  \bigg) \label{bound_rough_integral_psi_u}.
\end{align}
A possible choice of $\tilde{J}(t, u)$ is now 
\begin{equation}\label{definition_J_tilde}
\tilde{J}(t, u) = \int_t^T \tilde{f}(u_r) dr,
\end{equation}
for any continuous function $\tilde{f}(u) \geq f_0 |u|^q - C_0,$ where $f_0, C_0 > 0$ and  $q > \floor{p}p \vee \frac{\kappa}{\kappa-1}$.

\begin{remark}\label{remark_allan_6_4}
Let us compare more precisely our admissibility condition with that of
\cite{allan2020pathwise}. Let \(\boldsymbol{\zeta}\) have finite
\(p_0\)-variation for some \(p_0\in[2,3)\). Allan and Cohen consider
controls \(\gamma\in W^{1,q}\), with \(u=\dot{\gamma}\), and require the
penalty exponent \(q\) in \eqref{definition_J_tilde} to satisfy
\[
    q>2(p_0+1).
\]
In our framework, \(\gamma\in AC^\alpha\), with
\(u=D_{0+}^{\alpha}(\gamma-\gamma_0)\), and the corresponding condition is
\[
    q>
    2\bar p
    \vee
    \frac{\kappa}{\kappa-1},
    \qquad
    1<\kappa\leq
    \frac{1}{1-\alpha+2/\bar p},
\]
where we may use any \(\bar p\in[p_0,3)\), since
\(\boldsymbol{\zeta}\) also has finite \(\bar p\)-variation.

Choose \(\bar p=\max\{p_0,5/2\}\) and take \(\kappa\) at its maximal
admissible value. As \(\alpha\uparrow1\),
\[
    \frac{\kappa}{\kappa-1}
    \to
    \frac{\bar p}{\bar p-2}
    \leq 2\bar p<6.
\]
Consequently, for \(\alpha\) sufficiently close to \(1\), one may choose
\[
    2\bar p<q<6,
\]
whereas \cite{allan2020pathwise} requires
\(q>2(p_0+1)\geq6\).

The improvement comes from the greedy-partition estimate in
Lemma~\ref{lemma_3_6}. It avoids the additional \(n^{1/\bar p}\) loss
arising when the local estimates are concatenated using the partition
bound employed in \cite{allan2020pathwise}. This lowers the
rough-dynamics contribution to the required growth exponent from
\(2(\bar p+1)\) to \(2\bar p\), while the freedom to increase the
roughness exponent to \(\bar p\) prevents the fractional integrability
requirement from obscuring this gain.
\end{remark}

From this point onward, in order to streamline the notation, we
incorporate the pseudo-control penalisation into the running cost.
More precisely, with a slight abuse of notation, we redefine
\(
f(x,y,u):=f(x,y)+\widetilde f(u).
\)
Accordingly, the cost functional \(J\) already includes the
contribution of \(\widetilde J\).

Before writing down explicitly the new cost functional and proceeding to the proof of local boundedness of the value functional, we make the following observation. 
If $\alpha < 1$, it is well known that the  fractional integral is not a local operator. Hence, in order to obtain the value of $\gamma \in AC^{\alpha}([0,T], \mathbb{R}^\ell)$ at a point $t \in [0,T]$, it is necessary to provide the full path of its fractional derivative. Analogously, when extending a path  $\gamma$ defined on $[0, r]$ to a path $\tilde{\gamma}$ defined on $[0, z]$ by using the fractional derivative one must know the values of the fractional derivative of the former path up to the concatenation point. This justifies the choice to introduce the path $\nu^{r, \gamma, z, u}: [0, z] \rightarrow \mathbb{R}^\ell$ to denote the unique path that agrees with $\gamma$ up to $r$ and has fractional derivative $u$ from time $r$ to $z \leq T$. From this characterisation, $\nu^{r, \gamma, z, u}$ satisfies the integral equation
    \begin{equation*}
        \nu^{r, \gamma, z, u}_t = \gamma_0 + \frac{1}{\Gamma(\alpha)}\int_0^r \frac{D^{\alpha}_{0^+}(\gamma - \gamma_0)(s)}{(t-s)^{1-\alpha}} ds + \frac{1}{\Gamma(\alpha)}\int_r^t \frac{u_s}{(t-s)^{1-\alpha}} ds, 
    \end{equation*}
with $u \in L^{\infty}([r,z], \mathbb{R}^\ell)$ and $\Gamma$ is the Gamma function.
An important property of this new functional is the inequality
\begin{equation}
    |\nu_s^{ r, \gamma, t, u} - \nu_s^{ r, \tilde{\gamma}, t, u}| \leq 2\|\gamma - \tilde{\gamma}\|_{\infty; [0, r]} ,\label{Lipschitz_continuity_nu}  
\end{equation}
which holds for any $0 \leq r \leq s \leq t \leq T$ , $u \in  L^\infty([0, T ], \mathbb{R}^\ell)$ and $\gamma, \tilde{\gamma} \in AC^{\alpha}([0, T ], \mathbb{R}^\ell)$ and can be deduced from Lemma 7.2 in \cite{gomoyunov2020dynamic}.

The introduction of the path $\nu^{r,\gamma,z,u}$ allows us to consistently extend trajectories and controls beyond their original time horizon. With this notation in place, we can now formulate the extended cost functional in a way that is compatible with the fractional framework. Specifically, we define
\begin{align}\label{extended_cost_functional}
J(t,x,\gamma,u) = \int_t^T f\big(X^{t,x,\nu^{t, \gamma, T,u},u}_r, \nu^{t,\gamma,T,u}_r, u_r\big)\,dr 
+ g\big(X^{t,x,\nu^{t, \gamma, T,u},u}_T\big),
\end{align}
and the associated value functional
\begin{equation}\label{extended_value_functional}
    v(t,x,\gamma) 
    = \inf_{u \in L^{\infty}([t,T],\mathbb{R}^\ell)} J(t,x,\gamma,u).
\end{equation}
Following the preceding observation, we define the value functional $v$ on
triples $(t,x,\gamma)$, where $\gamma$ represents the relevant control
trajectory. Thus, $v$ depends on the path $\gamma$ itself, rather than only on
the pointwise value $\gamma_t$.\newline

To improve readability, whenever no ambiguity can arise, we write
\(X^{t,x,\nu,u}_r\) in place of
\(X^{t,x,\nu^{t,\gamma,T,u},u}_r\).\newline

In order to analyse regularity properties of $v$, we equip the space $[0,T]\times \mathbb{R}^e \times AC^{\alpha}([0,T],\mathbb{R}^\ell)$ with the product metric induced by the norm
\begin{equation}\label{product_metric_definition}
    (t,x,\gamma) \mapsto \max\Big\{t,\, |x|,\,  |\gamma_0| + \|\gamma\|_{p/\floor p;[0,T]}\Big\}.
\end{equation}

The first key property we establish is the well-posedness of the value
functional. This is precisely the property for which the analysis in
Sections~\ref{section_rough} and~\ref{section_recovering} was developed: the estimate
\eqref{bound_rough_integral_psi_u} quantifies the growth of the controlled
dynamics, while the penalisation term $\tilde J$ defined in
\eqref{definition_J_tilde} is chosen to offset this growth. Their combination
yields the desired local boundedness.

\begin{proposition}\label{boundedness_val_fun}
For every $k\geq 1$, the value functional $v$ defined in
\eqref{extended_value_functional} is bounded on
\(
    [0,T]\times \overline{B(0,k)}\times AC^\alpha_k,
\)
where $\overline{B(0,k)} \subset \mathbb R^e$ denotes the closed ball of radius $k$ and $AC^\alpha_k$ is defined as in Appendix \ref{fundamentals_fractional_integrations}.
Equivalently,
\[
    \sup_{\substack{t\in[0,T],\, |x|\leq k\\
                    \gamma\in AC^\alpha_k}}
    |v(t,x,\gamma)|<\infty.
\]
\end{proposition}

\begin{proof}
The upper bound is verified using standard techniques. In fact, consider the control $u\equiv 0$. By definition of $v$ we have
\begin{align*}
v(t,x,\gamma)
&\leq J(t,x,\gamma,u)\\
&=
\int_t^T
f\big(X^{t,x,\nu,u}_r,\nu^{t,\gamma,T,u}_r,u_r\big)\,dr
+
g\big(X^{t,x,\nu,u}_T\big)\\
&\leq
C_f
+
C_{\lambda,b,p,L,T,g,k,\alpha}
\bigg(
1+
\Big(
\int_0^t
\big|
D_{0^+}^{\alpha}(\gamma^{a,u}-a)(r)
\big|^{\frac{\kappa}{\kappa-1}}
\,dr
\Big)^{\frac{\floor{p}p(\kappa-1)}{\kappa}}
\bigg),
\end{align*}
where in the last line we used the definition and boundedness of $f$ alongside the estimate
\eqref{bound_rough_integral_psi_u} for the term
$g\big(X^{t,x,\nu,u}_T\big)$, and we used $|x|\leq k$ to absorb the dependence on the initial state into the constant.
We can now conclude by recalling that, for every $k\geq 1$, $r\in[0,T]$, and
$\gamma\in AC^\alpha_k$,
\(
\big|
D_{0^+}^{\alpha}(\gamma^{a,u}-a)(r)
\big|
\leq k.
\)
Hence the preceding estimate gives an upper bound which is uniform on
\(
[0,T]\times\overline{B(0,k)}\times AC^\alpha_k.
\)

For the lower bound, we again use the same two properties. Given any
\(\tilde u\in L^\infty([t,T];\mathbb{R}^\ell)\), set
\[
u_s
=
D^\alpha_{0^+}(\gamma-\gamma_0)_s\mathbbm{1}_{\{s<t\}}
+
\tilde u_s\mathbbm{1}_{\{s\geq t\}}.
\]
Then, thanks to \eqref{bound_rough_integral_psi_u} and the definition of
$\tilde{J}$ in \eqref{definition_J_tilde},
\begin{align*}
J(t,x,\gamma,u)
&\geq
-C_f
+
\int_t^T f_0|u_r|^q\,dr\\
&\quad
-
C_{\lambda,b,p,L,T,g,k,\alpha}
\bigg(
1
+
\int_t^T
|u_r|^{p\floor{p}\vee\frac{\kappa}{\kappa-1}}
\,dr
+
\Big(
\int_0^t
|u_r|^{\frac{\kappa}{\kappa-1}}
\,dr
\Big)^{\frac{\floor{p}p(\kappa-1)}{\kappa}}
\bigg).
\end{align*}
Rearranging the terms in the last line and re-defining the second constant gives
\begin{align*}
J(t,x,\gamma,u)
&\geq
\int_t^T
\left(
f_0|u_r|^q
-
C_{\lambda,b,p,L,T,g,k,\alpha}
|u_r|^{p\floor{p}\vee\frac{\kappa}{\kappa-1}}
\right)\,dr\\
&\quad
-
C_{\lambda,b,p,L,T,g,k,\alpha,f}
\bigg(
1
+
\Big(
\int_0^t
|u_r|^{\frac{\kappa}{\kappa-1}}
\,dr
\Big)^{\frac{\floor{p}p(\kappa-1)}{\kappa}}
\bigg)\\
&\geq
-C_{q,p,\kappa,k}(T-t)
-
C_{\lambda,b,p,L,T,g,k,\alpha,f}
\bigg(
1
+
\Big(
\int_0^t
|u_r|^{\frac{\kappa}{\kappa-1}}
\,dr
\Big)^{\frac{\floor{p}p(\kappa-1)}{\kappa}}
\bigg),
\end{align*}
where in the last line we used the fact that
\[
f_0>0,
\qquad
q>
p\floor{p}\vee\frac{\kappa}{\kappa-1}.
\]
Since $\gamma\in AC^\alpha_k$, we have
\[
|u_r|
=
\big|
D^\alpha_{0^+}(\gamma-\gamma_0)(r)
\big|
\leq k,
\qquad r\in[0,t].
\]
Consequently, the right-hand side above is bounded below by a constant
depending only on $k$ and the coefficients of the problem. Taking the
infimum over the admissible pseudo-controls therefore yields a lower
bound for $v$ which is uniform on
\(
[0,T]\times\overline{B(0,k)}\times AC^\alpha_k.
\)
Together with the upper bound obtained above, this concludes the proof.
\end{proof}

\section{Properties of the value functional}\label{section_properties}

Building on the boundedness result established in Proposition
\ref{boundedness_val_fun}, we now turn to further properties of the value functional
$v(t,x,\gamma)$: the Dynamic Programming Principle, non-anticipativity, and continuity.

\begin{lemma}[Dynamic Programming Principle (DPP)]\label{dpp_principle}
For any $(t,x,\gamma) \in [0,T]\times\mathbb{R}^e\times
AC^{\alpha}([0,T],\mathbb{R}^\ell)$ and $t\leq s\leq T$, the following identity holds:
\begin{equation*}
v(t,x,\gamma)
=
\inf_{u\in L^{\infty}([t,T],\mathbb{R}^\ell)}
\Bigg(
\int_t^s
f\big(X^{t,x,\nu,u}_r,\nu^{t,\gamma,T,u}_r,u_r\big)\,dr
+
v\big(s,X^{t,x,\nu,u}_s,\nu^{t,\gamma,T,u}\big)
\Bigg).
\end{equation*}
\end{lemma}

\begin{proof}
The proof follows from a standard concatenation argument.\newline
Fix any triple $(t,x,\gamma)$ and pick any control
$u\in L^{\infty}([t,T],\mathbb{R}^\ell)$. From the definition of the value function, it follows that for any $t\leq s\leq T$
\begin{align*}
&v(t,x,\gamma)\\
&\leq
\int_t^T
f\big(X^{t,x,\nu,u}_r,\nu^{t,\gamma,T,u}_r,u_r\big)\,dr
+
g\big(X^{t,x,\nu,u}_T\big)\\
&=
\int_t^s
f\big(X^{t,x,\nu,u}_r,\nu^{t,\gamma,T,u}_r,u_r\big)\,dr
+
\int_s^T
f\big(X^{t,x,\nu,u}_r,\nu^{t,\gamma,T,u}_r,u_r\big)\,dr
+
g\big(X^{t,x,\nu,u}_T\big).
\end{align*}
Fix $\tilde u\in L^\infty([t,s],\mathbb R^\ell)$. Since the preceding
inequality holds for every continuation
$\hat u\in L^\infty((s,T],\mathbb R^\ell)$, taking the infimum over
$\hat u$ yields
\begin{align*}
v(t,x,\gamma)
&\leq
\int_t^s
f\big(X^{t,x,\nu,u}_r,\nu^{t,\gamma,T,u}_r,\tilde{u}_r\big)\,dr
+
v\big(s,X_s^{t,x,\nu,\tilde u},\nu^{t,\gamma,T,\tilde{u}}\big).
\end{align*}
Notice that from the definition of $\nu$, the path
$\nu^{t,\gamma,T,\tilde{u}}$ does not depend on the full trajectory of
$\gamma$, but only on $\gamma|_{[0,t]}$. Similarly, since
\begin{align*}
&v\big(
s,
X_s^{t,x,\nu^{t,\gamma,T,\tilde{u}},\tilde u},
\nu^{t,\gamma,T,\tilde{u}}
\big)\\
&=
\inf_{u\in L^{\infty}([s,T],\mathbb{R}^\ell)}
\Bigg(
\int_s^T
f\big(
X^{s,X_s^{t,x,\nu^{t,\gamma,T,\tilde{u}},\tilde u},
\nu^{t,\gamma,T,\tilde{u}}_s,u}_r,
\nu^{t,\nu^{t,\gamma,T,\tilde{u}},T,u}_r,
u_r
\big)\,dr\\
&\qquad\qquad\qquad
+
g\big(
X^{s,X_s^{t,x,\nu^{t,\gamma,T,\tilde{u}},\tilde u},u}_T
\big)
\Bigg),
\end{align*}
the value function depends on $\tilde{u}$ only through its restriction to the interval $[t,s]$.\newline
From this we can proceed to take the infimum over the paths
$u\in L^{\infty}([t,T],\mathbb{R}^\ell)$ to conclude, after appropriate relabelling of $\tilde u$, that
\[
v(t,x,\gamma)
\leq
\inf_{u\in L^{\infty}([t,T],\mathbb{R}^\ell)}
\Bigg(
\int_t^s
f\big(X^{t,x,\nu,u}_r,\nu^{t,\gamma,T,u}_r,u_r\big)\,dr
+
v\big(s,X^{t,x,\nu,u}_s,\nu^{t,\gamma,T,u}\big)
\Bigg).
\]
For the converse inequality, we fix $\epsilon>0$ and pick a control
$\tilde{u}$ such that
\begin{align*}
v(t,x,\gamma)+\epsilon
&\geq
J(t,x,\gamma,\tilde{u})\\
&\geq
\int_t^s
f\big(
X^{t,x,\nu,\tilde{u}}_r,
\nu^{t,\gamma,T,\tilde{u}}_r,
\tilde{u}_r
\big)\,dr
+
v\big(
s,
X_s^{t,x,\nu,\tilde u},
\nu^{t,\gamma,T,\tilde{u}}
\big)\\
&\geq
\inf_{u\in L^{\infty}([t,T],\mathbb{R}^\ell)}
\Bigg(
\int_t^s
f\big(
X^{t,x,\nu,u}_r,
\nu^{t,\gamma,T,u}_r,
u_r
\big)\,dr
+
v\big(
s,
X^{t,x,\nu,u}_s,
\nu^{t,\gamma,T,u}
\big)
\Bigg),
\end{align*}
where, in the above inequalities, we used the same considerations about the dependence of the value function and the cost functional on the control $\tilde{u}$ and the path $\gamma$ that we outlined before.\newline
Since $\epsilon$ was chosen to be arbitrary, the two inequalities we just derived are enough to conclude.
\end{proof}

In the proof above, we used the fact that the construction of
\(\nu^{t,\gamma,T,u}\) is non-anticipative in \(\gamma\): it depends on
\(\gamma\) only through its restriction \(\gamma|_{[0,t]}\). Consequently, the
cost functional
\[
J(t,x,\gamma,u)
=
\int_t^T f\big(X^{t,x,\nu,u}_r,
\nu^{t,\gamma,T,u}_r,u_r\big)\,dr
+
g\big(X^{t,x,\nu,u}_T\big)
\]
is also non-anticipative in the path variable \(\gamma\). It follows that the
value functional is non-anticipative in \(\gamma\). This property is formalized
in the definition below.
\begin{definition}
A functional
\(
    \phi:[0,T]\times V\times C([0,T],U)\to\mathbb R
\)
is said to be non-anticipative if, for every \(t\in[0,T]\), every \(x\in V\),
and every \(\gamma,\nu\in C([0,T],U)\) such that
\(
    \gamma(r)=\nu(r),\ 0\leq r\leq t,
\)
one has
\(
    \phi(t,x,\gamma)=\phi(t,x,\nu).
\)
\end{definition} 

Under the additional assumption that $g$ is Lipschitz continuous, we can recover local Lipschitz continuity of the value functional $v(t,x,\gamma)$ in the state variables $(x,\gamma)$ and local uniform continuity in the time variable $t$.

\begin{proposition}\label{regularity_value_functional}
The value functional \(v\) defined in \eqref{extended_value_functional} is
continuous on its domain with respect to the product metric introduced in
\eqref{product_metric_definition}.
\end{proposition}
\begin{proof}
Following the method in Theorem 2.1 in \cite{bardi1997bellman}, we start by showing that the value function is locally Lipschitz continuous in the state variables, uniformly in the time variable.
First, fix two sets of initial data $(x,\gamma)$ and
$(\tilde x,\tilde\gamma)$. Given $\epsilon>0$, choose a control
$\tilde u \in L^{\infty}([t,T],\mathbb{R}^{\ell})$ such that
\[
    v(t,\tilde{x},\tilde{\gamma})
    \geq
    J(t,\tilde{x},\tilde{\gamma},\tilde{u})-\epsilon .
\]
By Proposition 4.3 and the coercivity of the pseudo-control cost, arguing as below, we may restrict attention to \(\varepsilon\)-optimal controls \(\tilde u\) satisfying

$$
\int_t^T |\tilde u_r|^q\,dr\le C,
$$

where \(C>0\) is independent of \(\varepsilon\). By the fractional reconstruction estimate, this yields a uniform bound

$$
\|\nu^{t,\gamma,T,\tilde u}\|_{p/\lfloor p\rfloor\text{-var};[0,T]}
+
\|\nu^{t,\tilde\gamma,T,\tilde u}\|_{p/\lfloor p\rfloor\text{-var};[0,T]}
\le M
$$
for some \(M>0\) independent of \(\varepsilon\). We may therefore apply the stability estimate of Proposition 3.1 with a constant independent of \(\varepsilon\). Using this we can find that there exists a $C_{x, a, \lambda, b, L, T, f, g}>0$ 
\begin{align*}
&v(t, x, \gamma) - v(t, \tilde{x}, \tilde{\gamma}) \leq   J(t, x, \gamma, \tilde{u}) -  J(t, \tilde{x}, \tilde{\gamma}, \tilde{u}) + \epsilon \\
&= \int_{t}^{T} \Big( f(X^{t, x, \nu, \tilde{u}}_r, \nu^{t, \gamma, T, \tilde{u}}_r, \tilde{u}_r) - f(X^{t, \tilde{x}, \tilde{\nu}, \tilde{u}}_r, \nu^{t, \tilde{\gamma},T, \tilde{u}}_r, \tilde{u}_r) \Big) dr \\
&\quad + g(X^{t, x, \nu, \tilde{u}}_T) - g(X^{t, \tilde{x}, \tilde{\nu}, \tilde{u}}_T) +\epsilon\\
&\le C_{x, a, \lambda, b, L, T, f, g}(T+1)\left(|x - \tilde{x}| + |\gamma_0 - \tilde{\gamma}_0| + \|\gamma - \tilde{\gamma}\|_{p / \floor p; [0, t]}\right) + \epsilon,
\end{align*}
where in the last step we used the Lipschitz continuity of  $f$ and $g$ in the state variables, alongside estimate in Lemma \ref{lem:fractional_continuation_estimates} and the stability estimates in Proposition \ref{existence_uniqueness_stability_control_process}. The proof of this part is then concluded by the fact that the previous inequality holds for any $\epsilon >0$ and, by appropriately choosing the multiplicative constant in the last line, is symmetric with respect to the pairs $(x, \gamma)$ and $(\tilde{x}, \tilde{\gamma})$.\newline
For the second part, we establish the continuity of the value functional with
respect to the time variable.

Fix $k\geq 1$ and set
\[
    K:=[0,T]\times \overline{B(0,k)}\times AC^\alpha_k.
\]
By Proposition \ref{boundedness_val_fun}, the value functional is bounded on
$K$. This observation, together with the coercivity of the auxiliary cost
functional $\tilde J$, allows us to restrict the set of admissible controls to
those satisfying
\begin{equation}
    \int_0^T |u_r|^q\,dr\leq \tilde C,
    \label{coercivity_estimate}
\end{equation}
where $\tilde C$ depends only on $K$ and the coefficients of the problem.

Indeed, let
\[
    v_K:=\sup_{(t,x,\gamma)\in K}v(t,x,\gamma)<\infty
\]
and set
\(
    \tilde p:=p\lfloor p\rfloor\vee\frac{\kappa}{\kappa-1}.
\)
Fix $(t,x,\gamma)\in K$ and $\epsilon\in(0,1)$, and let
$u\in L^\infty([t,T],\mathbb R^\ell)$ be an $\epsilon$-optimal control.
Then
\[
    J(t,x,\gamma,u)
    \leq
    v(t,x,\gamma)+\epsilon
    \leq
    v_K+1.
\]
On the other hand, by the coercive lower bound on $f$ and
\eqref{bound_rough_integral_psi_u},
\begin{align*}
    J(t,x,\gamma,u)
    &\geq
    \int_t^T f_0|u_r|^q\,dr
    -C_f
    -C_{\lambda,b,p,L,T,g,k}
    \left(
        1+\int_t^T|u_r|^{\tilde p}\,dr
    \right) \\
    &\geq
    \int_t^T f_0|u_r|^q\,dr
    -C_f
    -C_{\lambda,b,p,L,T,g,k}
    \left(
        1+
        T^{1-\frac{\tilde p}{q}}
        \left(
            \int_t^T|u_r|^q\,dr
        \right)^{\frac{\tilde p}{q}}
    \right),
\end{align*}
where the second inequality follows from H\"older's inequality. Consequently,
\[
    v_K+1
    \geq
    \int_t^T f_0|u_r|^q\,dr
    -C_f
    -C_{\lambda,b,p,L,T,g,k}
    \left(
        1+
        T^{1-\frac{\tilde p}{q}}
        \left(
            \int_t^T|u_r|^q\,dr
        \right)^{\frac{\tilde p}{q}}
    \right).
\]
Since $q>\tilde p$, the right-hand side tends to $+\infty$ as
$\int_t^T|u_r|^q\,dr\to\infty$. Hence there exists $\tilde C>0$,
independent of $(t,x,\gamma)\in K$ and of $\epsilon\in(0,1)$, such that every
$\epsilon$-optimal control satisfies
\[
    \int_t^T|u_r|^q\,dr\leq\tilde C.
\]
Thus, without changing the value functional, the infimum may be restricted
to controls satisfying \eqref{coercivity_estimate}, where controls defined on
$[t,T]$ are identified with their zero extensions to $[0,T]$.
Then, by H\"older's inequality, for any $s\geq t$,
\begin{align}
|\nu_s^{t,\gamma,T,u}-\gamma_s|
&\leq
\frac{1}{\Gamma(\alpha)}
\int_t^s
\frac{|u_r|}{(s-r)^{1-\alpha}}\,dr
+
\frac{1}{\Gamma(\alpha)}
\int_t^s
\frac{|D^\alpha_{0^+}(\gamma^{a,u}-a)(r)|}
     {(s-r)^{1-\alpha}}\,dr
\notag\\
&\leq
\frac{1}{\Gamma(\alpha)}
\left(\int_t^s|u_r|^q\,dr\right)^{1/q}
\left(
\int_t^s
(s-r)^{\frac{q(\alpha-1)}{q-1}}\,dr
\right)^{\frac{q-1}{q}}
+
\frac{2k}{\Gamma(\alpha+1)}(s-t)^\alpha
\notag\\
&\leq
\frac{\tilde C^{1/q}}{\Gamma(\alpha)}
\left(
\frac{q-1}{q\alpha-1}
\right)^{\frac{q-1}{q}}
(s-t)^{\frac{q\alpha-1}{q}}
+
\frac{2k}{\Gamma(\alpha+1)}(s-t)^\alpha .
\label{bound_nu}
\end{align}
Analogously, for any $z>s\geq t$ and any admissible control $u$
satisfying \eqref{coercivity_estimate}, since $\alpha-1<0$, we have
\(
    (z-r)^{\alpha-1}
    \leq
    (s-r)^{\alpha-1},
\) for $r\in[t,s]$.
Consequently,
\begin{align}
|\nu_z^{s,\gamma,T,u}-\nu_z^{t,\gamma,T,u}|
&\leq
\frac{1}{\Gamma(\alpha)}
\int_t^s
\frac{|D^\alpha_{0^+}(\gamma^{a,u}-a)(r)|}
     {(z-r)^{1-\alpha}}\,dr
+
\frac{1}{\Gamma(\alpha)}
\int_t^s
\frac{|u_r|}
     {(z-r)^{1-\alpha}}\,dr
\notag\\
&\leq
\frac{2k}{\Gamma(\alpha+1)}(s-t)^\alpha
+
\frac{\tilde C^{1/q}}{\Gamma(\alpha)}
\left(
    \frac{q-1}{q\alpha-1}
\right)^{\frac{q-1}{q}}
(s-t)^{\frac{q\alpha-1}{q}}.
\label{bound_nu_2}
\end{align}
From this last result, the Lipschitz continuity of the value functional in the state variables derived above and the regularity in time of RDE solutions (Proposition \ref{existence_uniqueness_stability_control_process}), it follows that
\begin{equation}
    |v(s, X^{t, x, \nu, u}_s, \nu_s^{t, \gamma, T, u}) - v(s, x, \gamma)|\leq C_{\tilde{C}, \alpha, q, \lambda, b, L, T, f, g}\left((s-t)^{\frac{q\alpha - 1}{q}} + \|\boldsymbol{\zeta}\|_{p; [t,s]}  \right). \label{continuity_value_function_1}
\end{equation}
Using again the time regularity of RDE solutions with the estimate \eqref{bound_nu} we obtain that for any admissible control $u$ such that $\int_t^T |u_r|^qdr < \tilde{C}$ and every \hbox{$(t, x, \gamma) \in K$}  there exists a positive constant $M$ such that
\[
\|X^{t, x, \gamma, u}\|_{\infty; [t,T]}, \|\nu^{ t, \gamma, T, u}\|_{\infty; [t, T]} \leq M.\]

On the other hand, applying the DPP with the control $u=0$ on $[t,s]$
gives
\begin{align*}
v(t,x,\gamma)
&\leq
\int_t^s
f\big(
    X_r^{t,x,\nu,0},
    \nu_r^{t,\gamma,T,0},
    0
\big)\,dr
+
v\big(
    s,
    X_s^{t,x,\nu,0},
    \nu_s^{t,\gamma,T,0}
\big).
\end{align*}
Since $f$ is bounded in the state variables when $u=0$, it follows that
\begin{align}
&v(t,x,\gamma)
-
v\big(
    s,
    X_s^{t,x,\nu,0},
    \nu_s^{t,\gamma,T,0}
\big)
\leq
C_f(s-t).
\label{continuity_value_function_2}
\end{align}
Moreover, for any $\epsilon>0$ the DPP allows us to find a control $u$ such that 
\begin{equation*}
     v(t,x,\gamma)  \geq \int_t^s f(X^{t,x, \nu, u}_r, \nu_r^{t, \gamma, T, u}, u_r) dr  -\epsilon + v(s,X^{t, x, \nu, u}_s,\nu^{t, \gamma, T, u}).
\end{equation*}
Combining this last expression with the fact that $f$ is bounded below in $K$, allows us to get for $C_{\tilde{C}, \alpha, q, \lambda, b, L, T, f, g} >0 $
\begin{equation}
     v(t,x,\gamma)   -  v(s,X^{t,x, \nu, u}_s,\nu_s^{ t, \gamma, T, u}) \geq -C_{\tilde{C}, \alpha, q, \lambda, b, L, T, f, g}\left((s-t)^{\frac{q\alpha - 1}{q}} +  \|\boldsymbol{\zeta}\|_{p; [t,s]}\right)-\epsilon .\label{continuity_value_function_3} 
\end{equation}
Combining the continuity in the state variables, uniformly with respect to
time, with \eqref{continuity_value_function_1},
\eqref{continuity_value_function_2} and
\eqref{continuity_value_function_3}, we conclude that the value functional
is jointly continuous with respect to the product metric induced by
\eqref{product_metric_definition}, concluding the proof.\end{proof}

\begin{remark}
The previous result, together with the boundedness established in
Proposition~\ref{boundedness_val_fun}, implies that the value functional is locally bounded
on its domain with respect to the product metric
\eqref{product_metric_definition}.
\end{remark}

Following the same argument as in the previous proof, we can also prove that the value functional is locally Lipschitz continuous with respect to the driving rough path in equation \eqref{smooth_fractional_control_dynamics}. 
\begin{proposition}\label{continuity_in_noise}
Let \(v^{\boldsymbol{\zeta}}\) denote the value functional
associated with the controlled rough differential equation
\eqref{smooth_fractional_control_dynamics} and the cost functional
\eqref{extended_cost_functional}. For every \(k\in\mathbb{N}\) and
\(L>0\), there exists a constant \(C_{k,L}>0\) such that, for every
\(\boldsymbol{\zeta},\boldsymbol{\eta}\in
\mathscr{C}^p([0,T],\mathbb{R}^d)\) satisfying
\(
    \|\boldsymbol{\zeta}\|_{p;[0,T]},
    \|\boldsymbol{\eta}\|_{p;[0,T]}
    \leq L,
\)
one has
\[
    \sup_{(t,x,\gamma)\in
    [0,T]\times B(0,k)\times AC^\alpha_k}
    \left|
        v^{\boldsymbol{\zeta}}(t,x,\gamma)
        -
        v^{\boldsymbol{\eta}}(t,x,\gamma)
    \right|
    \leq
    C_{k,L}
    \|\boldsymbol{\zeta}-\boldsymbol{\eta}\|_{p;[0,T]}.
\]
In particular, the value functional is locally uniformly continuous
with respect to the \(p\)-variation rough path topology, uniformly on
\([0,T]\times B(0,k)\times AC^\alpha_k\).
\end{proposition}

\begin{proof}
Fix \(k\in\mathbb{N}\) and \(L>0\), and set
\(
    K_k:=[0,T]\times B(0,k)\times AC^\alpha_k.
\)
Fix \((t,x,\gamma)\in K_k\) and \(\epsilon\in(0,1)\).
Proceeding as in the proof of Proposition
\ref{regularity_value_functional}, there exists a control
\(\tilde{u}\in L^\infty([t,T],\mathbb{R}^\ell)\) such that
\[
    v^{\boldsymbol{\eta}}(t,x,\gamma)
    \geq
    J^{\boldsymbol{\eta}}(t,x,\gamma,\tilde{u})
    -
    \epsilon.
\]

Moreover, the boundedness result in Proposition
\ref{boundedness_val_fun}, together with the coercivity of the
pseudo-control cost, implies, exactly as in the proof of Proposition
\ref{regularity_value_functional}, that there exists a constant
\(\widetilde C_{k,L}>0\), independent of
\((t,x,\gamma)\in K_k\), \(\epsilon\in(0,1)\), and
\(\boldsymbol{\eta}\) satisfying
\(\|\boldsymbol{\eta}\|_{p;[0,T]}\leq L\), such that
\[
    \int_t^T |\tilde{u}_r|^q\,dr
    \leq
    \widetilde C_{k,L}.
\]
By Proposition~\ref{cost_function_fractional}, this yields a uniform
bound on the \(p/\lfloor p\rfloor\)-variation norm of the reconstructed
control \(\nu^{t,\gamma,T,\tilde{u}}\).

Let \(X^{t,x,\nu,\tilde{u}}\) and
\(Y^{t,x,\nu,\tilde{u}}\) denote the trace of the  solutions of
\eqref{smooth_fractional_control_dynamics} driven by
\(\boldsymbol{\zeta}\) and \(\boldsymbol{\eta}\), respectively, with
the same reconstructed control
\(\nu^{t,\gamma,T,\tilde{u}}\). Then
\begin{align*}
&v^{\boldsymbol{\zeta}}(t,x,\gamma)
-
v^{\boldsymbol{\eta}}(t,x,\gamma)\\
&\leq
J^{\boldsymbol{\zeta}}(t,x,\gamma,\tilde{u})
-
J^{\boldsymbol{\eta}}(t,x,\gamma,\tilde{u})
+
\epsilon\\
&=
\int_t^T
\Big(
f(X^{t,x,\nu,\tilde{u}}_r,
  \nu^{t,\gamma,T,\tilde{u}}_r,
  \tilde{u}_r)
-
f(Y^{t,x,\nu,\tilde{u}}_r,
  \nu^{t,\gamma,T,\tilde{u}}_r,
  \tilde{u}_r)
\Big)\,dr\\
&\quad
+
g(X^{t,x,\nu,\tilde{u}}_T)
-
g(Y^{t,x,\nu,\tilde{u}}_T)
+
\epsilon.
\end{align*}
Using the fact that \(g\) is Lipschitz continuous and that the running
cost \(f\) is the sum of a function which is Lipschitz continuous in
the state variables and a function depending exclusively on the
pseudo-control, we obtain
\begin{align*}
&v^{\boldsymbol{\zeta}}(t,x,\gamma)
-
v^{\boldsymbol{\eta}}(t,x,\gamma)\\
&\leq
C_{f,g}
\left(
\int_t^T
\left|
X^{t,x,\nu,\tilde{u}}_r
-
Y^{t,x,\nu,\tilde{u}}_r
\right|\,dr
+
\left|
X^{t,x,\nu,\tilde{u}}_T
-
Y^{t,x,\nu,\tilde{u}}_T
\right|
\right)
+
\epsilon.
\end{align*}
The stability property for solutions to RDEs in Proposition
\ref{existence_uniqueness_stability_control_process}, together with
the uniform bound on the reconstructed control obtained above, implies
that there exists a constant \(C_{k,L}>0\), independent of
\((t,x,\gamma)\in K_k\) and \(\epsilon\in(0,1)\), such that
\[
    v^{\boldsymbol{\zeta}}(t,x,\gamma)
    -
    v^{\boldsymbol{\eta}}(t,x,\gamma)
    \leq
    C_{k,L}
    \|\boldsymbol{\zeta}-\boldsymbol{\eta}\|_{p;[0,T]}
    +
    \epsilon.
\]
Letting \(\epsilon\to0\) and exchanging the roles of
\(\boldsymbol{\zeta}\) and \(\boldsymbol{\eta}\), we obtain
\[
    \left|
        v^{\boldsymbol{\zeta}}(t,x,\gamma)
        -
        v^{\boldsymbol{\eta}}(t,x,\gamma)
    \right|
    \leq
    C_{k,L}
    \|\boldsymbol{\zeta}-\boldsymbol{\eta}\|_{p;[0,T]}.
\]
Since the constant is independent of
\((t,x,\gamma)\in K_k\), taking the supremum over \(K_k\) concludes the
proof.
\end{proof}

We conclude this section with a brief application of the preceding continuity
result to stochastic pathwise control. The purpose of this discussion is not to
develop a full theory of stochastic control with fractional dynamics, but rather
to explain how the deterministic pathwise formulation introduced above can be
combined with a random rough path driver in such a way that the resulting pathwise problem admits a natural stochastic interpretation.

We emphasize that the control problem studied in the previous sections is
pathwise: the driving rough path is fixed, and the value
\(v^{\boldsymbol\zeta}(t,x,\gamma)\) is defined for each deterministic rough
path \(\boldsymbol\zeta\). When the driving signal is random, there are two
distinct stochastic objects one may associate with this pathwise problem.

The first one is the randomized pathwise value, obtained by freezing the
realization of the random rough path \(\mathbf B\), solving the corresponding
deterministic pathwise control problem, and then averaging:
\[
    \mathbb E\big[v^{\mathbf B}(t,x,\gamma)\big]
    =
    \mathbb E\left[
        \inf_{u\in L^\infty([t,T],\mathbb R^\ell)}
        J^{\mathbf B}(t,x,\gamma,u)
    \right].
\]
This is not, in general, the same as the classical stochastic control value
with adapted controls. Indeed, in the stochastic
control problem the control must be chosen non-anticipatively whilst no such restriction is imposed in this case. Thus the
quantity above should be understood as a randomized pathwise value, rather than
as the value of the corresponding adapted stochastic control problem.

In particular, when the driving signal is given by a Gaussian rough path
\(\mathbf B\), the pathwise value functional gives rise to a well-defined
expected value. The next proposition makes this precise.
\begin{proposition}\label{prop_5_6}
Let \((\Omega,\mathcal F,\mathbb P)\) be a filtered probability space, and let
\(B=(B^1,\ldots,B^d)\) be a continuous centred Gaussian process satisfying the
Friz--Victoir assumptions \cite{friz2010differential} ensuring the existence of its geometric \(p\)-rough
path lift \(\mathbf B\). Then the stochastic pathwise problem \begin{align*}
\overline{v}(t,x,\gamma) &= \mathbb{E}\left[ \inf_{u \in L^{\infty}\left([t, T], \mathbb{R}^\ell \right)} J^{\boldsymbol{\zeta}}(t,x,\gamma, u)\Bigg|_{\boldsymbol{\zeta} = \mathbf{B}(\omega)} \right]\\
&= \mathbb{E}\left[ \left.v^{\boldsymbol{\zeta}}(t,x,\gamma)\right|_{\boldsymbol{\zeta} = \mathbf{B}(\omega)} \right]\\
\end{align*}
is well-defined.
\end{proposition}
\begin{proof}
Since \(\omega\mapsto \mathbf B(\omega)\) is a measurable map into the space of geometric rough paths (Section 6.3 in \cite{friz2010differential}) and, for fixed
\((t,x,\gamma)\), the map
\[
\boldsymbol\zeta\mapsto v^{\boldsymbol\zeta}(t,x,\gamma)
\]
is continuous, hence Borel measurable, the composition
\[
\omega\mapsto v^{\mathbf B(\omega)}(t,x,\gamma)
\]
is measurable. Moreover, repeating the proof of
Proposition~\ref{boundedness_val_fun} while keeping explicit the
dependence on the driving rough path in the estimates of
Lemma~\ref{lemma_3_6}, one obtains, for fixed $(t,x,\gamma)$, constants
$C>0$ and $m\in\mathbb{N}$, independent of
$\boldsymbol{\zeta}$, such that
\[
    \left|v^{\boldsymbol{\zeta}}(t,x,\gamma)\right|
    \leq
    C\left(
        1+\|\boldsymbol{\zeta}\|_{p;[0,T]}^m
    \right).
\] Applying this with
\(\boldsymbol\zeta=\mathbf B(\omega)\), and using Fernique's theorem for
Gaussian rough paths (Corollary 30 in \cite{friz2010differential}), gives
\[
\mathbb E\left[
\left|v^{\mathbf B}(t,x,\gamma)\right|
\right]<\infty.
\]
Therefore the expected value is well-defined.\end{proof}

The corresponding non-anticipative stochastic value is
\[
    \inf_{\nu\in\mathcal A}
    \mathbb E\left[J^{\mathbf B}(t,x,\gamma,\nu)\right],
\]
where \(\mathcal A\) denotes the class of progressively measurable bounded
pseudo-controls taking values in $\mathbb{R}^\ell$. In general,
\[
    \mathbb E\left[
        \inf_{u\in L^\infty([t,T],\mathbb R^\ell)}
        J^{\mathbf B}(t,x,\gamma,u)
    \right]
    \leq
    \inf_{\nu\in\mathcal A}
    \mathbb E\left[J^{\mathbf B}(t,x,\gamma,\nu)\right].
\]
This inequality follows from the interchange between pathwise minimisation and
expectation, together with the restriction to the admissible class
\(\mathcal A\).

We record the following dual representation as an illustration of how the
pathwise formulation interfaces with adapted stochastic control. It shows that
the adapted stochastic value can be recovered from pathwise minimisation by
introducing suitable penalties whose expectation is non-positive for adapted
controls. These penalties offset the
advantage gained by optimizing path by path over controls outside the admissible
class \(\mathcal A\). This can be viewed as an analogue of the duality principle of
\cite{diehl2017stochastic} in the present fractional-control framework.
In \cite{diehl2017stochastic}, the control enters only the drift coefficient;
the authors point out that allowing an arbitrary control to enter the
coefficient of the rough driver may lead to a degenerate pathwise control
problem. Controlled diffusion coefficients have also been considered
by \cite{henrylabordere2016dual} in the classical Brownian setting. Their
approach relies on a time discretisation of the admissible controls and
recovers the stochastic value as the limit of the corresponding dual
problems, rather than through a continuous pathwise rough control problem.

More recently, \cite{bank2026duality} established dual representations for
controlled It\^o diffusions using martingale penalties and the associated
Hamilton--Jacobi--Bellman equation; they also observe that the pathwise
stochastic optimization appearing in one of their representations may
degenerate when the diffusion coefficient is controlled.
The recent work \cite{lauriere2026dual} instead considers diffusion
coefficients independent of the control and focuses on the numerical
solution of the resulting dual problem using SPDE and pathwise approximation
methods.

The purpose of Proposition~\ref{prop_dual_representation} is different.
Here the duality is formulated directly at the level of the rough-pathwise
solution map. In particular, the pseudo-control acts through the fractional
control path $\gamma$, so that the coefficient of the rough driver may depend
on $\gamma$, while the resulting dynamics remain well defined pathwise for
each $\boldsymbol\zeta\in\mathscr C^p$. The proposition therefore records
how the information-relaxation argument of \cite{diehl2017stochastic}
interfaces with the particular controlled rough dynamics considered in this
paper, without appealing to a time discretisation or to It\^o calculus.
\begin{proposition}[Dual representation]\label{prop_dual_representation}
Let the filtered probability space \((\Omega,\mathcal F,\mathbb P)\) and \(\mathbf B\) be as above. For
\(t\in[0,T]\), let \(\mathcal A_t\) denote the set of progressively measurable
pseudo-controls \(w:[t,T]\to\mathbb R^\ell\) which are essentially bounded.

Let \(\mathcal Z_{\mathcal F}\) be the class of all maps
\[
z:\mathscr C^{p}([0,T],\mathbb R^d)
\times L^\infty([t,T],\mathbb R^\ell)
\to \mathbb R
\]
such that \(z\) is bounded, measurable, and continuous in
\(\boldsymbol\zeta\in\mathscr C^{p}([0,T],\mathbb R^d)\), uniformly over
\(w\in L^\infty([t,T],\mathbb R^\ell)\), and such that
\[
\mathbb E\bigl[z(\mathbf B,w)\bigr]\le 0
\]
for every adapted pseudo-control \(w\in\mathcal A_t\). Then
\[
\inf_{w \in\mathcal A_t}
\mathbb E\left[
 J^{\mathbf{B}}(t,x,\gamma, w)
\right]
=
\sup_{z\in\mathcal Z_{\mathcal F}}
\mathbb E\left[
\left.
\inf_{u\in L^\infty([t,T],\mathbb R^\ell)}
\left\{
J^{\boldsymbol{\zeta}}(t,x,\gamma, u)
+
z(\boldsymbol\zeta,u)
\right\}
\right|_{\boldsymbol\zeta=\mathbf B(\omega)}
\right].
\]
\end{proposition}
\begin{proof}
The result follows from Theorem~12 from \cite{diehl2017stochastic}. In the notation of that
theorem, the underlying noise is the Gaussian
rough path lift \(\mathbf B\). The pathwise cost is
\[
    (\boldsymbol\zeta,u)
    \mapsto
    J^{\boldsymbol\zeta}(t,x,\gamma,u),
\]
and the admissible class is the set \(\mathcal A_t\) of progressively
measurable essentially bounded pseudo-controls. The measurability and
integrability assumptions required in Theorem~12 from \cite{diehl2017stochastic} follow from
Proposition~\ref{prop_5_6}. With these identifications, the class \(\mathcal Z_F\) is
precisely the corresponding class of dual feasible penalties, and the proof is
unchanged.
\end{proof}

\section{Well-posedness of the rough fractional HJB equation}\label{section_well}

In this section we associate a Hamilton–Jacobi–Bellman (HJB) equation to the value functional $v^{\boldsymbol{\zeta}}(t,x,\gamma)$ and introduce an appropriate notion of solution. 
To establish these results, we temporarily restrict to the case where the driving path belongs to $C^\infty$. To avoid confusion with the rough path $\boldsymbol{\zeta}$, we denote such a path by $\eta$ and the corresponding value functional as $v^\eta(t, x, \gamma)$. 
In this setting we can rely on the viscosity solution theory for systems driven by fractional dynamics, developed in the series of works \cite{gomoyunov2020dynamic, gomoyunov2020theory, gomoyunov2021viscosity}, which we will extend to our framework. 

A key distinction between viscosity theory for fractional dynamics and the classical viscosity framework is the need for a derivative  applicable when one or more state variables are paths in the class $AC^\alpha$. 
For this purpose we follow \cite{gomoyunov2020dynamic}, where the concept of co-invariant derivatives of fractional type is introduced. These are functional derivatives that originate in the stability theory of retarded functional differential equations and have been extensively analysed in \cite{kim1999functional}. 
Their defining property is that, when evaluated at a point $(t,\gamma) \in [0,T] \times C([0,T],\mathbb{R}^\ell)$, the co-invariant derivative of a functional depends only on the restriction of $\gamma$ to $[0,t]$. Formally,  
\begin{definition}\label{ci_smooth_definition}
Let $t \in [0, T), x,y \in \mathbb{R}^e$ and $\gamma \in AC^{\alpha}([0,T], \mathbb{R}^\ell)$. A functional $\varphi: [0, T] \times \mathbb{R}^e \times AC^{\alpha}([0, T], \mathbb{R}^\ell) \rightarrow \mathbb{R}$ is said to be ci-differentiable of  order $\alpha$ at \hbox{($t, x, \gamma$)} if for every $\nu \in AC^{\alpha}([0, T], \mathbb{R}^\ell)$ such that $\nu(s) = \gamma(s)$ for every $s \in [0, t]$ there exist   $\frac{\partial^\alpha}{\partial t} \varphi(t, x, \gamma) \in \mathbb{R}, \nabla_x^{\alpha} \varphi(t,x, \gamma) \in \mathbb{R}^e, \nabla_\gamma^{\alpha} \varphi(t,x, \gamma) \in \mathbb{R}^\ell$ such that the following holds for any $z \in (t, T)$
\begin{align*}
    \varphi(z, y,  \nu)  - \varphi(t, x, \gamma) &= \frac{\partial^\alpha}{\partial t} \varphi(t, x, \gamma) (z - t)+  \langle \nabla^\alpha_x \varphi(t, x, \gamma), (y-x)\rangle  \\ 
    &\quad +  \langle\nabla^{\alpha}_\gamma \varphi(t,x, \gamma),\left(I_{0^+}^{1- \alpha} (\nu - \gamma_0)(z) - I_{0^+}^{1-\alpha} (\gamma - \gamma_0) (t)\right)\rangle + \small o (z-t + \|x-y\|), 
\end{align*}
where the remainder is allowed to depend on $(t, x, \gamma, \nu)$.
\end{definition}
The preceding definition determines the derivatives uniquely. As we shall prove, the combination of this definition with the Dynamic Programming Principle allows us to associate to $v^\eta$ the following problem

\begin{equation} \label{HJB_equation}
    \begin{cases}
    - \frac{\partial^\alpha}{\partial t} v^\eta(t, x, \gamma) - \langle \nabla^\alpha_x v^\eta(t, x, \gamma) , b(x,\gamma_t) + \lambda(x,\gamma_t)\dot{\eta}_t \rangle\\
    \quad\quad\quad + H(x, \gamma_t, \nabla^\alpha_\gamma v^\eta(t, x, \gamma))  = 0, & (t,x,\gamma)\in [0,T)\times\mathbb R^e\times AC^\alpha, \\
    v^\eta(T, x, \gamma) = g(x), &  (x,\gamma)\in\mathbb R^e\times AC^\alpha, 
    \end{cases}
\end{equation}
where $H(x, \gamma, \phi) = \sup\limits_{u \in \mathbb{R}^\ell}\{ - \langle \phi, u  \rangle - f(x, \gamma, u)\}$. \newline

Equation~\eqref{HJB_equation} should be viewed as a fractional,
path-dependent analogue of the classical Hamilton--Jacobi--Bellman equation
(see, for instance, \cite{bardi1997optimal}). Its nonlocal character comes
from the fractional representation of the control. In particular, the value
function is not defined only in terms of the current variables
\((t,x,\gamma_t)\), but depends on the whole past control path \(\gamma\).
Variations in the control are therefore measured through $\nu$, so that the equation retains
memory of the past evolution of the control. This contrasts with the classical
HJB equation, where the value function depends only on the current state and
the associated differential operator is local in time.

This places \eqref{HJB_equation} in the broader literature on path-dependent
Hamilton--Jacobi equations. In this setting, the value function depends not
only on the current state, but also on an underlying trajectory, and the
corresponding dynamic programming equation is formulated on a space of paths;
see, for instance,
\cite{gomoyunov2021path,BAYRAKTAR20182096,bayraktar2025viscosity}.

A particularly relevant subclass is given by fractional Hamilton--Jacobi
equations, where fractional derivatives encode memory effects and naturally
lead to value functionals on path spaces. This perspective has been developed,
for example, in
\cite{gomoyunov2020dynamic,gomoyunov2021viscosity}.

The equation derived in this paper is related to this fractional
path-dependent framework, but arises from a different mechanism. Here the
fractional structure is introduced through the regularity imposed on the
admissible controls, rather than through a fractional state equation.\newline

We now turn to the notion of solution for the Hamilton--Jacobi--Bellman
equation derived above. Since the value functional is not expected, in general,
to possess the regularity required to satisfy the equation pointwise, we seek
solutions in the viscosity sense. To this end, we first specify a class of test
functions adapted to the fractional, path-dependent structure of the problem.
These test functions play the role of the smooth functions used in the
classical viscosity solution framework. The choice is the class of  
ci-smooth functions of order $\alpha$, which are the ci-differentiable functions of order $\alpha$ for which the following conditions are met:
\begin{enumerate}
    \item $\varphi$ is ci-differentiable at every point $(t, x, \gamma) \in [0,T) \times \mathbb{R}^e \times AC^\alpha([0, T], \mathbb{R}^\ell)$,
    \item $\varphi$ and the functionals $\frac{\partial^\alpha}{\partial t} \varphi: [0, T) \times \mathbb{R}^e \times AC^{\alpha}([0,T], \mathbb{R}^\ell) \rightarrow \mathbb{R}$, $\nabla_\gamma^{\alpha} \varphi: [0, T) \times \mathbb{R}^e \times  AC^{\alpha}([0,T], \mathbb{R}^\ell) \rightarrow \mathbb{R}^\ell$ and $\nabla_x^{\alpha} \varphi: [0, T) \times \mathbb{R}^e \times  AC^{\alpha}([0,T], \mathbb{R}^\ell) \rightarrow \mathbb{R}^e$ are continuous with respect to the product metric defined above.
\end{enumerate}
It is important to notice that, by construction, any ci-differentiable functional is non-anticipative.

The notion of viscosity solution to the fractional HJB equation
\eqref{HJB_equation}, adapted from
\cite{gomoyunov2021viscosity} to our framework, is given below.
Since $\eta$ is fixed, for every $n,k\in\mathbb N$, set
\[
m_{n,k}:=n\vee k,
\qquad
R_k(t)
:=
k+\|b\|_\infty t
+\|\lambda\|_\infty\int_0^t|\dot\eta_r|\,dr,
\]
and define
\[
\mathcal K_k^n
:=
\left\{
(t,x,\gamma)\in
[0,T]\times\mathbb R^e\times AC^\alpha:
|x|\leq R_k(t),\quad
\gamma\in AC_{m_{n,k}}^\alpha
\right\}.
\]
For every fixed $n$, the family
$(\mathcal K_k^n)_{k\in\mathbb N}$ is a compact exhaustion of the
time-state space and is invariant under stopped continuations
generated by controls satisfying $|u|\leq n$.

\begin{definition}[Viscosity solution]
\label{definition_viscosity_solution}

For every $n\in\mathbb N$, define
\[
H_n(x,y,\phi)
:=
\sup_{|u|\leq n}
\bigl\{
-\langle \phi,u\rangle-f(x,y,u)
\bigr\}.
\]
A continuous functional
\(
v_n:
[0,T]\times\mathbb R^e
\times AC^\alpha([0,T],\mathbb R^\ell)
\to \mathbb R
\)
is called a viscosity solution of the $n$-truncated problem \eqref{HJB_equation} if
\[
v_n(T,x,\gamma)=g(x)
\]
and the following conditions hold. For every $k\in\mathbb N$ and every ci-smooth functional
\(
\varphi:
[0,T]\times\mathbb R^e
\times AC^\alpha([0,T],\mathbb R^\ell)
\to \mathbb R,
\)
if $v_n-\varphi$ attains a local maximum on
$\mathcal K_k^n$ at some point
\(
(t,x,\gamma)\in \mathcal K_k^n,
\) with $t < T$,
then
\[
-\frac{\partial^\alpha}{\partial t}\varphi(t,x,\gamma)
-
\left\langle
\nabla_x^\alpha\varphi(t,x,\gamma),
b(x,\gamma_t)+\lambda(x,\gamma_t)\dot\eta_t
\right\rangle
+
H_n\left(
x,\gamma_t,\nabla_\gamma^\alpha\varphi(t,x,\gamma)
\right)
\leq0.
\]
If instead $v_n-\varphi$ attains a local minimum on
$\mathcal K_k^n$ at some point
\(
(t,x,\gamma)\in\mathcal K_k^n,
\) with $t < T$, 
then
\[
-\frac{\partial^\alpha}{\partial t}\varphi(t,x,\gamma)
-
\left\langle
\nabla_x^\alpha\varphi(t,x,\gamma),
b(x,\gamma_t)+\lambda(x,\gamma_t)\dot\eta_t
\right\rangle
+
H_n\left(
x,\gamma_t,\nabla_\gamma^\alpha\varphi(t,x,\gamma)
\right)
\geq0.
\]

A continuous functional
\(
v:
[0,T]\times\mathbb R^e
\times AC^\alpha([0,T],\mathbb R^\ell)
\to \mathbb R
\)
is called a viscosity solution of \eqref{HJB_equation} if there exists a sequence $(v_n)_{n\in\mathbb N}$
such that:

\begin{enumerate}
    \item for every $n\in\mathbb N$, $v_n$ is a viscosity solution of
    the $n$-truncated problem;

    \item $v_n\to v$ locally uniformly.
\end{enumerate}

In particular, $v(T,x,\gamma)=g(x)$ for every
$(x,\gamma)\in
\mathbb R^e\times AC^\alpha([0,T],\mathbb R^\ell)$.
\end{definition}

Viscosity subsolutions and supersolutions of the \(n\)-truncated
problem are defined analogously by retaining, respectively, the
local-maximum inequality and the local-minimum inequality above,
together with the corresponding terminal inequality.

With this definition in place, we can now give a precise meaning to the following problem associated with the original value function $v^{\boldsymbol{\zeta}}$: 
\begin{equation} \label{rough_HJB_equation}
    \begin{cases}
    - \frac{\partial^\alpha}{\partial t} v^{\boldsymbol{\zeta}}(t, x, \gamma)dt - \langle \nabla^\alpha_x v^{\boldsymbol{\zeta}}(t, x, \gamma) , b(x,\gamma_t)dt + \lambda(x,\gamma_t)d\boldsymbol{\zeta}_t \rangle\\
    \quad\quad\quad + H(x, \gamma_t, \nabla^\alpha_\gamma v^{\boldsymbol{\zeta}}(t, x, \gamma))dt  = 0, & (t,x,\gamma)\in [0,T)\times\mathbb R^e\times AC^\alpha, \\
    v^{\boldsymbol{\zeta}}(T, x, \gamma) = g(x), & (x,\gamma)\in\mathbb R^e\times AC^\alpha, 
    \end{cases}
\end{equation}
where 
\[
H(x, \gamma, \phi) = \sup_{u \in \mathbb{R}^\ell}\{ - \langle \phi, u \rangle - f(x, \gamma, u)\}.
\]

This interpretation is formalized in the following definition, which we extend to our framework from the work \cite{caruana2011rough}.

\begin{definition}[Rough viscosity solution]
\label{rough_viscosity_solution}

Let
\(
v^{\boldsymbol{\zeta}}:
[0,T]\times\mathbb R^e
\times AC^\alpha([0,T];\mathbb R^\ell)
\to\mathbb R
\)
be continuous with respect to the topology induced by
\eqref{product_metric_definition}.

We say that \(v^{\boldsymbol{\zeta}}\) is a rough viscosity solution of
\eqref{rough_HJB_equation} if, for every sequence of smooth paths
\((\zeta_m)_{m\in\mathbb N}\) whose canonical lifts converge to
\(\boldsymbol{\zeta}\) in the \(p\)-variation topology, there exists a
sequence \((v^m)_{m\in\mathbb N}\) such that \(v^m\) is a viscosity
solution, in the sense of Definition~\ref{definition_viscosity_solution}, of the
fractional HJB equation driven by \(\zeta_m\), and
\[
\sup_{(t,x,\gamma)\in
[0,T]\times B(0,k)\times AC^\alpha_k}
|v^m(t,x,\gamma)-v^{\boldsymbol{\zeta}}(t,x,\gamma)|
\to0
\]
for every \(k\in\mathbb N\).
\end{definition}

This definition makes precise the sense in which
\eqref{rough_HJB_equation} combines the rough HJB equations studied in
\cite{diehl2017stochastic,allan2020pathwise} with a fractional, nonlocal
structure induced by the control variable.
\newline

To summarize, the discussion in this section provides a clear roadmap for linking the value functional to its associated HJB equations in the rough setting. By first introducing the HJB equation corresponding to the approximated value problem and clarifying the appropriate notion of viscosity solution, we establish the framework needed to connect the value functional with its PDE characterisation. This step is crucial: in order to ultimately address well-posedness for the full rough problem, we must first formulate and analyse the HJB equation for the approximated problem and prove well-posedness there. Once this foundation is in place, the value functional associated with
the rough driver can be identified as a rough viscosity solution of the
rough HJB equation, using the local uniform continuity with respect to the
driving rough path established in
Proposition~\ref{continuity_in_noise}.

\subsection{Fractional HJB equation for the approximated value functional}

As anticipated above, this part of the work is devoted to proving the well-posedness of the HJB equation associated with the approximated value problem and to showing that the value functional admits a characterisation as the viscosity solution of a fractional,
path-dependent HJB equation. Our approach relies on several preparatory results, after which we establish the main result and discuss its implications for the rough HJB equation.

\medskip

The next lemma builds on the fact that the definition of ci-differentiability can be rephrased in the form
\begin{equation} \label{ci-differentiable Taylor}
\begin{aligned}
    \varphi(z, y,  \nu)  - \varphi(t, x,  \gamma) 
    = & \frac{\partial^\alpha}{\partial t} \varphi(t, x, \gamma) (z - t) 
    +  \left\langle \nabla^\alpha_x \varphi(t, x, \gamma), (y-x)\right\rangle  \\
    & + \left \langle \nabla^\alpha_\gamma \varphi(t, x, \gamma), \int_{t}^{z} D^{\alpha}_{0^+}(\nu - \nu_0)(r) dr \right\rangle  
    + o (z-t + \|y-x\|).
\end{aligned}
\end{equation}
When combined with the non-anticipative property of the functional, this reformulation allows us to deduce Lipschitz continuity for functionals of the form
\[
\vartheta(t) := \varphi(t, X_t, \nu^t),
\]
where for any $\gamma \in  C([0, T], \mathbb{R}^\ell)$, the notation $\gamma^t$ denotes the path that coincides with $\gamma$ up to time $t$. The next lemma is a close variant of Lemma~9.2 in \cite{gomoyunov2020dynamic}, adapted to the present framework, and establishes the regularity of $\vartheta$.
\begin{lemma}\label{classic_differentiability_ci_functional} Let \( \varphi:[0,T]\times\mathbb R^e \times AC^\alpha([0,T],\mathbb R^\ell)\to\mathbb R \) be a ci-smooth functional, let \(X\in C^1([0,T],\mathbb R^e)\), and let \(\gamma\in AC^\alpha([0,T],\mathbb R^\ell)\). For each \(t\in[0,T)\), let \(\nu^t\in AC^\alpha([0,T],\mathbb R^\ell)\) be such that \( \nu^t_r=\gamma_r, \, 0\leq r\leq t. \) Then the map \( t\mapsto \vartheta(t) := \varphi(t,X_t,\nu^t) \) is locally Lipschitz continuous on \([0,T)\). \end{lemma}
\begin{proof}
Since every ci-differentiable functional is non-anticipative and
\(\nu^t_r=\gamma_r\) for \(0\leq r\leq t\), we have
\[
    \vartheta(t)
    =
    \varphi(t,X_t,\nu^t)
    =
    \varphi(t,X_t,\gamma),
    \qquad t\in[0,T).
\]
Thus it is enough to study the map
\(
    t\mapsto \varphi(t,X_t,\gamma).
\)

Fix \(T'<T\). Since \(\varphi\) is ci-smooth and
\(t\mapsto (t,X_t,\gamma)\) has compact image on \([0,T']\), there
exists \(M_{T'}>0\) such that, for every \(t\in[0,T']\),
\[
    \left|
        \frac{\partial^\alpha}{\partial t}
        \varphi(t,X_t,\gamma)
    \right|,
    \qquad
    \left|
        \nabla^\alpha_x\varphi(t,X_t,\gamma)
    \right|,
    \qquad
    \left|
        \nabla^\alpha_\gamma\varphi(t,X_t,\gamma)
    \right|
    \leq M_{T'}.
\]

Using the expansion \eqref{ci-differentiable Taylor} at the point
\((t,X_t,\gamma)\), with
\(
    (z,y,\nu)
    =
    (t+h,X_{t+h},\gamma),
\)
we obtain, for \(h>0\) sufficiently small,
\begin{align*}
    \frac{\vartheta(t+h)-\vartheta(t)}{h}
    &=
    \frac{\partial^\alpha}{\partial t}
    \varphi(t,X_t,\gamma)
    +
    \left\langle
        \nabla^\alpha_x\varphi(t,X_t,\gamma),
        \frac{X_{t+h}-X_t}{h}
    \right\rangle
    \\
    &\quad+
    \left\langle
        \nabla^\alpha_\gamma\varphi(t,X_t,\gamma),
        \frac{1}{h}
        \int_t^{t+h}
        D^\alpha_{0+}(\gamma-\gamma_0)(r)\,dr
    \right\rangle
    +
    \frac{o(h)}{h}.
\end{align*}
Since \(X\in C^1([0,T],\mathbb R^e)\),
\[
    \left|
        \frac{X_{t+h}-X_t}{h}
    \right|
    \leq
    \|\dot X\|_{\infty;[0,T']},
\]
and
\[
    \left|
        \frac{1}{h}
        \int_t^{t+h}
        D^\alpha_{0+}(\gamma-\gamma_0)(r)\,dr
    \right|
    \leq
    \|D^\alpha_{0+}(\gamma-\gamma_0)\|_{\infty;[0,T']}.
\]
Therefore, for every \(t\in[0,T')\),
\[
    \limsup_{h\downarrow0}
    \left|
        \frac{\vartheta(t+h)-\vartheta(t)}{h}
    \right|
    \leq
    L_{T'},
\]
where
\[
    L_{T'}
    :=
    M_{T'}
    \left(
        1+\|\dot X\|_{\infty;[0,T']}
        +
        \|D^\alpha_{0+}(\gamma-\gamma_0)\|_{\infty;[0,T']}
    \right).
\]
Hence the upper right Dini derivatives of both \(\vartheta\) and
\(-\vartheta\) are bounded above by \(L_{T'}\). By the standard Dini
derivative criterion, this implies
\[
    |\vartheta(t)-\vartheta(s)|
    \leq
    L_{T'}|t-s|,
    \qquad s,t\in[0,T'].
\]
Thus \(\vartheta\) is Lipschitz continuous on \([0,T']\). Since
\(T'<T\) was arbitrary, \(\vartheta\) is locally Lipschitz continuous
on \([0,T)\).
\end{proof}
\begin{proposition}\label{v_viscosity_sol}
The value functional $v^\eta$ is a viscosity solution, in the sense of
Definition \ref{definition_viscosity_solution}, of
\[
\begin{cases}
-\dfrac{\partial^\alpha}{\partial t}v(t,x,\gamma)
-\left\langle
\nabla_x^\alpha v(t,x,\gamma),
b(x,\gamma_t)+\lambda(x,\gamma_t)\dot\eta_t
\right\rangle
+H\bigl(x,\gamma_t,\nabla_\gamma^\alpha v(t,x,\gamma)\bigr)=0,
\\
\hfill
(t,x,\gamma)\in[0,T)\times\mathbb R^e\times AC^\alpha,
\\[1mm]
v(T,x,\gamma)=g(x),
\qquad
(x,\gamma)\in\mathbb R^e\times AC^\alpha,
\end{cases}
\]
where
\[
H(x,y,\phi)
:=
\sup_{u\in\mathbb R^\ell}
\bigl\{
-\langle\phi,u\rangle-f(x,y,u)
\bigr\}.
\]
\end{proposition}

\begin{proof}
For every $n\in\mathbb N$, set
\[
v_n^\eta(t,x,\gamma)
:=
\inf_{u\in L^\infty([t,T],\overline{B(0,n)})}
J(t,x,\gamma,u)
\]
and
\[
H_n(x,y,\phi)
:=
\sup_{|u|\leq n}
\bigl\{
-\langle\phi,u\rangle-f(x,y,u)
\bigr\}.
\]
We first prove that $v_n^\eta$ is a viscosity solution of the
$n$-truncated problem. Since $\eta$ is fixed, we omit it from the
notation.

Let $v_n-\varphi$ attain a local maximum on
$\mathcal K_k^n$ at
\(
(t,x,\gamma),
\) with $t < T$
and fix $\bar u\in\overline{B(0,n)}$. Take $u_r=\bar u$ on
$[t,z]$. By the invariance property of $\mathcal K_k^n$, the
corresponding stopped continuation belongs to $\mathcal K_k^n$ for
every $z>t$ sufficiently close to $t$. Therefore, the contact
inequality and the dynamic programming principle give
\[
\begin{aligned}
&\varphi(t,x,\gamma)
-\varphi\bigl(
z,X_z^{t,x,\nu,u},\nu^{t,\gamma,T,u}
\bigr)
\\
&\qquad\leq
v_n(t,x,\gamma)
-v_n\bigl(
z,X_z^{t,x,\nu,u},\nu^{t,\gamma,T,u}
\bigr)
\\
&\qquad\leq
\int_t^z
f\bigl(
X_r^{t,x,\nu,u},
\nu_r^{t,\gamma,T,u},
\bar u
\bigr)\,dr.
\end{aligned}
\]
Dividing by $z-t$, applying the ci-chain rule and letting
$z\to t$, we obtain
\[
\begin{aligned}
&-\frac{\partial^\alpha}{\partial t}\varphi(t,x,\gamma)
-\left\langle
\nabla_x^\alpha\varphi(t,x,\gamma),
b(x,\gamma_t)+\lambda(x,\gamma_t)\dot\eta_t
\right\rangle
\\
&\qquad
-\left\langle
\nabla_\gamma^\alpha\varphi(t,x,\gamma),\bar u
\right\rangle
-f(x,\gamma_t,\bar u)
\leq0.
\end{aligned}
\]
Taking the supremum over $|\bar u|\leq n$ proves the subsolution
inequality with $H_n$.

Suppose instead that $v_n-\varphi$ attains a local minimum on
$\mathcal K_k^n$ at
\(
(t,x,\gamma)
\) with $t < T$.
For $h>0$, choose an $\varepsilon h$-optimal control for the dynamic
programming principle on $[t,t+h]$. Again, its stopped continuation
belongs to $\mathcal K_k^n$. Local minimality and the dynamic
programming principle yield
\[
\begin{aligned}
&\varphi(t,x,\gamma)
-\varphi\bigl(
t+h,X_{t+h}^{t,x,\nu,u},\nu^{t,\gamma,T,u}
\bigr)
\\
&\qquad\geq
\int_t^{t+h}
f\bigl(
X_r^{t,x,\nu,u},
\nu_r^{t,\gamma,T,u},
u_r
\bigr)\,dr
-\varepsilon h.
\end{aligned}
\]
By the definition of the ci-differentiability,
\[
\begin{aligned}
-\varepsilon
\leq\frac{1}{h}\int_t^{t+h}\Bigg[
&-\frac{\partial^\alpha}{\partial t}
\varphi\bigl(
r,X_r^{t,x,\nu,u},\nu^{t,\gamma,T,u}
\bigr)
\\
&-\left\langle
\nabla_x^\alpha
\varphi\bigl(
r,X_r^{t,x,\nu,u},\nu^{t,\gamma,T,u}
\bigr),
b\bigl(
X_r^{t,x,\nu,u},\nu_r^{t,\gamma,T,u}
\bigr)
\right.
\\
&\hspace{35mm}\left.
+\lambda\bigl(
X_r^{t,x,\nu,u},\nu_r^{t,\gamma,T,u}
\bigr)\dot\eta_r
\right\rangle
\\
&-\left\langle
\nabla_\gamma^\alpha
\varphi\bigl(
r,X_r^{t,x,\nu,u},\nu^{t,\gamma,T,u}
\bigr),
u_r
\right\rangle
\\
&-f\bigl(
X_r^{t,x,\nu,u},
\nu_r^{t,\gamma,T,u},
u_r
\bigr)
\Bigg]\,dr.
\end{aligned}
\]
Since $|u_r|\leq n$, the last two terms involving $u_r$ are bounded
above by the corresponding Hamiltonian $H_n$. The continuity of
$H_n$, the ci-derivatives of $\varphi$ and the controlled
trajectories therefore allows us to let $h\to 0$, obtaining
\[
\begin{aligned}
&-\frac{\partial^\alpha}{\partial t}\varphi(t,x,\gamma)
-\left\langle
\nabla_x^\alpha\varphi(t,x,\gamma),
b(x,\gamma_t)+\lambda(x,\gamma_t)\dot\eta_t
\right\rangle
\\
&\qquad
+H_n\bigl(
x,\gamma_t,\nabla_\gamma^\alpha\varphi(t,x,\gamma)
\bigr)
\geq-\varepsilon.
\end{aligned}
\]
Letting $\varepsilon\to 0$ proves the supersolution property.
The terminal condition follows directly from the definition of
$v_n$. Thus $v_n$ is a viscosity solution of the $n$-truncated
problem.

It remains to prove convergence. Since the admissible control sets
increase with $n$,
\(
v_n^\eta\geq v_{n+1}^\eta\geq v^\eta.
\)
Conversely, given $\varepsilon>0$, choose
$u\in L^\infty([0,T],\mathbb R^\ell)$ such that
\(
J(t,x,\gamma,u)
\leq v^\eta(t,x,\gamma)+\varepsilon.
\)
For every $n\geq\|u\|_\infty$,
\[
v^\eta(t,x,\gamma)
\leq v_n^\eta(t,x,\gamma)
\leq J(t,x,\gamma,u)
\leq v^\eta(t,x,\gamma)+\varepsilon.
\]
Hence
\(
v_n^\eta\to  v^\eta
\)
pointwise. Since $v_n^\eta$ and $v^\eta$ are continuous, Dini's
theorem implies that the convergence is uniform on every compact
subset of
\(
[0,T]\times\mathbb R^e
\times AC^\alpha([0,T],\mathbb R^\ell).
\)
Moreover,
\(
H_n\to H,
\)
and the coercivity of $f$ makes this convergence locally uniform.

Consequently, $v^\eta$ is a viscosity solution in the sense of
Definition \ref{definition_viscosity_solution}.
\end{proof}
As in the classical uniqueness theory for path-dependent HJB equations,
comparison arguments rely on suitable auxiliary functionals; see, for instance,
Sections~II.5 and~III.3 of \cite{bardi1997optimal}. These functionals are the
counterpart of the penalisation terms in the doubling-of-variables method: they
force the two arguments to remain close at the maximum point. In the fractional
setting, this closeness must be measured not only at the level of the current
state, but also at the level of the histories entering the fractional dynamics.
The auxiliary functionals below are designed precisely for this purpose. In the fractional path-dependent
setting, the corresponding functionals were introduced in \cite{gomoyunov2021viscosity}. We use
the first of these functionals in the form given there, while the second one is
slightly adapted to the present formulation.

The first auxiliary functional we introduce is
\begin{equation}\label{auxiliary_functional_a}
a(\,\cdot\,|\, s, \gamma) : [0,T] \to \mathbb{R}^\ell, \quad
a(t\,|\, s, \gamma) := 
\begin{cases}
\gamma_t, & t \in [0,s],\\[2mm]
\gamma_0 + \frac{1}{\Gamma(\alpha)} \int_0^s \frac{D^\alpha_{0^+}(\gamma - \gamma_0)(r)}{(t-r)^{1-\alpha}} dr, & t \in (s,T],
\end{cases}
\end{equation}
for $(s,\gamma) \in [0,T] \times AC^\alpha([0,T], \mathbb{R}^\ell)$.

Comparing the definitions of $\nu^{s,\gamma,z,u}_t$ and $a(t\,|\, s, \gamma)$, one obtains the identity
\begin{equation}\label{relation_nu_a}
\nu_t^{s,\gamma,z,u} - \nu_t^{s,\tilde{\gamma},z,u} = \gamma_0 - \tilde{\gamma}_0 + \frac{1}{\Gamma(\alpha)} \int_0^s \frac{D^\alpha_{0^+}(\gamma-\gamma_0)(r) - D^\alpha_{0^+}(\tilde{\gamma}-\tilde{\gamma}_0)(r)}{(t-r)^{1-\alpha}} dr
= a(t\,|\, s, \gamma) - a(t\,|\, s, \tilde{\gamma}),
\end{equation}
which holds for all $s \le t \le z \le T$ and $\gamma, \tilde{\gamma} \in AC^\alpha([0,T],\mathbb{R}^\ell)$. Moreover from the proof of Lemma 3 in \cite{gomoyunov2020theory}  if \(\gamma\in AC_k^\alpha\), then
\(a(\cdot\mid t,\gamma)\in AC_k^\alpha\) for every \(t\in[0,T]\).\newline

The second auxiliary functional is a modification of the functional introduced in \cite{gomoyunov2021viscosity}, adapted to our setting:
\begin{equation}\label{auxiliary_functional}
\varpi_\epsilon(t, \gamma, \tau, \nu)
= (\epsilon^{\frac{2}{c-1}} + |a(T\,|\, t, \gamma) - a(T\,|\, \tau, \nu)|^2)^{\frac{c}{2}}
+ \int_0^T \frac{(\epsilon^{\frac{2}{c-1}} + |a(r\,|\, t, \gamma) - a(r\,|\, \tau, \nu)|^2)^{\frac{c}{2}}}{(T-r)^{(1-\alpha-\beta)c}} dr
- C_1 \epsilon^{\frac{c}{c-1}},
\end{equation}
where $c = 2/(2-\alpha)$, $0<\beta<\min(1-\alpha, \alpha/2)$, and $C_1 = 1 + \frac{T^{1-(1-\alpha-\beta)c}}{1-(1-\alpha-\beta)c}$.  
This functional extends Gomoyunov’s construction to handle the additional state variable and the unboundedness of the control-dependent Hamiltonian. The role of \(\varpi_\epsilon\) is to provide a distance-like penalisation
between two fractional histories. More precisely, it compares the extensions
\(a(\cdot\,|\,t,\gamma)\) and \(a(\cdot\,|\,\tau,\nu)\), and therefore measures
the discrepancy between the current states and the memory variables associated
with the two trajectories. This is the fractional analogue of the quadratic
penalisation used in the classical doubling-of-variables argument.\newline

We now specify the class in which uniqueness holds. A viscosity
solution $v$ belongs to this class if it admits an approximating
sequence $(v_n)_{n\in\mathbb N}$, as in Definition
\ref{definition_viscosity_solution}, such that every $v_n$ satisfies
property {\rm (L)} from \cite{gomoyunov2021viscosity}. More precisely,
for every $n,k\in\mathbb N$, there exists a constant $C_{n,k}>0$ such
that, for every $t\in[0,T]$, $x,y\in\overline{B(0,k)}$, and
$\gamma,\nu\in AC_k^\alpha$,
\[
\begin{aligned}
|v_n(t,x,\gamma)-v_n(t,y,\nu)|
\leq C_{n,k}\Bigg(
&|x-y|
+\left|a(T\mid t,\gamma)-a(T\mid t,\nu)\right|
\\
&+\int_0^T
\frac{|a(r\mid t,\gamma)-a(r\mid t,\nu)|}
{(T-r)^{1-\alpha}}\,dr
\Bigg).
\end{aligned}
\]
No uniformity of $C_{n,k}$ with respect to $n$ is required.
\begin{proposition}
For every $n\in\mathbb N$, the truncated value functional
\[
v_n^\eta(t,x,\gamma)
:=
\inf_{u\in L^\infty([t,T],\overline{B(0,n)})}
J(t,x,\gamma,u)
\]
satisfies property {\rm (L)}.
\end{proposition}

\begin{proof}
The proof is analogous to that of Proposition
\ref{regularity_value_functional}. Indeed, when comparing the values
at two initial conditions, one may use the same admissible control in
both problems. Since $|u|\leq n$, all the estimates are uniform over
the admissible controls defining $v_n^\eta$. In particular,
\[
\left|
X_t^{s,x,\gamma,u}
-
X_t^{s,\widetilde{x},\nu,u}
\right|
\leq
C_{\lambda,b,\eta,T}
\left(
|x-\widetilde{x}|
+
\int_s^t
|a(r\mid s,\gamma)-a(r\mid s,\nu)|\,dr
\right).
\]
Combining this estimate with the corresponding continuity estimate
for the cost functional and taking the infimum over
$\|u\|_\infty\leq n$ proves property {\rm (L)} for $v_n^\eta$.
The constants may depend on $n$.
\end{proof}

Finally, we prove a comparison principle for viscosity solutions of
\eqref{HJB_equation} whose approximating sequences satisfy property
{\rm (L)}. This yields uniqueness within the class introduced above. While the proof is a variant of the proof of Theorem 5.1 in  \cite{gomoyunov2021viscosity}, it requires significant adaptations to account for the additional state variable, which is not present in the  framework analysed in that work.

\begin{theorem}[Comparison principle]
\label{comparison_thm}
Let \(\varphi_1\) and \(\varphi_2\) be continuous functionals for which
there exist sequences
\((\varphi_{1,n})_{n\in\mathbb N}\) and
\((\varphi_{2,n})_{n\in\mathbb N}\)
converging locally uniformly to \(\varphi_1\) and \(\varphi_2\),
respectively. Assume that, for every \(n\in\mathbb N\),
\(\varphi_{1,n}\) is a viscosity subsolution of the \(n\)-truncated
problem, \(\varphi_{2,n}\) is a viscosity supersolution of the
\(n\)-truncated problem, and both \(\varphi_{1,n}\) and
\(\varphi_{2,n}\) satisfy property {\rm (L)}. Assume moreover that
\(
\varphi_{1,n}(T,x,\gamma)
= 
\varphi_{2,n}(T,x,\gamma) = g(x)
\)
for every \(n\in\mathbb N\) and every
\(
(x,\gamma)\in
\mathbb R^e\times AC^\alpha([0,T],\mathbb R^\ell).
\)
Then
\(
\varphi_1\leq\varphi_2
\)
on
\(
[0,T]\times\mathbb R^e
\times AC^\alpha([0,T],\mathbb R^\ell).
\)
\end{theorem}\begin{proof}
Fix $n\in\mathbb N$. For simplicity, throughout the comparison
argument we write
\(
\underline{\varphi}:=\varphi_{1,n},
\overline{\varphi}:=\varphi_{2,n}.
\)
We first prove that
\(
\underline{\varphi}\leq\overline{\varphi}.
\)

Fix $k\in\mathbb N$ and set
\(
K:=\mathcal K_k^n.
\)
Suppose, by contradiction, that
\[
\kappa
:=
\max_{(t,x,\gamma)\in  K}
\bigl(
\underline{\varphi}(t,x,\gamma)
-
\overline{\varphi}(t,x,\gamma)
\bigr)
>0.
\]
Since $K$ is compact, there exists
$(t_0,x_0,\gamma_0)\in K$ such that
\(
\underline{\varphi}(t_0,x_0,\gamma_0)
-
\overline{\varphi}(t_0,x_0,\gamma_0)
=
\kappa.
\)
The goal is to show that for any $k \in \mathbb{N}$ and for the viscosity subsolution \(\underline{\varphi}\) and the
viscosity supersolution \(\overline{\varphi}\) of the problem above, we have  
\begin{equation*}
\underline{\varphi}(t, x,\gamma) \leq \overline{\varphi}(t, x,\gamma), \quad \forall (t, x, \gamma) \in K.
\end{equation*}

Define the auxiliary functional $\Phi_\epsilon: [0,T] \times  AC^\alpha([0,T], \mathbb{R}^\ell) \times \mathbb{R}^e  \times [0,T]  \times AC^\alpha([0,T], \mathbb{R}^\ell) \times \mathbb{R}^e \to \mathbb{R}$ by
\begin{align*}
\Phi_\epsilon(t, \gamma, x, \tau, \nu, y) 
&:= \underline{\varphi}(t, x, \gamma) - \overline{\varphi}(\tau, y, \nu) - (2T - t - \tau)\tilde{\kappa} 
- \frac{(t-\tau)^2}{\epsilon^{3/\alpha}} \\
&\quad - \frac{\varpi_\epsilon(t, \gamma, \tau, \nu)}{\epsilon} 
- \frac{((x-y)^2 + \epsilon^{\frac{2}{c}(3/\alpha + 2)})^{c/2}}{c\epsilon},
\end{align*}
where $\tilde{\kappa} = \kappa/(4T)$ and $\varpi_\epsilon$ is defined in \eqref{auxiliary_functional}.  
The functional \(\Phi_\epsilon\) is the doubled-variable penalisation used in
the comparison argument. It compares the subsolution \(\underline{\varphi}\) and the
supersolution \(\overline{\varphi}\) at two possibly distinct points, while forcing
\((t,\gamma,x)\) and \((\tau,\nu,y)\) to approach one another as
\(\epsilon\to 0\). The separate penalisation terms control, respectively,
the time variables, the fractional histories, and the finite-dimensional state
variables. The regularising powers of \(\epsilon\) smooth the functional near
the diagonal and balance the fractional estimates, while
\((2T-t-\tau)\tilde\kappa\) provides the strictness needed for the final
contradiction.

We maximise $\Phi_\epsilon$ on the compact set
\(
K \times K
\). By continuity of
\(\varpi_\epsilon\) (Lemma~\ref{lemma_5_4_g21}) and of the remaining terms in
\(\Phi_\epsilon\), there exists a maximiser
\[
(t_\epsilon,\gamma^\epsilon,x_\epsilon,
 \tau_\epsilon,\nu^\epsilon,y_\epsilon)
\]
such that
\begin{equation*}
\Phi_\epsilon(t_\epsilon,\gamma^\epsilon,x_\epsilon,
\tau_\epsilon,\nu^\epsilon,y_\epsilon)
:=
\max_{(t,x, \gamma),(\tau,y, \nu)\in K}
\Phi_\epsilon(t,\gamma,x,\tau,\nu,y).
\end{equation*}
The rest of the argument consists in extracting the consequences of this
maximality as \(\epsilon\to 0\). We first show that the two time variables
coalesce in the limit.
Recalling that since $\underline{\varphi}(T, x_\epsilon, \gamma^{\epsilon}) = \overline{\varphi}(T, x_\epsilon, \gamma^{\epsilon}) = g(x_{\epsilon})$, then
\begin{align*}
    0 &\leq \Phi_\epsilon(t_\epsilon, \gamma^\epsilon, x_\epsilon, \tau_\epsilon, \nu^\epsilon, y_\epsilon) - \Phi_\epsilon(T, \gamma^\epsilon, x_\epsilon, T, \gamma^\epsilon, x_\epsilon)\\
    &= \underline{\varphi}(t_\epsilon, x_\epsilon, \gamma^\epsilon) - \overline{\varphi}(\tau_\epsilon, y_\epsilon, \nu^\epsilon) -(2T- t_\epsilon- \tau_\epsilon)\tilde \kappa - \frac{(t_\epsilon- \tau_\epsilon)^2}{\epsilon^{3/\alpha}} - \varpi_\epsilon(t_\epsilon, \gamma^\epsilon, \tau_\epsilon, \nu^\epsilon) \\
     & \quad - \frac{((x_\epsilon-y_\epsilon)^2 + \epsilon^{\frac{2}{c}(3/\alpha + 2)})^{c/2}}{c\epsilon} + \frac{\epsilon^{3/\alpha + 1}}{c}, 
\end{align*}
where we used Lemma \ref{lemma_5_4_g21} in the first identity to replace $\varpi_\epsilon(T, \gamma^\epsilon, x_\epsilon, T, \gamma^\epsilon, x_\epsilon)$ with $0$.\newline
Now we can use the fact that $\varpi_\epsilon$  is non negative (Lemma \ref{lemma_5_4_g21}) and by defining
\[
\kappa_1 := \max_{(t,x,\gamma),(\tau,y,\nu) \in  K} (\underline{\varphi}(t,x,\gamma) - \overline{\varphi}(\tau, y, \nu)),
\]
it follows that
\begin{align*}
0 &\leq \kappa_1  - \frac{(t_\epsilon- \tau_\epsilon)^2}{\epsilon^{3/\alpha}}   + \frac{\epsilon^{3/\alpha + 1}}{c}. 
\end{align*}

From the above, one obtains
\begin{equation}\label{convergence_t_eps_tau_eps}
|t_\epsilon-\tau_\epsilon|
\leq C\epsilon^{3/(2\alpha)}.
\end{equation}
We next show that the state variables also coalesce as \(\epsilon \to 0\). Assume without loss of generality that $\tau_\epsilon \ge t_\epsilon$ and $\epsilon \in (0,1)$. 
As we recalled above for every $k > 0$ if $\gamma \in AC_k$, then $a(\cdot \mid t,\gamma) \in AC_k$ for every $t \in [0, T]$.\newline
From this it is possible to define 
$\Phi_\epsilon(t_\epsilon, \gamma^\epsilon, x_\epsilon, \tau_\epsilon, a(\cdot | t_\epsilon, \gamma^\epsilon), x_\epsilon)$. Moreover, by the definition of $(t_\epsilon, \gamma^\epsilon, x_\epsilon, \tau_\epsilon, \nu^\epsilon, y_\epsilon)$, we have 
\[\Phi_\epsilon(t_\epsilon, \gamma^\epsilon, x_\epsilon, \tau_\epsilon, \nu^\epsilon, y_\epsilon) \ge \Phi_\epsilon(t_\epsilon, \gamma^\epsilon, x_\epsilon, \tau_\epsilon, a(\cdot | t_\epsilon, \gamma^\epsilon), x_\epsilon),\]
so that
\begin{align*}
    & - \overline{\varphi}(\tau_\epsilon, y_\epsilon, \nu^\epsilon) - \varpi_\epsilon(t_\epsilon, \gamma^\epsilon, \tau_\epsilon, \nu^\epsilon) - \frac{((x_\epsilon-y_\epsilon)^2 + \epsilon^{\frac{2}{c}(3/\alpha + 2)})^{c/2}}{c\epsilon} 
     \\
     &- \left( - \overline{\varphi}(\tau_\epsilon, x_\epsilon, a(\cdot|t_\epsilon, \gamma^\epsilon)) - \varpi_\epsilon(t_\epsilon, \gamma^\epsilon, x_\epsilon, \tau_\epsilon, a(\cdot|t_\epsilon, \gamma^\epsilon), x_\epsilon) - \frac{\epsilon^{(3/\alpha + 1)}}{c}\right) \geq 0.
\end{align*}

By using  the property (L) and subsequently Lemma \ref{lemma_5_5_g21} to find (for a possibly changing constant $C > 0$) the estimate 
\begin{align*}
\overline{\varphi}(\tau_\epsilon, x_\epsilon, a(\cdot|t_\epsilon, \gamma^\epsilon))- \overline{\varphi}(\tau_\epsilon, y_\epsilon, \nu^\epsilon) &\leq C_k\left( |x_\epsilon - y_\epsilon| + \varpi_{\epsilon}(t_\epsilon, \gamma^\epsilon, \tau_\epsilon, \nu^\epsilon) + C_1\epsilon^{\frac{c}{c-1}} \right)^{\frac{1}{c}}\\
&\leq C\Big(\varpi_\epsilon(t_\epsilon, \gamma^\epsilon, \tau_\epsilon, \nu^\epsilon) + C_1 \epsilon^{c/(c-1)} + \frac{((x_\epsilon - y_\epsilon)^2 + \epsilon^2)^{c/2}}{c}\Big)^{1/c}.
\end{align*}
Using Lemma \ref{lemma_5_4_g21} to write
\[
\varpi_\epsilon(t_\epsilon, \gamma^\epsilon, x_\epsilon, \tau_\epsilon, a(\cdot|t_\epsilon, \gamma^\epsilon), x_\epsilon) = \varpi_\epsilon(t_\epsilon, \gamma^\epsilon, x_\epsilon, t_\epsilon,  \gamma^\epsilon, x_\epsilon) = 0,
\]
we obtain, after basic algebraic manipulations and for $C > 0$, the inequality
\begin{align*}
&\frac{((x_\epsilon - y_\epsilon)^2 + \epsilon^{\frac{2}{c}(3/\alpha + 2)})^{c/2}}{c\epsilon} + \frac{\varpi_\epsilon(t_\epsilon, \gamma^\epsilon, \tau_\epsilon, \nu^\epsilon)}{\epsilon} \\
&\quad \le C\Big(\varpi_\epsilon(t_\epsilon, \gamma^\epsilon, \tau_\epsilon, \nu^\epsilon) + C_1 \epsilon^{c/(c-1)} + \frac{((x_\epsilon - y_\epsilon)^2 + \epsilon^2)^{c/2}}{c}\Big)^{1/c}.
\end{align*}
Since all terms on the left are positive, $\epsilon \in (0, 1)$ and $1< c < 2$, this implies
\begin{equation}\label{convergence_x_eps_y_eps}
|x_\epsilon - y_\epsilon| \le C_3 \epsilon^{2/c}, \quad \varpi_\epsilon(t_\epsilon, \gamma^\epsilon, \tau_\epsilon, \nu^\epsilon) \le C_3 \epsilon^{2}.
\end{equation}

By Lemma \ref{lemma_5_6_g21}, equiboundedness and equicontinuity of $AC^\alpha_k$ functions, one also has by Arzelà–Ascoli theorem
\begin{equation}\label{ascoli_arzela_a}
\|a(\cdot| t_\epsilon, \gamma^\epsilon) - a(\cdot| \tau_\epsilon, \nu^\epsilon)\|_\infty \to 0, \quad
\|a(t_\epsilon| t_\epsilon, \gamma^\epsilon) - a(\tau_\epsilon| \tau_\epsilon, \nu^\epsilon)\|_\infty \to 0, \quad \text{as } \epsilon \to 0^+.
\end{equation}
Consequently, since by definition  
\(
\gamma^\epsilon_{t_\epsilon} = a(t_\epsilon \mid t_\epsilon, \gamma^\epsilon), \nu^\epsilon_{\tau_\epsilon} = a(\tau_\epsilon \mid \tau_\epsilon, \nu^\epsilon), 
\)
then 
\begin{equation}\label{convergence_gamma_eps_nu_eps}
\|\gamma^\epsilon(t_\epsilon) - \nu^\epsilon(\tau_\epsilon)\| \to 0, \quad \text{as } \epsilon \to 0^+.
\end{equation}
The next step in the proof consists in showing that there exist
$0<z<T$ and $\epsilon^*\in(0,1)$ such that, for every
$\epsilon\in(0,\epsilon^*]$,
\(
    t_\epsilon,\tau_\epsilon\in[0,T-z].
\)
This ensures that the ci-derivatives of the penalisation functionals
used below are well defined at the maximising points.

By the maximality of
$(t_\epsilon,\gamma^\epsilon,x_\epsilon,
\tau_\epsilon,\nu^\epsilon,y_\epsilon)$, we have
\begin{align*}
    &\Phi_\epsilon
    (t_\epsilon,\gamma^\epsilon,x_\epsilon,
    \tau_\epsilon,\nu^\epsilon,y_\epsilon)
    \geq
    \Phi_\epsilon
    (t_0,\gamma_0,x_0,t_0,\gamma_0,x_0)=
    \kappa
    -2(T-t_0)\tilde\kappa
    -\frac{\epsilon^{3/\alpha+1}}{c}.
\end{align*}
Since $\tilde\kappa=\kappa/(4T)$,
\[
    2(T-t_0)\tilde\kappa
    \leq
    2T\tilde\kappa
    =
    \frac{\kappa}{2}.
\]
Hence, after reducing $\epsilon^*>0$ if necessary so that
\[
    \frac{\epsilon^{3/\alpha+1}}{c}
    \leq
    \frac{\kappa}{8},
    \qquad
    \epsilon\in(0,\epsilon^*],
\]
we obtain
\[
    \Phi_\epsilon
    (t_\epsilon,\gamma^\epsilon,x_\epsilon,
    \tau_\epsilon,\nu^\epsilon,y_\epsilon)
    \geq
    \frac{3\kappa}{8}.
\]
Since all the penalisation terms appearing with a negative sign in the
definition of $\Phi_\epsilon$ are non-negative, it follows that
\begin{equation}\label{away_from_terminal_lower_bound}
    \underline{\varphi}
    (t_\epsilon,x_\epsilon,\gamma^\epsilon)
    -
    \overline{\varphi}
    (\tau_\epsilon,y_\epsilon,\nu^\epsilon)
    \geq
    \frac{3\kappa}{8}.
\end{equation}

Suppose, by contradiction, that no such $z>0$ exists. Then there is a
sequence $\epsilon_j\to 0$ such that
\(
    \max\{t_{\epsilon_j},\tau_{\epsilon_j}\}\to T.
\)
Since
\(
    |t_\epsilon-\tau_\epsilon|\to0,
\)
it follows that
\(
    t_{\epsilon_j}\to T,
    \,
    \tau_{\epsilon_j}\to T.
\)
By the uniform continuity of $\underline{\varphi}$ and
$\overline{\varphi}$ on $K$, together with the terminal condition
\(
    \underline{\varphi}(T,x,\gamma)
    =
    \overline{\varphi}(T,x,\gamma)
    =
    g(x),
\)
we obtain
\(
    \left|
        \underline{\varphi}
        (t_{\epsilon_j},x_{\epsilon_j},\gamma^{\epsilon_j})
        -
        g(x_{\epsilon_j})
    \right|
    \to0
\)
and
\(
    \left|
        \overline{\varphi}
        (\tau_{\epsilon_j},y_{\epsilon_j},\nu^{\epsilon_j})
        -
        g(y_{\epsilon_j})
    \right|
    \to0.
\)
Moreover, since
\(
    |x_{\epsilon_j}-y_{\epsilon_j}|\to0
\)
and $g$ is continuous,
\(
    |g(x_{\epsilon_j})-g(y_{\epsilon_j})|
    \to0.
\)
Consequently,
\(
    \underline{\varphi}
    (t_{\epsilon_j},x_{\epsilon_j},\gamma^{\epsilon_j})
    -
    \overline{\varphi}
    (\tau_{\epsilon_j},y_{\epsilon_j},\nu^{\epsilon_j})
    \to0,
\)
which contradicts \eqref{away_from_terminal_lower_bound}.

Therefore, there exist $0<z<T$ and $\epsilon^*\in(0,1)$ such that
\(
    t_\epsilon,\tau_\epsilon\in[0,T-z]
\)
for every $\epsilon\in(0,\epsilon^*]$.
We now enter the comparison step. The estimates above show that, along a
subsequence if necessary, the maximising points of \(\Phi_\epsilon\) concentrate
near the diagonal as \(\epsilon\to 0\). This allows us to use the
penalisation terms in \(\Phi_\epsilon\) as test functions for
\(\underline{\varphi}\) and \(\overline{\varphi}\) at the points
\[
    (t_\epsilon,\gamma^\epsilon,x_\epsilon),
    \qquad
    (\tau_\epsilon,\nu^\epsilon,y_\epsilon).
\]
Applying the viscosity subsolution inequality to \(\underline{\varphi}\) and the
viscosity supersolution inequality to \(\overline{\varphi}\), and then subtracting the
two inequalities, we obtain the key estimate leading to a contradiction.
For every $\epsilon \in (0, \epsilon^*)$, define
\[
\psi_1(t,x,\gamma) := \overline{\varphi}(\tau_\epsilon, y_\epsilon, \nu^\epsilon) + (2T - t - \tau_\epsilon)\tilde{\kappa} + \frac{(t - \tau_\epsilon)^2}{\epsilon^{3/\alpha}} + \frac{\mu_\epsilon^{(\tau_\epsilon, \nu^\epsilon)}(t,\gamma)}{\epsilon} + \frac{((x - y_\epsilon)^2 + \epsilon^{\frac{2}{c}(3/\alpha + 2)})^{c/2}}{c\epsilon},
\]
with $\mu_\epsilon^{(\tau,\nu)}(t, \gamma) := \varpi_\epsilon(t, \gamma, \tau,\nu )$. Since $\mu_\epsilon^{(\tau,\nu)}(t, \gamma)$ is ci-smooth of order $\alpha$ (Lemma \ref{lemma_5_7_g21}), then it follows that $\psi_1$ is ci-smooth of order $\alpha$ with
\[
\frac{\partial^\alpha}{\partial t} \psi_1 = -\tilde{\kappa} + 2\frac{t - \tau_\epsilon}{\epsilon^{3/\alpha}}, \quad
\nabla^\alpha_x \psi_1 = \frac{(x - y_\epsilon)((x - y_\epsilon)^2 + \epsilon^{\frac{2}{c}(3/\alpha + 2)})^{c/2 - 1}}{\epsilon}, \quad
\nabla^\alpha_\gamma \psi_1 = \frac{\nabla^\alpha_\gamma \mu_\epsilon^{(\tau_\epsilon, \nu^\epsilon)}(t,\gamma)}{\epsilon}.
\]
Notice that for every $(t,x,\gamma) \in  K$ we have
\begin{equation}\label{phi_psi_is_sub}
\begin{aligned}
    \underline{\varphi}(t,x,\gamma) - \psi_1(t,x,\gamma) &= \Phi_\epsilon(t, \gamma, x, \tau_\epsilon, \nu^\epsilon, y_\epsilon)\\
    &\leq  \Phi_\epsilon(t_\epsilon, \gamma^\epsilon, x_\epsilon, \tau_\epsilon, \nu^\epsilon, y_\epsilon)\\
    &= \underline{\varphi}(t_\epsilon, \gamma^\epsilon, x_\epsilon) - \psi_1(t_\epsilon,  x_\epsilon, \gamma^\epsilon),
\end{aligned}
\end{equation}
so that, by the definition of viscosity subsolution and the fact that $t_\epsilon < T$, one obtains
\begin{align}
&\tilde{\kappa} - 2\frac{t_\epsilon - \tau_\epsilon}{\epsilon^{3/\alpha}} 
- \langle \nabla^\alpha_x \psi_1, b(x_\epsilon, \gamma^\epsilon_{t_\epsilon}) + \lambda(x_\epsilon, \gamma^\epsilon_{t_\epsilon}) \dot{\eta}_{t_\epsilon} \rangle
+ H_n(x_\epsilon, \gamma^\epsilon_{t_\epsilon}, \nabla^\alpha_\gamma \psi_1) \le 0. \label{HJB_1_uniqueness}
\end{align}

Similarly, define
\[
\psi_2(\tau,y,\nu) := \underline{\varphi}(t_\epsilon, x_\epsilon, \gamma^\epsilon) - (2T - t_\epsilon - \tau)\tilde{\kappa} - \frac{(t_\epsilon - \tau)^2}{\epsilon^{3/\alpha}} - \frac{\mu_\epsilon^{(t_\epsilon, \gamma^\epsilon)}(\tau,\nu)}{\epsilon} - \frac{((x_\epsilon - y)^2 + \epsilon^{\frac{2}{c}(3/\alpha + 2)})^{c/2}}{c\epsilon},
\] 
which, thanks to Lemma \ref{lemma_5_7_g21} is ci-smooth  of order $\alpha$ with
\[
\frac{\partial^\alpha}{\partial t} \psi_2 = \tilde{\kappa} + 2\frac{t_\epsilon - \tau}{\epsilon^{3/\alpha}}, \quad
\nabla^\alpha_x \psi_2 = -\frac{(x_\epsilon - y)((x_\epsilon - y)^2 + \epsilon^{\frac{2}{c}(3/\alpha + 2)})^{c/2 - 1}}{\epsilon}, \quad
\nabla^\alpha_\gamma \psi_2 = -\frac{\nabla^\alpha_\gamma \mu_\epsilon^{(t_\epsilon, \gamma^\epsilon)}(\tau,\nu)}{\epsilon}.
\]
Moreover, with an argument similar to the one shown in \eqref{phi_psi_is_sub} we have that $\overline{\varphi} - \psi_2$ attains a local minimum at $(\tau_\epsilon, y_\epsilon, \nu_\epsilon)$ so that, thanks to the definition of viscosity supersolution and $\tau_\epsilon < T$, we have
\begin{align}
&-\tilde{\kappa} - 2\frac{t_\epsilon - \tau_\epsilon}{\epsilon^{3/\alpha}} 
- \langle \nabla^\alpha_x \psi_2, b(y_\epsilon, \nu^\epsilon_{\tau_\epsilon}) + \lambda(y_\epsilon, \nu^\epsilon_{\tau_\epsilon}) \dot{\eta}_{\tau_\epsilon} \rangle
+ H_n(y_\epsilon, \nu^\epsilon_{\tau_\epsilon}, \nabla^\alpha_\gamma \psi_2) \ge 0. \label{HJB_2_uniqueness}
\end{align}

Combining \eqref{HJB_1_uniqueness} and \eqref{HJB_2_uniqueness}, we obtain

\begin{align}
2\widetilde{\kappa}
\leq{}&
\left\langle
\nabla_x^\alpha\psi_1,
b(x_\epsilon,\gamma^\epsilon_{t_\epsilon})
+\lambda(x_\epsilon,\gamma^\epsilon_{t_\epsilon})
\dot\eta_{t_\epsilon}
-b(y_\epsilon,\nu^\epsilon_{\tau_\epsilon})
-\lambda(y_\epsilon,\nu^\epsilon_{\tau_\epsilon})
\dot\eta_{\tau_\epsilon}
\right\rangle
\notag\\
&+
H_n\left(
y_\epsilon,\nu^\epsilon_{\tau_\epsilon},
\nabla_\gamma^\alpha\psi_2
\right)
-
H_n\left(
x_\epsilon,\gamma^\epsilon_{t_\epsilon},
\nabla_\gamma^\alpha\psi_1
\right),
\label{final_inequality_uniqueness}
\end{align}
Since \(1<c<2\), via \eqref{convergence_x_eps_y_eps} we get for some $C>0$
\begin{align*}
\left|\nabla_x^\alpha\psi_1\right|
&=
\frac{
|x_\epsilon - y_\epsilon|
\left(
|x_\epsilon - y_\epsilon|^2+
\epsilon^{\frac{2}{c}(3/\alpha+2)}
\right)^{(c-2)/2}
}{\epsilon}\\
&
\leq
C\epsilon^{\frac{2(c-1)}{c}-1}
=
C\epsilon^{\alpha-1}.
\end{align*}
Therefore, using the Lipschitz continuity of \(b\) and
\(\lambda\), the smoothness of $\eta$, Lemma~\ref{lemma_5_5_g21}, and
\(2/c=2-\alpha\), we obtain for some $C> 0$ independent of $\epsilon$
\begin{align*}
 &\left|\langle \nabla^\alpha_x \psi_1, b(x_\epsilon, \gamma^\epsilon_{t_\epsilon}) - b(y_\epsilon, \nu^\epsilon_{\tau_\epsilon}) + \lambda(x_\epsilon, \gamma^\epsilon_{t_\epsilon})\dot{\eta}_{t_\epsilon} - \lambda(y_\epsilon, \nu^\epsilon_{\tau_\epsilon})\dot{\eta}_{\tau_\epsilon} \rangle \right| \\
&\leq C \left|\nabla_x^\alpha\psi_1\right|
\left(
|x_\epsilon-y_\epsilon|
+
|\gamma^\epsilon_{t_\epsilon}
-\nu^\epsilon_{\tau_\epsilon}| + |t_\epsilon - \tau_\epsilon|
\right)\\
&\qquad\leq
C\epsilon^{\alpha-1}
\left(
\epsilon^{2/c}
+
\epsilon^{\frac{2\alpha}{1+\alpha c}} 
\right)\\
&\qquad=
C\epsilon
+
C\epsilon^{
\alpha-1+\frac{2\alpha}{1+\alpha c}
}.
\end{align*}
Now
\[
    \alpha-1+\frac{2\alpha}{1+\alpha c}
    =
    \frac{-\alpha^2+5\alpha-2}{2+\alpha}.
\]
If
\(
    \alpha > \lfloor p\rfloor/p
\)
for some non-integer \(p \in [2, \infty) \setminus \mathbb N\), then \(\alpha>1/2\). Since
\[
    -\left(\frac12\right)^2
    +5\left(\frac12\right)-2
    =\frac14>0
\]
and \(5-2\alpha>0\) for \(\alpha\in(0,1)\), we have
\(
    -\alpha^2+5\alpha-2>0.
\)
Consequently,
\begin{equation}\label{lambda_b_to_zero}
\left|\langle \nabla^\alpha_x \psi_1, b(x_\epsilon, \gamma^\epsilon_{t_\epsilon}) - b(y_\epsilon, \nu^\epsilon_{\tau_\epsilon}) + \lambda(x_\epsilon, \gamma^\epsilon_{t_\epsilon})\dot{\eta}_{t_\epsilon} - \lambda(y_\epsilon, \nu^\epsilon_{\tau_\epsilon})\dot{\eta}_{\tau_\epsilon} \rangle \right|
\to 0.
\end{equation}

Since the control set defining $H_n$ is compact, the local Lipschitz
continuity of $f$ gives, on the compact under consideration,
\[
\begin{aligned}
|H_n(x,y,\phi)-H_n(\widetilde x,\widetilde y,\widetilde \phi)|
&\leq
\sup_{|u|\leq n}
\Big(
|\langle \phi-\widetilde \phi,u\rangle|
+
|f(x,y,u)-f(\widetilde x,\widetilde y,u)|
\Big)
\\
&\leq
n|\phi-\widetilde \phi|
+
C_{n,k}
\bigl(
|x-\widetilde x|
+
|y-\widetilde y|
\bigr).
\end{aligned}
\]
Therefore, Lemma \ref{lemma_5_8_g21}, together with
\eqref{convergence_t_eps_tau_eps},
\eqref{convergence_x_eps_y_eps}, and
\eqref{ascoli_arzela_a}, implies that
\[
\begin{aligned}
&H_n\left(
y_\epsilon,\nu^\epsilon_{\tau_\epsilon},
\nabla_\gamma^\alpha\psi_2
\right)
-
H_n\left(
x_\epsilon,\gamma^\epsilon_{t_\epsilon},
\nabla_\gamma^\alpha\psi_1
\right)
\to 0.
\end{aligned}
\] 
Letting $\epsilon\to0$ in
\eqref{final_inequality_uniqueness} gives
\(
2\widetilde{\kappa}\leq0,
\)
which is a contradiction. Therefore,
\(
\varphi_{1,n}\leq\varphi_{2,n}
\)
on $\mathcal K_k^n$. Since $k$ is arbitrary and
$(\mathcal K_k^n)_{k\in\mathbb N}$ exhausts the time-state space,
\(
\varphi_{1,n}\leq\varphi_{2,n}
\)
everywhere.

Finally, letting $n\to\infty$ and using the locally uniform
convergences
\(
\varphi_{1,n}\to \varphi_1,
\varphi_{2,n}\to \varphi_2,
\)
we obtain
\(
\varphi_1\leq\varphi_2.
\)
\end{proof}

By Proposition~\ref{v_viscosity_sol}, the value functional $v$ is a
viscosity solution of \eqref{HJB_equation} in the sense of
Definition~\ref{definition_viscosity_solution}, with approximating
sequence $(v_n)_{n\in\mathbb N}$. Moreover, each $v_n$ satisfies
property {\rm (L)}.

Let $w$ be another viscosity solution admitting an approximating
sequence $(w_n)_{n\in\mathbb N}$ such that each $w_n$ satisfies
property {\rm (L)}. Since
\(
v_n(T,x,\gamma)=g(x)=w_n(T,x,\gamma)
\)
for every $n\in\mathbb N$ and every
\(
(x,\gamma)\in
\mathbb R^e\times AC^\alpha([0,T],\mathbb R^\ell),
\)
Theorem~\ref{comparison_thm}, applied with
$\varphi_{1,n}=v_n$ and $\varphi_{2,n}=w_n$, gives
\(
v\leq w.
\)
Applying the same theorem with the roles of $v_n$ and $w_n$
interchanged gives
\(
w\leq v.
\)
Consequently,
\(
v=w
\)
on
\(
[0,T]\times\mathbb R^e
\times AC^\alpha([0,T],\mathbb R^\ell).
\)
We record this conclusion in the following corollary.

\begin{corollary}\label{uniqueness_thm_viscosity}
Let
\[
H(x,y,\phi)
:=
\sup_{u\in\mathbb R^\ell}
\bigl\{
-\langle\phi,u\rangle-f(x,y,u)
\bigr\}.
\]
The value functional $v$ is the unique viscosity solution, in the
sense of Definition \ref{definition_viscosity_solution}, of
\[
\begin{cases}
-\dfrac{\partial^\alpha}{\partial t}v(t,x,\gamma)
-\left\langle
\nabla_x^\alpha v(t,x,\gamma),
b(x,\gamma_t)+\lambda(x,\gamma_t)\dot\eta_t
\right\rangle
+H\bigl(x,\gamma_t,\nabla_\gamma^\alpha v(t,x,\gamma)\bigr)
=0,
\\
\hfill
(t,x,\gamma)\in
[0,T)\times\mathbb R^e
\times AC^\alpha([0,T],\mathbb R^\ell),
\\[1mm]
v(T,x,\gamma)=g(x),
\qquad
(x,\gamma)\in
\mathbb R^e\times AC^\alpha([0,T],\mathbb R^\ell),
\end{cases}
\]
among viscosity solutions whose approximating sequences satisfy
property {\rm (L)}.
\end{corollary}

This uniqueness result, together with the local uniform continuity of
the value functional with respect to the driving path, yields the main
result of this section.

\begin{theorem}\label{final_theorem_viscosity}
The value functional \eqref{extended_value_functional} is a rough
viscosity solution of \eqref{rough_HJB_equation} in the sense of
Definition~\ref{rough_viscosity_solution}. Moreover, it is unique
within the class of rough viscosity solutions whose smooth
approximating viscosity solutions admit truncated approximating
sequences satisfying property \emph{(L)}.
\end{theorem}
\begin{example}
We now go back to the example presented in the introduction and show that the assumptions presented in the paper allow us to recover well-posedness. Recall that in the original formulation the trading rate is bounded by \(M>0\). Fix \(\overline M>M\) and choose a function $\lambda \in C_b^\infty(\mathbb R)$ with
\[
    \lambda(y) = y, \quad  |y| \leq \overline M.
\]
The wealth process then satisfies
\begin{equation*}
    dX^{0,x,\gamma,u}_t
    =
    \sigma\lambda(\gamma_t^{0,u})\,d\mathbf{P}_t,
    \qquad
    X_0=x\in\mathbb{R},
\end{equation*}
where, for any non-integer $p>1/(1/2-\rho)$, the path $P$ has finite
$p$-variation. Since $P$ is real-valued, it admits a canonical geometric
$p$-rough path lift $\mathbf P$; see \cite{friz2010differential}. Fix
\[
    \frac{\lfloor p\rfloor}{p}<\alpha<1,
    \qquad
    1<\kappa<
    \frac{1}{1-\alpha+\lfloor p\rfloor/p},
\]
and let $q\in\mathbb N$ be such that
\[
    2q>
    \lfloor p\rfloor p
    \vee
    \frac{\kappa}{\kappa-1}.
\]
We assume that $\gamma\in AC^\alpha([0,T],\mathbb R)$.

Introducing the fractional control \(u\), we obtain
\begin{equation}\label{example_dynamics}
\begin{aligned}
&dX^{0,x,\gamma,u}_t
=
\sigma\lambda(\gamma^{0,u}_t)\,d\mathbf{P}_t,
\qquad
X_0=x\in\mathbb{R},
\\
&D^\alpha_{0^+}(\gamma^{0,u})(t)=u_t.
\end{aligned}
\end{equation}
The cost functional for our problem is
\begin{equation}\label{example_cost_functional}
J(t,x,\gamma,u)
=
X^{t,x,\gamma,u}_T
-
c\int_t^T |u_r|^{2q}\,dr.
\end{equation}
Here \(q\in\mathbb N\) is chosen as above. When considering the
explicit representation below, we retain the bound from the
motivating example and restrict attention to initial histories
satisfying
\(
    \|\gamma\|_{\infty;[0,t]}\leq M.
\)
The continuation then satisfies
\(    \|a(\cdot\mid t,\gamma)\|_{\infty;[t,T]}\leq M.
\)
By the coercivity argument used in the proof of
Proposition~\ref{boundedness_val_fun}, together with the fractional
reconstruction estimate, and keeping the dependence on \(c\)
explicit, the pseudo-controls relevant to the optimization can be
made sufficiently small by choosing \(c>0\) sufficiently large.
Since \(\overline M-M>0\), we may therefore choose \(c\) so that the
optimization may be restricted, without changing its value, to
pseudo-controls for which
\[
    \sup_{s\in[t,T]}
    \left|
        \gamma_s^{t,\gamma,u}
        -
        a(s\mid t,\gamma)
    \right|
    \leq
    \overline M-M.
\]
Consequently,
\[
    \|\gamma^{t,\gamma,u}\|_{\infty;[t,T]}
    \leq \overline M,
\]
and hence, along all control trajectories relevant to the optimization,
\[
    \lambda(\gamma_s^{t,\gamma,u})
    =
    \gamma_s^{t,\gamma,u},
    \qquad s\in[t,T].
\]

For every \(n\in\mathbb{N}\), define the truncated value functional
\[
v_n^\eta(t,x,\gamma)
:=
\sup_{\substack{u\in L^\infty([t,T])\\
\|u\|_\infty\leq n}}
J(t,x,\gamma,u).
\]
The corresponding \(n\)-truncated HJB equation is obtained from the
equation below by replacing the infimum over \(u\in\mathbb{R}\) with
the infimum over \(|u|\leq n\). By Proposition
\ref{v_viscosity_sol}, \(v_n^\eta\) is a viscosity solution of this
truncated equation and
\(
v_n^\eta\to v^\eta
\)
locally uniformly. The HJB equation corresponding to the limiting
problem is therefore
\[
\begin{cases}
\begin{aligned}
-\frac{\partial^\alpha}{\partial t}v(t,x,\gamma)
&-
\sigma\lambda(\gamma_t)
\nabla^\alpha_xv(t,x,\gamma)\dot\eta_t
\\
&\quad
+\inf_{u\in\mathbb{R}}
\left(
-u\nabla^\alpha_\gamma v(t,x,\gamma)
+c|u|^{2q}
\right)
=0,
\end{aligned}
&
[0,T)\times\mathbb{R}\times
AC^\alpha([0,T],\mathbb{R}),
\\
v(T,x,\gamma)=x,
&
\mathbb{R}\times AC^\alpha([0,T],\mathbb{R}).
\end{cases}
\]

For initial histories satisfying
\(
    \|\gamma\|_{\infty;[0,t]}\leq M,
\)
the truncation is inactive along the optimizing trajectories and the
limiting value functional can be computed explicitly. It is given by
\begin{align*}
v^\eta(t,x,\gamma)
&=
x
+\sigma\int_t^T a(r\mid t,\gamma)\,d\eta_r
+c_1\int_t^T
\left|I^\alpha_{T^-}(\dot\eta)(r)\right|^{\frac{2q}{2q-1}}\,dr,
\end{align*}
where
\[
c_1
=
\frac{2q-1}{2q}
\left(\frac{1}{2qc}\right)^{\frac{1}{2q-1}}
|\sigma|^{\frac{2q}{2q-1}},
\]
and the right-sided integral \(I^\alpha_{T^-}\) is defined as in
Appendix \ref{fundamentals_fractional_integrations}. Indeed, the
constant \(c_1\) follows from
\[
\sup_{u\in\mathbb{R}}
\left\{
\sigma zu-c|u|^{2q}
\right\}
=
c_1|z|^{\frac{2q}{2q-1}}.
\]
Consequently, \(v^\eta\) is a viscosity solution of the preceding HJB
equation in the sense of Definition
\ref{definition_viscosity_solution}.

For the term
\[
A(t,\gamma)
=
\int_t^T a(r\mid t,\gamma)\,d\eta_r,
\]
we have
\[
\frac{\partial^\alpha}{\partial t}A(t,\gamma)
=
-\lambda(\gamma_t)\dot\eta_t,
\qquad
\nabla^\alpha_\gamma A(t,\gamma)
=
I^\alpha_{T^-}(\dot\eta)(t),
\qquad
\nabla^\alpha_xA(t,\gamma)=0.
\]
In fact, for any path \(\nu\in AC^\alpha\) such that
\(\nu(s)=\gamma(s)\) for every \(s\in[0,t]\), we have
\begin{align}
&A(z,\nu)-A(t,\gamma)
\nonumber\\
&=
\int_z^T a(r\mid z,\nu)\,d\eta_r
-
\int_t^T a(r\mid t,\gamma)\,d\eta_r
\nonumber\\
&=
\int_z^T
\bigl(
a(r\mid z,\nu)-a(r\mid t,\gamma)
\bigr)\,d\eta_r
-
\int_t^z a(r\mid t,\gamma)\,d\eta_r
\nonumber\\
&=
\frac{1}{\Gamma(\alpha)}
\int_z^T\int_t^z
\frac{D^\alpha_{0^+}\nu(\xi)}
{(r-\xi)^{1-\alpha}}
\,d\xi\,d\eta_r
-
\int_t^z a(r\mid t,\gamma)\,d\eta_r.
\label{example_last_line_ci_derivative}
\end{align}
Here we used the definition of the functional \(a\), together with
the fact that \(\nu=\gamma\) on \([0,t]\).

For the first term in
\eqref{example_last_line_ci_derivative}, we have
\begin{align*}
&\frac{1}{\Gamma(\alpha)}
\int_z^T\int_t^z
\frac{D^\alpha_{0^+}\nu(\xi)}
{(r-\xi)^{1-\alpha}}
\,d\xi\,d\eta_r
\\
&\qquad
=
I^\alpha_{T^-}(\dot\eta)(t)
\int_t^zD^\alpha_{0^+}\nu(\xi)\,d\xi
+o(z-t).
\end{align*}
Indeed, for every \(\xi\in[t,z]\),
\[
\left|
\frac{1}{\Gamma(\alpha)}
\int_z^T
\frac{\dot\eta_r}{(r-\xi)^{1-\alpha}}\,dr
-
I^\alpha_{T^-}(\dot\eta)(t)
\right|
\leq
C_{\eta,\alpha}|z-t|^\alpha.
\]
Consequently, by Fubini's theorem,
\begin{align*}
&\left|
\frac{1}{\Gamma(\alpha)}
\int_z^T\int_t^z
\frac{D^\alpha_{0^+}\nu(\xi)}
{(r-\xi)^{1-\alpha}}
\,d\xi\,d\eta_r
-
I^\alpha_{T^-}(\dot\eta)(t)
\int_t^zD^\alpha_{0^+}\nu(\xi)\,d\xi
\right|
\\
&\qquad
\leq
C_{\eta,\nu,\alpha}|z-t|^{1+\alpha}
=
o(z-t).
\end{align*}

For the second term in
\eqref{example_last_line_ci_derivative}, we have
\[
\int_t^z a(r\mid t,\gamma)\,d\eta_r
=
\lambda(\gamma_t)\dot\eta_t(z-t)
+o(z-t),
\]
where we used
\[
    a(t\mid t,\gamma)=\gamma_t=\lambda(\gamma_t),
\]
which follows from
\(
|\gamma_t|\leq M<\overline M
\).
Therefore,
\[
A(z,\nu)-A(t,\gamma)
=
I^\alpha_{T^-}(\dot\eta)(t)
\int_t^zD^\alpha_{0^+}\nu(\xi)\,d\xi
-
\lambda(\gamma_t)\dot\eta_t(z-t)
+o(z-t),
\]
which proves the preceding expressions for the ci-derivatives of \(A\).

With these calculations at hand, for initial histories satisfying
\[
    \|\gamma\|_{\infty;[0,t]}\leq M,
\]
the value functional for the control problem associated with the
dynamics \eqref{example_dynamics} and cost functional
\eqref{example_cost_functional},
\[
v^{\boldsymbol{\zeta}}(t,x,\gamma)
=
\sup_{u\in L^\infty([t,T])}
J(t,x,\gamma,u),
\]
admits the representation
\begin{align*}
v^{\boldsymbol{\zeta}}(t,x,\gamma)
&=
x
+\sigma\int_t^T a(r\mid t,\gamma)\,d\boldsymbol{\zeta}_r
+c_1\int_t^T
|Z_r|^{\frac{2q}{2q-1}}\,dr,
\end{align*}
where
\[
c_1
=
\frac{2q-1}{2q}
\left(\frac{1}{2qc}\right)^{\frac{1}{2q-1}}
|\sigma|^{\frac{2q}{2q-1}}
\]
and
\[
Z_r
=
\frac{1}{\Gamma(\alpha)}
\int_r^T
\frac{1}{(s-r)^{1-\alpha}}\,d\boldsymbol{\zeta}_s.
\]
The latter integral is understood in the pathwise sense; for its
construction, we refer the reader to
\cite{harang2021volterra}.

Indeed, let \((\eta_m)_{m\in\mathbb{N}}\) be any sequence of smooth
paths whose canonical lifts converge to \(\boldsymbol{\zeta}\) in the
\(p\)-variation topology. The continuity of the rough integral and of
the corresponding Volterra integral, together with the explicit
representations above, implies that
\[
\sup_{(t,x,\gamma)\in\mathcal C}
\left|
v^{\eta_m}(t,x,\gamma)
-
v^{\boldsymbol{\zeta}}(t,x,\gamma)
\right|
\to0
\]
for every compact subset
\[
\mathcal C
\subset
\left\{
(t,x,\gamma)\in
[0,T]\times\mathbb{R}\times AC^\alpha([0,T],\mathbb{R})
:
\|\gamma\|_{\infty;[0,t]}\leq M
\right\}.
\]
Thus the explicit representation is compatible with the smooth
approximations on the region corresponding to the original bounded
trading-rate problem. Independently of this explicit representation,
since \(\lambda\) is bounded and globally Lipschitz, the general
rough-viscosity result established above applies globally. Hence
\(v^{\boldsymbol{\zeta}}\) is a rough viscosity solution in the sense
of Definition \ref{rough_viscosity_solution}.
\end{example}

\section{Concluding remarks and future directions}

In this paper we considered a pathwise control problem driven by a signal of
arbitrary low regularity, for which the controlled dynamics are interpreted in
the sense of rough paths. We allowed the diffusion coefficient to depend on the
control. This dependence is the source of a degeneracy phenomenon: without
additional restrictions, the control may exploit the irregularity of the driving
signal. Following the idea that this degeneracy can be removed by imposing
suitable regularity on the control, together with an appropriate penalisation in
the cost functional, we introduced a fractional formulation in which the control
is represented through its fractional derivative.

The main point of this reformulation is that the control path is no longer
treated as an arbitrary input. Instead, it is reconstructed from a
pseudo-control through fractional integration. This leads to extended dynamics
which keep track of the memory effects induced by the fractional derivative,
while recovering the additive structure needed for dynamic programming. In this
setting we proved well-posedness estimates for the controlled dynamics and the
cost functional, established the dynamic programming principle, and derived the
continuity properties of the value functional.

We then used these results to obtain the HJB characterisation of the value
functional. The resulting equation is fractional and path-dependent, due to
the nonlocal reconstruction of the control, and rough, due to the dependence
on the driving signal. For each smooth approximation of the driving rough
path, the corresponding value functional is characterised as the unique
viscosity solution of the approximating fractional HJB equation within the
class of viscosity solutions whose truncated approximating sequences satisfy
property~{\rm (L)}. Using the stability of the value functional with respect
to the driving rough path, we then showed that the value functional associated
with the rough driver is a rough viscosity solution of the rough fractional
HJB equation in the sense of Definition~\ref{rough_viscosity_solution}.

There are several directions in which the present framework could be extended.
A first direction is suggested by a generalisation of the toy example introduced
in the Introduction. If the informed trader does not have full information about
the external market impact, then the dynamics of his wealth could be modelled by
a rough stochastic differential equation; we refer to \cite{friz2021rough} for a detailed
account of RSDEs. The corresponding control problem has been investigated in
the case where the diffusion coefficient multiplying the rough signal is
independent of the control in \cite{friz2024controlled}. In the setting that
extends the dynamics considered in the present paper, it would be natural to
investigate which conditions on the admissible controls and on the cost
functional ensure that the corresponding control problem is well posed.

A second question concerns the sharpness of the penalisation. As discussed in
Remark~\ref{remark_allan_6_4}, the estimates obtained in this paper improve on
the corresponding estimate in \cite{allan2020pathwise}. It remains unclear,
however, whether this penalisation is optimal, or whether the same
well-posedness results can be obtained under a weaker penalty.

A related, but distinct, question concerns the regularity required of the
control. In this paper we impose enough regularity to ensure complementary
Young regularity with respect to the driver, motivated by the fact that the resulting dynamics is much more natural. The well-posedness of the control
problem when the control has lower regularity remains open.
A further direction is the development of numerical schemes based on the
fractional extended dynamics, with the aim of approximating the associated
pathwise value function and rough fractional HJB equation. Approximation
methods for Caputo derivatives are studied in
\cite{gomoyunov2019approximation}, while approximation schemes for HJB systems
with Riemann--Liouville derivatives are considered in \cite{DB09}. It would be
interesting to investigate whether the schemes in the latter work can be
adapted to the fractional derivatives considered in the former, in order to
solve the equation studied here numerically.

\section*{Acknowledgements}
The authors thank the anonymous referees for their helpful comments and
suggestions, which improved the clarity of the exposition.

\appendix
\begin{appendices}
\section{Existence and stability of the rough control process}
This part of the paper is devoted to the proof of Proposition
\ref{existence_uniqueness_stability_control_process}, which establishes the
global well-posedness and stability of the controlled dynamics
\[
    X_t
    =
    x_0
    +
    \int_0^t b(X_r,\gamma_r)\,dr
    +
    \int_0^t \lambda(\overline X,\gamma)_r\,d\boldsymbol{\zeta}_r,
    \qquad t\in[0,T].
\]
The proof is based on a contraction mapping argument, relying on the
invariance properties of the rough integration map, Lemma
\ref{invariance_of_integration_map}, and on its stability estimates, Lemma
\ref{stability_of_integration_map}.

The main point requiring care is the role of the control \(\gamma\). In the
present framework, \(\gamma\) is itself regarded as a controlled path, and the
solution \(X\) must be estimated in terms of the controlled-path norm of
\(\gamma\). Tracking this dependence quantitatively is important, since the
sharpness of later results in the paper relies on obtaining estimates with the
correct dependence on \(\gamma\). We first prove that the integration map
preserves the relevant controlled-path space, then establish quantitative
stability estimates for this map, and finally use these estimates to obtain
local well-posedness by contraction. The global result follows by
concatenating the local solutions.

Throughout, the notation $A\lesssim_x B$ means that
$A\leq C_x B$ for some constant $C_x>0$ depending only on $x$, whose value may
change from line to line.

\begin{lemma}\label{invariance_of_integration_map} Let \(\boldsymbol{\zeta}\), \(b\), and \(\gamma\) be as in
\eqref{control_process}. Let
\[
    \psi\in C_b^{\floor p + 1}
    (\mathbb R^e\times\mathbb R^\ell,
    \mathcal L(\mathbb R^d,\mathbb R^e))
\]
and let \(\overline X\in\mathcal D_{\boldsymbol{\zeta}}(\mathbb R^e)\).
Define
\[
    Z_t:=z_0+\int_0^t b(X_s,\gamma_s)\,ds
    +\int_0^t \psi(\overline X,\gamma)_s\,d\boldsymbol{\zeta}_s .
\]
Then, for every \(0\leq s<t\leq T\), the controlled rough path lift
\(\overline Z\) satisfies the following remainder estimates:
\begin{align*}
    &\left\|R^{Z}\right\|_{{\frac{p}{\floor{p}}}; [s,t]}  \leq C_{p,b, \psi, \left\|\overline{X}\right\|} \left(t -s + \left\|\boldsymbol{\zeta}\right\|_{p; [s,t ]}\right)\left(1 + \left \|\gamma\right\|_{\frac{p}{\floor{p}}; [s,t ]}  + \sum\limits_{|\delta|=0}^{\floor{p} - 1} \left\|R^{X, \delta}\right\|_{\frac{p}{\floor{p} - |\delta|}; [s, t]} \right), \\ 
    & \left\|R^{Z, \beta}\right\|_{\frac{p}{\floor{p} - |\beta| +1}; [s,t]} \leq C_{\psi, p, \left\|\overline{X}\right\|} \bigg(\left\|\gamma\right\|_{\frac{p}{\floor{p}}; [s,t ]}  + \sum\limits_{|\delta|=0}^{|\beta| - 1} \left\|R^{X, \delta}\right\|_{\frac{p}{\floor{p} - |\delta|}; [s, t]} \bigg), \quad |\beta| \geq 1.
\end{align*}

\end{lemma}
\begin{proof}
    For the trace, using Remark 4.15 in \cite{friz2018differential} conjunction with the inequality
    \begin{equation} \label{bound_increment_controlled_norm}
    \left|X_{st}\right| \leq \sum_{|\tau|=1}^{\floor{p} - 1} |\overline{X}_{\tau; s}\boldsymbol{\zeta}^{\tau}_{st}| + \left|R^X_{st}\right| \leq \left(1 + \left\|\boldsymbol{\zeta}\right\|_{p; [s,t]}\right)\left( |\overline{X}_s| + \sum_{|\tau|=0}^{\floor{p} - 1} |R^{X, \tau}_{st}| \right) \leq C_{L, \left\|\overline{X}\right\|},
    \end{equation}
    yields immediately that whenever $|\beta|=1$,
    \begin{alignat*}{1}
    &\left|R^{Z, \beta}_{st}\right| = \left|R^{\psi(X, \gamma), \beta}_{st}\right| \\ &\leq \sum_{|k| \leq \floor{p}} \bigg| \frac{|k|}{k!} \int_0^1\frac{\partial^{|k|} \psi(X_s + rX_{st}, \gamma_s + r\gamma_{st})_{\beta}}{\partial x^{(k_2, ..., k_{j})} \partial \gamma^{k_{1}}}X^{k_2}_{st}...X^{k_{j-1}}_{st}\gamma^{k_1}_{st} (1-r)^{|k|-1} dr\bigg|\\
    &\qquad + \sum_{|k| = \floor{p}} \bigg|\frac{|k|}{k!}\int_0^1 \frac{\partial^{|k|} \psi(X_s + rX_{st}, \gamma_s + r\gamma_{st})_{\beta}}{\partial x^{k}}X^{k}_{st} (1-r)^{|k|-1} dr\bigg| \\ 
    &\hspace{1cm} + \sum_{{\substack{1\leq n \leq |k| \leq \floor{p}-1 \\|\tau_i| \leq \floor{p} - 1 - |\beta|}}} \bigg|\frac{1}{k!} \frac{\partial^{|k|} \psi(X_{\xi_k}, \gamma_{\xi_k})_{\beta}}{\partial x^k} \overline{X}^{k_1}_{\tau_1, s}...\overline{X}^{k_{n-1}}_{\tau_{n-1}, s}R^{X}_{st}X^{k_{n+1}}_{st}...X^{k_j}_{st} \boldsymbol{\zeta}^{\tau_1}_{st}...\boldsymbol{\zeta}^{\tau_{n-1}}_{st} \bigg| \\
    &\lesssim_{\psi,p, L, \left\|\overline{X}\right\|}(1 + |X_{st}|^{\floor{p}})\bigg(\left\|\gamma\right\|_{\frac{p}{\floor{p}}; [s,t ]}   + \left\|R^{X}\right\|_{\frac{p}{\floor{p}}; [s,t ]} \bigg) \\
    &\lesssim_{\psi,p, L, \left\|\overline{X}\right\|}\bigg(\left\|\gamma\right\|_{\frac{p}{\floor{p}}; [s,t ]}   + \left\|R^{X}\right\|_{\frac{p}{\floor{p}}; [s,t ]} \bigg),
    \end{alignat*}
    where $\xi_k \in (s,t)$.\newline
    If $|\beta| > 1$, relying again on the definition of rough integral and on Remark 4.15 in \cite{friz2018differential} (see also the proof of the second inequality in Proposition \ref{fundamental inequalities}), we have 
    \begin{alignat*}{2}
    &\left|R^{Z, \beta}_{st}\right| \leq  \sum_{{\substack{(\beta^1, ..., \beta^m) \in \overline{Sh}_1^{-1}(\beta^-) \\ 1 \leq n \leq |\beta|\ \\|\tau_i| \leq \floor{p}-1 - |\beta^i|}}} \bigg| \frac{\partial^{|k|} \psi(X_s, \gamma_s)^{\beta^.}}{\partial x^k} \overline{X}^{k_1}_{(\tau_1,\beta^1), s} ...\overline{X}^{k_{n-1}}_{(\tau_{n-1}, \beta^{n-1}), s}  R^{X, \beta^n}_{st}\overline{X}^{k_{n+1}}_{\beta^{n+1}, t} ... \overline{X}^{k_m}_{\beta^m, t} \boldsymbol{\zeta}^{\tau_1}_{st} ... \boldsymbol{\zeta}^{\tau_{n-1}}_{st}\bigg|   \\
    & \quad   + \sum_{{\substack{(\beta^1, ..., \beta^m) \in \overline{Sh}_1^{-1}(\beta^-)  \\  |\tau_i| \leq \floor{p} - 1 - |\beta^i|  \\ 1\leq |k'| \leq \floor{p} - m}}} \bigg| \frac{|k'|}{k'!} \int_0^1 \frac{\partial^{|k|+|k'|} \psi(X_s + rX_{st}, \gamma_s + r\gamma_{st})^{\beta^.}}{\partial x^{(k, k'_1, .., k'_{j-1})} \partial \gamma^{k'_j}} \overline{X}^{k_1}_{(\tau_1, \beta^1), s} ... \overline{X}^{k_m}_{(\tau_m, \beta^m), s} X^{k'^-}_{st}\gamma^{k'^\cdot}_{st} \boldsymbol{\zeta}^{\tau_1}_{st}...\boldsymbol{\zeta}^{\tau_m}_{st}(1-r)^{|k'|-1}dr\bigg| \\ 
    & \quad  +  \sum_{{\substack{(\beta^1, ..., \beta^m) \in \overline{Sh}_1^{-1}(\beta^-)  \\  |\tau_i| \leq \floor{p} - 1 - |\beta^i| \\  m + |k'| = \floor{p}}}}  \bigg| \frac{|k'|}{k'!} \int_0^1 \frac{\partial^{|k|+ |k'|} \psi(X_s + rX_{st}, \gamma_{s} + r\gamma_{st})^{\beta^.}}{\partial x^{(k, k')}} \overline{X}^{k_1}_{(\tau_1, \beta^1), s} ... \overline{X}^{k_m}_{(\tau_m, \beta^m), s} X^{k'}_{st}\boldsymbol{\zeta}^{\tau_1}_{st}...\boldsymbol{\zeta}^{\tau_m}_{st}(1-r)^{|k'|-1}\bigg|  \\ 
    & \quad \begin{aligned} +\sum_{{\substack{(\beta^1, ..., \beta^m) \in \overline{Sh}_1^{-1}(\beta^-) \\  |\tau_i| \leq \floor{p} - 1 - |\beta^i| \\ 1 \leq n \leq |k'| \leq \floor{p} - 1 - m \\  1\leq \delta_i \leq \floor{p}-1}}}  \bigg| \frac{1}{k'!}&\frac{\partial^{|k|+ |k'|} \psi(X_s, \gamma_s)^{\beta^.}}{\partial x^{(k, k')}} \overline{X}^{k_1}_{(\tau_1, \beta^1), s} ... \overline{X}^{k_m}_{(\tau_m, \beta^m), s} \overline{X}^{k'_1}_{\delta_1, s} ... \overline{X}^{k'_{n -1}}_{\delta_{n -1}, s} \\
     &  \hspace{-1cm}\times R^{X}_{st} X^{k'_{n +1}}_{st}...X^{k'_j}_{st} \boldsymbol{\zeta}^{\tau_1}_{st}...\boldsymbol{\zeta}^{\tau_m}_{st}\boldsymbol{\zeta}^{\delta_1}_{st}...\boldsymbol{\zeta}^{\delta_{n-1}}_{st}\bigg|.
     \end{aligned}
    \end{alignat*}
    Applying inequality \eqref{bound_increment_controlled_norm} 
     to the previous inequality yields
    \begin{align*}
        \left|R^{Z, \beta}_{st}\right| &\lesssim_{\psi,p,  L, \left\|\overline{X}\right\|} \left(1 +  |X_{st}|^{\floor{p}}\right) \bigg( \left\| \gamma\right \|_{\frac{p}{\floor{p}}; [s,t]} +  \sum_{|\delta|=0}^{|\beta| - 1} \left\|R^{X, \delta}\right\|_{\frac{p}{\floor{p} - |\delta|}; [s, t]} \bigg) \\
        & \lesssim_{\psi,p,  L, ||\overline{X}||} \bigg( \left\|\gamma\right\|_{\frac{p}{\floor{p}}; [s,t]} +  \sum_{|\delta|=0}^{|\beta| - 1} \left\|R^{X, \delta}\right\|_{\frac{p}{\floor{p} - |\delta|}; [s, t]} \bigg). 
    \end{align*}
    That concludes the proof for the case $|\beta| \geq 1$.\newline
    For the first estimate one has by inequality \eqref{continuity_rough_integral}
    \begin{align*}
     \left|R_{st}^{Z}\right| &= \bigg|\int_{s}^{t} b(X_r, \gamma_r) dr +  \int_{s}^{t} \psi(\overline{X}, \gamma)_r d\boldsymbol{\zeta}_r -  \sum_{|\beta| = 1}^{\floor{p} -1}  \psi(\overline{X}, \gamma)_{\beta; s} \boldsymbol{\zeta}^{\beta}_{st} \bigg| \\
     & \leq \bigg| \int_{s}^{t} b(X, \gamma)_r dr \bigg| +  \bigg| \int_{s}^{t} \psi(\overline{X}, \gamma)_r d\boldsymbol{\zeta}_r  -  \sum_{|\beta| = 1}^{\floor{p}}  \psi(\overline{X}, \gamma)_{\beta; s} \boldsymbol{\zeta}^{\beta}_{st} \bigg|  + \bigg|\sum \limits_{|\beta| = \floor{p}}  \psi(\overline{X}, \gamma)_{\beta; s} \boldsymbol{\zeta}^{\beta}_{st}  \bigg| \\
     &\lesssim_{b, p, \psi} (t-s) + \sum_{|\beta| = 1}^{\floor{p}}\left\|\boldsymbol{\zeta}^\beta\right\|_{\frac{p}{|\beta|}, [s,t]}\left\|R^{\psi, \beta}\right\|_{\frac{p}{\floor{p} - |\beta| + 1}, [s,t]} + \sum_{|\beta| = \floor{p}} \left\|\boldsymbol{\zeta}^{\beta}\right\|_{\frac{p}{\floor{p}}; [s,t]}.
    \end{align*}
    This implies 
    \[
    \left\|R_{st}^{Z}\right\|_{\frac{p}{\floor{p}}; [s,t]} \lesssim _{b, p, \psi} (t-s) + \sum_{|\beta| = 1}^{\floor{p}}\left\|\boldsymbol{\zeta}^\beta\right\|_{\frac{p}{|\beta|}, [s,t]}\left\|R^{\psi, \beta}\right\|_{\frac{p}{\floor{p} - |\beta| + 1}, [s,t]} + \sum_{|\beta| = \floor{p}} \left\|\boldsymbol{\zeta}^{\beta}\right\|_{\frac{p}{\floor{p}}; [s,t]}.
    \]
    Substituting in the bounds obtained for $R^{X, \beta}$, $|\beta|\geq 1$ yields
    \begin{align*}
    \left\|R_{st}^{Z}\right\|_{\frac{p}{\floor{p}}; [s,t]} \lesssim_{p, b, \psi, \left\|\overline{X}\right\|}\left(t -s + \left\|\boldsymbol{\zeta}\right\|_{p; [s,t ]}\right)\bigg(1+ \left\|\gamma\right\|_{\frac{p}{\floor{p}}; [s,t ]}  + \sum\limits_{|\delta|=0}^{\floor{p} - 1} \left\|R^{X, \delta}\right\|_{\frac{p}{\floor{p} - |\delta|}; [s, t]} \bigg), 
    \end{align*}
    which concludes the proof
\end{proof}

\begin{lemma}
[Stability estimates for the integration map] \label{stability_of_integration_map}
Let $\boldsymbol{\zeta}, \boldsymbol{\eta} \in \mathscr{C}^{p}([0,T], \mathbb{R}^d)$ with $\|\boldsymbol{\zeta}\|_{p; [0, T]}, \|\boldsymbol{\eta}\|_{p; [0, T]} \le L$, $X \in \mathcal{D}_{\boldsymbol{\zeta}}(\mathbb{R}^e)$, $Y \in \mathcal{D}_{\boldsymbol{\eta}}(\mathbb{R}^e)$ and $\gamma, \nu \in  \mathcal{C}^{p / \floor{p}}([0,T], \mathbb{R}^\ell)$ satisfying $\left\|\gamma\right\|_{\frac{p}{\floor{p}}; [0,T]}, \left\|\nu\right\|_{\frac{p}{\floor{p}}; [0, T]} < M$. Define the two rough integrals
\begin{equation*}
\begin{aligned}
Z_{t} &= z_0 + \int_{0}^{t} b(X_r, \gamma_r) dr +  \int_{0}^{t} \psi(\overline{X}, \gamma)_r d\boldsymbol{\zeta}_r, \\
V_{t} &= v_0 + \int_{0}^{t} b(Y_r, \nu_r) dr +  \int_{0}^{t} \psi(\overline{Y}, \nu)_r d\boldsymbol{\eta}_r. 
\end{aligned}
\end{equation*}
Then the following inequalities are satisfied
\begin{equation*}
\begin{aligned}
    &\begin{aligned}
        \left\|R^Z_{st} - R^V_{st} \right\|_{\frac{p}{\floor{p}}; [s,t]} \leq  C_{b, p, L, M,  \psi, \left\|\overline{X}\right\|, \left\|\overline{Y}\right\|}\bigg(&(\left\|\gamma_s - \nu_s\right\| + \left\| \gamma - \nu \right\|_{\frac{p}{\floor{p}}; [s,t]} + d(\overline{X}, \overline{Y})_{[s,t]})(t-s + \left\|\boldsymbol{\zeta}\right\|_{p; [s,t]})\\
        &+  \left\|\boldsymbol{\zeta} - \boldsymbol{\eta}\right\|_{p; [s,t]} \bigg),
    \end{aligned}\\
     &\begin{aligned}\left\|R^{Z, \beta}_{st} - R^{V, \beta}_{st} \right\|_{\frac{p}{\floor{p} - |\beta| -1}; [s,t]} \leq  C_{p, L, M,  \psi, \left\|\overline{X}\right\|, \left\|\overline{Y}\right\|}\bigg(&\left\|\gamma_s - \nu_s\right\| + \left\| \gamma - \nu \right\|_{\frac{p}{\floor{p}}; [s,t]} +  \left\|\overline{X}_s - \overline{Y}_s\right\| \\
     & + \sum_{|\delta| \leq|\beta|-1} \left\|R^{X, \delta} -  R^{Y, \delta}\right\|_{\frac{p}{\floor{p} - |\delta|}; [s,t]} + \left\|\boldsymbol{\zeta} - \boldsymbol{\eta}\right\|_{p; [s,t]} \bigg), && |\beta| \geq 1.
     \end{aligned}
\end{aligned}
\end{equation*}
\end{lemma}
\begin{proof}
Using an expansion analogous to the one in the previous proof, we recover that when $|\beta|>1$  
\begin{align}
&|R^{Z, \beta}_{st} - R^{V, \beta}_{st}| = |R^{\psi(X, \gamma), \beta}_{st} - R^{\psi(Y, \nu), \beta}_{st}| \notag \\
&\begin{aligned}
\leq \sum_{{\substack{(\beta^1, ..., \beta^m) \in \overline{Sh}_1^{-1}(\beta^-) \\ 1 \leq n \leq m \\ |\tau_i| \leq \floor{p}-1 - |\beta^i|}}} 
\Bigg| & 
\frac{\partial^{|k|} \psi(X_t, \gamma_t)^{\beta^.}}{\partial x^k} 
\overline{X}^{k_1}_{(\tau_1,\beta^1), s} \cdots \overline{X}^{k_{n-1}}_{(\tau_{n-1}, \beta^{n-1}), s}  
R^{X, \beta^n}_{st} \overline{X}^{k_{n+1}}_{\beta^{n+1}, t} \cdots \overline{X}^{k_m}_{\beta^m, t} 
\boldsymbol{\zeta}^{\tau_1}_{st} \cdots \boldsymbol{\zeta}^{\tau_{n-1}}_{st} \\
& - 
\frac{\partial^{|k|} \psi(Y_t, \nu_t)^{\beta^.}}{\partial x^k} 
\overline{Y}^{k_1}_{(\tau_1,\beta^1), s} \cdots \overline{Y}^{k_{n-1}}_{(\tau_{n-1}, \beta^{n-1}), s}  
R^{Y, \beta^n}_{st} \overline{Y}^{k_{n+1}}_{\beta^{n+1}, t} \cdots \overline{Y}^{k_m}_{\beta^m, t} 
\boldsymbol{\eta}^{\tau_1}_{st} \cdots \boldsymbol{\eta}^{\tau_{n-1}}_{st} 
\Bigg|
\end{aligned} \label{first_line_gub_derivative} \\
&\begin{aligned}
&\quad + \sum_{{\substack{(\beta^1, ..., \beta^m) \in \overline{Sh}_1^{-1}(\beta^-) \\ |\tau_i| \leq \floor{p}-1-|\beta^i| \\ |k'| \leq \floor{p}-m}}} 
\Bigg| \frac{|k'|}{k'!} \int_0^1 
\Bigl( 
\frac{\partial^{|k|+|k'|} \psi(X_s + rX_{st}, \gamma_s + r\gamma_{st})^{\beta^.}}{\partial x^{(k, k'_2, ..., k'_j)} \partial \gamma^{k'_1}} 
\overline{X}^{k_1}_{(\tau_1, \beta^1), s} \cdots \overline{X}^{k_m}_{(\tau_m, \beta^m), s} X^{k'_2}_{st} \cdots X^{k'_j}_{st} \gamma^{k'_1}_{st} \boldsymbol{\zeta}^{\tau_1}_{st} \cdots \boldsymbol{\zeta}^{\tau_m}_{st} \\
&\qquad - 
\frac{\partial^{|k|+|k'|} \psi(Y_s + rY_{st}, \nu_s + r\nu_{st})^{\beta^.}}{\partial x^{(k, k'_2, ..., k'_j)} \partial \gamma^{k'_1}} 
\overline{Y}^{k_1}_{(\tau_1, \beta^1), s} \cdots \overline{Y}^{k_m}_{(\tau_m, \beta^m), s} Y^{k'_2}_{st} \cdots Y^{k'_j}_{st} \nu^{k'_1}_{st} \boldsymbol{\eta}^{\tau_1}_{st} \cdots \boldsymbol{\eta}^{\tau_m}_{st} 
\Bigr) (1-r)^{|k'|-1} dr \Bigg|
\end{aligned} \label{second_line_gub_derivative} \\
&\begin{aligned}
&\quad + \sum_{{\substack{(\beta^1, ..., \beta^m) \in \overline{Sh}_1^{-1}(\beta^-) \\ |\tau_i| \leq \floor{p}-1-|\beta^i| \\ m+|k'| = \floor{p}}}} 
\Bigg| \frac{|k'|}{k'!} \int_0^1 
\Bigl(
\frac{\partial^{|k|+|k'|} \psi(X_s + rX_{st}, \gamma_s + r\gamma_{st})^{\beta^.}}{\partial x^{(k, k')}} 
\overline{X}^{k_1}_{(\tau_1, \beta^1), s} \cdots \overline{X}^{k_m}_{(\tau_m, \beta^m), s} X^{k'_1}_{st} \cdots X^{k'_j}_{st} \boldsymbol{\zeta}^{\tau_1}_{st} \cdots \boldsymbol{\zeta}^{\tau_m}_{st} \\
&\qquad - 
\frac{\partial^{|k|+|k'|} \psi(Y_s + rY_{st}, \nu_s + r\nu_{st})^{\beta^.}}{\partial x^{(k, k')}} 
\overline{Y}^{k_1}_{(\tau_1, \beta^1), s} \cdots \overline{Y}^{k_m}_{(\tau_m, \beta^m), s} Y^{k'_1}_{st} \cdots Y^{k'_j}_{st} \boldsymbol{\eta}^{\tau_1}_{st} \cdots \boldsymbol{\eta}^{\tau_m}_{st} 
\Bigr) (1-r)^{|k'|-1} dr \Bigg|
\end{aligned} \label{third_line_gub_derivative} \\
&\begin{aligned}
&\quad + \sum_{{\substack{(\beta^1, ..., \beta^m) \in \overline{Sh}_1^{-1}(\beta^-) \\ |\tau_i| \leq \floor{p}-1-|\beta^i| \\ 1 \leq n \leq |k'| \leq \floor{p}-1-m \\ 1 \leq \delta_i \leq \floor{p}-1}}} 
\Bigg| \frac{1}{k'!} \frac{\partial^{|k|+|k'|} \psi(X_s, \gamma_s)^{\beta^.}}{\partial x^{(k, k')}} 
\overline{X}^{k_1}_{(\tau_1, \beta^1), s} \cdots \overline{X}^{k_m}_{(\tau_m, \beta^m), s} \overline{X}^{k'_1}_{\delta_1, s} \cdots \overline{X}^{k'_{n-1}}_{\delta_{n-1}, s} \\
&\qquad \times R^{X}_{st} X^{k'_{n+1}}_{st} \cdots X^{k'_j}_{st} \boldsymbol{\zeta}^{\tau_1}_{st} \cdots \boldsymbol{\zeta}^{\tau_m}_{st} \boldsymbol{\zeta}^{\delta_1}_{st} \cdots \boldsymbol{\zeta}^{\delta_{n-1}}_{st} \\
&\qquad - \frac{1}{k'!} \frac{\partial^{|k|+|k'|} \psi(Y_s, \nu_s)^{\beta^.}}{\partial x^{(k, k')}} 
\overline{Y}^{k_1}_{(\tau_1, \beta^1), s} \cdots \overline{Y}^{k_m}_{(\tau_m, \beta^m), s} \overline{Y}^{k'_1}_{\delta_1, s} \cdots \overline{Y}^{k'_{n-1}}_{\delta_{n-1}, s} \\
&\qquad \times R^{Y}_{st} Y^{k'_{n+1}}_{st} \cdots Y^{k'_j}_{st} \boldsymbol{\eta}^{\tau_1}_{st} \cdots \boldsymbol{\eta}^{\tau_m}_{st} \boldsymbol{\eta}^{\delta_1}_{st} \cdots \boldsymbol{\eta}^{\delta_{n-1}}_{st} 
\Bigg|
\end{aligned} \label{fourth_line_gub_derivative}
\end{align}

We start from noticing that for every $\beta$
\begin{align*}
    |\overline{X}_{\beta; t} - \overline{Y}_{\beta; t}| \leq |R^{X, \beta}_{st} - R^{Y, \beta}_{st}| +  \sum_{|\tau|=0}^{\floor{p} -1 - |\beta|} \big| \overline{X}_{(\tau, \beta); s} - \overline{Y}_{(\tau, \beta); s} \big|\big|\boldsymbol{\zeta}_{st}^{\tau}\big|  + \big| \boldsymbol{\zeta}_{st}^{\tau} - \boldsymbol{\eta}_{st}^{\tau} \big| \big| \overline{Y}_{(\tau, \beta); s}\big| .
\end{align*}
A telescopic sum allows to estimate the term \eqref{first_line_gub_derivative} in the previous inequality with 
\begin{align*}
& \begin{aligned}
\sum_{{\substack{(\beta^1, ..., \beta^m) \in \overline{Sh}_1^{-1}(\beta^-) \\ 1 \leq n \leq m \\ |\tau_i| \leq \floor{p}-1 - |\beta^i|}}} \bigg|& \frac{\partial^k \psi(X_t, \gamma_t)^{\beta^.}}{\partial x^k} \overline{X}^{k_1}_{(\tau_1,\beta^1), s} ...\overline{X}^{k_{n-1}}_{(\tau_{n-1}, \beta^{n-1}), s}  R^{X, \beta^n}_{st}\overline{X}^{k_{n+1}}_{\beta_{n+1}, t} ... \overline{X}^{k_m}_{\beta^m, t} \boldsymbol{\zeta}^{\tau_1}_{st} ... \boldsymbol{\zeta}^{\tau_{n-1}}_{st}  \\[-1em]
& - \frac{\partial^k \psi(Y_t, \nu_t)^{\beta^.}}{\partial x^k} \overline{Y}^{k_1}_{(\tau_1,\beta^1), s} ...\overline{Y}^{k_{n-1}}_{(\tau_{n-1}, \beta^{n-1}), s}  R^{Y, \beta^n}_{st}\overline{Y}^{k_{n+1}}_{\beta_{n+1}, t} ... \overline{Y}^{k_m}_{\beta^m, t} \boldsymbol{\eta}^{\tau_1}_{st} ... \boldsymbol{\eta}^{\tau_{n-1}}_{st} \bigg| \end{aligned} \\
& \begin{aligned} \leq 
\sum_{{\substack{(\beta^1, ..., \beta^m) \in \overline{Sh}_1^{-1}(\beta^-) \\ 1 \leq n \leq m \\ |\tau_i| \leq \floor{p}-1 - |\beta^i|}}} \bigg|& \bigg(\frac{\partial^k \psi(X_t, \gamma_t)^{\beta^.}}{\partial x^k}  - \frac{\partial^k \psi(Y_t, \nu_t)^{\beta^.}}{\partial x^k}\bigg) \overline{X}^{k_1}_{(\tau_1,\beta^1), s} ...\overline{X}^{k_{n-1}}_{(\tau_{n-1}, \beta^{n-1}), s}  R^{X, \beta^n}_{st}\\[-1em]
&\times \overline{X}^{k_{n+1}}_{\beta_{n+1}, t} ... \overline{X}^{k_m}_{\beta^m, t} \boldsymbol{\zeta}^{\tau_1}_{st} ... \boldsymbol{\zeta}^{\tau_{n-1}}_{st} \bigg| 
     \end{aligned} \\ 
& \begin{aligned} + 
    \sum_{{\substack{(\beta^1, ..., \beta^m) \in \overline{Sh}_1^{-1}(\beta^-) 1 \\ \leq n \leq m \\ |\tau_i| \leq \floor{p}-1 - |\beta^i| \\ 1 \leq l \leq n-1}}} \bigg|& \frac{\partial^k \psi(Y_t, \nu_t)^{\beta^.}}{\partial x^k} \overline{Y}^{k_1}_{(\tau_1,\beta^1), s} ...(\overline{X}^{k_{l}}_{(\tau_{l}, \beta_{l}), s} - \overline{Y}^{k_{l}}_{(\tau_{l}, \beta_{l}), s} ) ...\overline{X}^{k_{n-1}}_{(\tau_{n-1}, \beta^{n-1}), s}  R^{X, \beta^n}_{st}\overline{X}^{k_{n+1}}_{\beta_{n+1}, t} \\[-1em]
      &  \times ... \overline{X}^{k_m}_{\beta^m, t} \boldsymbol{\zeta}^{\tau_1}_{st} ... \boldsymbol{\zeta}^{\tau_{n-1}}_{st} \bigg| 
     \end{aligned} \\
& \begin{aligned}+ 
    \sum_{{\substack{(\beta^1, ..., \beta^m) \in \overline{Sh}_1^{-1}(\beta^-)  \\ 1 \leq n \leq m \\ |\tau_i| \leq \floor{p}-1 - |\beta^i|}}} \bigg|& \frac{\partial^k \psi(Y_t, \nu_t)^{\beta^.}}{\partial x^k} \overline{Y}^{k_1}_{(\tau_1,\beta^1), s} ...\overline{Y}^{k_{n-1}}_{(\tau_{n-1}, \beta^{n-1}), s}  (R^{X, \beta^n}_{st} -  R^{Y, \beta^n}_{st})\overline{X}^{k_{n+1}}_{\beta_{n+1}, t} ... \overline{X}^{k_m}_{\beta^m, t} \boldsymbol{\zeta}^{\tau_1}_{st} ... \boldsymbol{\zeta}^{\tau_{n-1}}_{st} \bigg| 
     \end{aligned} \\ 
& \begin{aligned} + 
    \sum_{{\substack{(\beta^1, ..., \beta^m) \in \overline{Sh}_1^{-1}(\beta^-) \\ 1 \leq n \leq m \\ |\tau_i| \leq \floor{p}-1 - |\beta^i| \\  n+1 \leq l \leq m}}} \bigg|& \frac{\partial^k \psi(Y_t, \nu_t)^{\beta^.}}{\partial x^k} \overline{Y}^{k_1}_{(\tau_1,\beta^1), s} ...\overline{Y}^{k_{n-1}}_{(\tau_{n-1}, \beta^{n-1}), s}  R^{Y, \beta^n}_{st}\overline{Y}^{k_{n+1}}_{\beta_{n+1}, t} ...(\overline{X}^{k_{n+1}}_{\beta_{n+1}, t} -\overline{Y}^{k_{n+1}}_{\beta_{n+1}, t} )   \\[-1em]
       & \times ... \overline{X}^{k_m}_{\beta^m, t} \boldsymbol{\zeta}^{\tau_1}_{st} ... \boldsymbol{\zeta}^{\tau_{n-1}}_{st} \bigg| 
     \end{aligned} \\ 
& \begin{aligned} + 
    \sum_{{\substack{(\beta^1, ..., \beta^m) \in \overline{Sh}_1^{-1}(\beta^-) \\ 1 \leq n \leq m \\ |\tau_i| \leq \floor{p}-1 - |\beta^i| \\ 1 \leq l \leq n-1}}} \bigg|& \frac{\partial^k \psi(Y_t, \nu_t)^{\beta^.}}{\partial x^k} \overline{Y}^{k_1}_{(\tau_1,\beta^1), s} ...\overline{Y}^{k_{n-1}}_{(\tau_{n-1}, \beta^{n-1}), s}  R^{Y, \beta^n}_{st}\overline{Y}^{k_{n+1}}_{\beta_{n+1}, t} ... \overline{Y}^{k_m}_{\beta^m, t} \boldsymbol{\eta}^{\tau_1}_{st} ... (\boldsymbol{\zeta}^{\tau_l}_{st} - \boldsymbol{\eta}^{\tau_l}_{st}) ... \boldsymbol{\zeta}^{\tau_{n-1}}_{st} \bigg| ,
     \end{aligned}
    \end{align*}
    which implies from standard estimates that the $\frac{p}{\floor{p} - |\beta| +1}$-variation of this first part satisfies the bound
    \begin{align*}
     &\begin{aligned}
     \bigg|& \frac{\partial^k \psi(X_t, \gamma_t)^{\beta^.}}{\partial x^k} \overline{X}^{k_1}_{(\tau_1,\beta^1), s} ...\overline{X}^{k_{n-1}}_{(\tau_{n-1}, \beta^{n-1}), s}  R^{X, \beta^n}_{st}\overline{X}^{k_{n+1}}_{\beta_{n+1}, t} ... \overline{X}^{k_m}_{\beta^m, t} \boldsymbol{\zeta}^{\tau_1}_{st} ... \boldsymbol{\zeta}^{\tau_{n-1}}_{st}  \\
    & - \frac{\partial^k \psi(Y_t, \nu_t)^{\beta^.}}{\partial x^k} \overline{Y}^{k_1}_{(\tau_1,\beta^1), s} ...\overline{Y}^{k_{n-1}}_{(\tau_{n-1}, \beta^{n-1}), s}  R^{Y, \beta^n}_{st}\overline{Y}^{k_{n+1}}_{\beta_{n+1}, t} ... \overline{Y}^{k_m}_{\beta^m, t} \boldsymbol{\eta}^{\tau_1}_{st} ... \boldsymbol{\eta}^{\tau_{n-1}}_{st} \bigg| \end{aligned} \\
     &\lesssim_{p, L, \psi, |\left\|X\right\|, \left\|Y\right\|} \left\|\gamma - \nu\right\|_{\infty} +  \left\|\overline{X}_s - \overline{Y}_s\right\| + \left\|\boldsymbol{\zeta} - \boldsymbol{\eta}\right\|_{p; [s,t]} + \sum_{|\beta^i| < |\beta|}\left\|R^{X, \beta^i} - R^{Y, \beta^i}\right\|_{\frac{p}{\floor{p} - |\beta^i|}; [s, t]}.
    \end{align*}

Using the same procedure one can verify that the same bound holds for the $\frac{p}{p - |\beta| +1}$-variation of the remainders in \eqref{third_line_gub_derivative} and \eqref{fourth_line_gub_derivative}.\newline
A similar result holds for \eqref{second_line_gub_derivative}, with the only difference being the need to use the assumption $\left\| \gamma\right\|_{\frac{p}{\floor{p}}; [0, T]} < M$. 
\begin{align*}
&\begin{aligned}
\sum_{{\substack{(\beta^1, ..., \beta^m) \in \overline{Sh}_1^{-1}(\beta^-)  \\  |\tau_i| \leq \floor{p} - 1 - |\beta^i| \\ |k'| \leq \floor{p}-m }}} \bigg|& \frac{|k'|}{k'!} \int_0^1 \left(\frac{\partial^{|k|+ |k'|} \psi(X_s + rX_{st}, \gamma_s + r\gamma_{st})^{\beta^.}}{\partial x^{(k, k'_2, .., k'_{j})} \partial \gamma^{k'_1}} \overline{X}^{k_1}_{(\tau_1, \beta^1), s} ... \overline{X}^{k_m}_{(\tau_m, \beta^m), s} X^{k'_2}_{st}...X^{k'_{j}}_{st}\gamma^{k'_1}_{st} \boldsymbol{\zeta}^{\tau_1}_{st}...\boldsymbol{\zeta}^{\tau_m}_{st} \right.\\[-2em]
    &- \left. \frac{\partial^{|k|+ |k'|} \psi(Y_s + rY_{st}, \nu_s + r\nu_{st})^{\beta^.}}{\partial x^{(k, k'_2, .., k'_{j-1})} \partial \gamma^{k'_{j}}} \overline{Y}^{k_1}_{(\tau_1, \beta^1), s} ... \overline{Y}^{k_m}_{(\tau_m, \beta^m), s} Y^{k'_2}_{st}...Y^{k'_j}_{st}\nu^{k'_1}_{st} \boldsymbol{\eta}^{\tau_1}_{st}...\boldsymbol{\eta}^{\tau_m}_{st}\right) (1-r)^{|k'|-1} dr 
    \bigg| 
    \end{aligned} \\
     &\lesssim_{p, L,M, \psi, \left\|X\right\|, \left\|Y\right\|} \left\|\gamma - \nu\right\|_{\infty} + \left\|\gamma - \nu\right\|_{\frac{p}{\floor{p}}; [s,t]} + \left\|\overline{X}_s - \overline{Y}_s\right\| + \left\|\boldsymbol{\zeta} - \boldsymbol{\eta}\right\|_{p; [s,t]} +  \sum_{|\beta^i| < |\beta|}\left\|R^{X, \beta^i} - R^{Y, \beta^i}\right\|_{\frac{p}{\floor{p} - |\beta^i|}; [s, t]}. 
\end{align*}

One can easily see that the previous method can be extended to obtain the estimate in the case $|\beta| = 1$, where the quantity of interest is
\begin{align*}
    |R^{Z, \beta}_{st} - R^{Y, \beta}_{st}|
    &= |R^{\psi(X, \gamma), \beta}_{st} - R^{\psi(Y, \gamma), \beta}_{st}| \\
    &\leq \sum_{|k'|\leq \floor{p}} 
    \bigg|
    \frac{|k'|}{k'!} \int_0^1 
    \Bigl(
        \frac{\partial^{|k'|} \psi(X_s + rX_{st}, \gamma_s + r\gamma_{st})_{\beta, \xi_j}}
             {\partial x^{(k_2, ..., k_j)} \partial \gamma^{k_1}}
        X^{k_2}_{st} \cdots X^{k_j}_{st} \gamma^{k_1}_{st} \\
    &\qquad\quad - 
        \frac{\partial^{|k'|} \psi(Y_s + rY_{st}, \nu_s + r\nu_{st})_{\beta}}
             {\partial x^{(k_2, ..., k_j)} \partial \gamma^{k_1}}
        Y^{k_2}_{st} \cdots Y^{k_j}_{st} \nu^{k_1}_{st}
    \Bigr) (1-r)^{|k'|-1} dr
    \bigg| \\
    &\quad + \sum_{|k'| = \floor{p}} 
    \bigg|
    \frac{|k'|}{k!} \int_0^1
    \Bigl(
        \frac{\partial^{|k'|} \psi(X_s + rX_{st}, \gamma_s + r\gamma_{st})_{\beta}}
             {\partial x^{k}}
        X^{j_1}_{st} \cdots X^{j_{|j|}}_{st} \\
    &\qquad\quad - 
        \frac{\partial^{|k'|} \psi(Y_s + rY_{st}, \nu_s + r\nu_{st})_{\beta}}
             {\partial x^{j}}
        Y^{j_1}_{st} \cdots Y^{j_{|j|}}_{st}
    \Bigr) (1-r)^{|k'|-1} dr
    \bigg| \\
    &\quad + \sum_{{\substack{ 1 \leq n \leq |k| \leq \floor{p} \\ |\tau_i| \leq \floor{p} - 1 - |\beta| }}} 
    \bigg|
    \frac{1}{k!} \frac{\partial^{|k|} \psi(X, \gamma)_{\beta, \xi_k}}{\partial x^k} 
    \overline{X}^{j_1}_{\tau_1, s} \cdots \overline{X}^{j_{n-1}}_{\tau_{n-1}, s} 
    R^{X}_{st} X^{k_{n+1}}_{st} \cdots X^{k_j}_{st} 
    \boldsymbol{\zeta}^{\tau_1}_{st} \cdots \boldsymbol{\zeta}^{\tau_{n-1}}_{st} \\
    &\qquad - 
    \frac{1}{k!} \frac{\partial^{|k|} \psi(Y, \nu)_{\beta, \xi_j}}{\partial x^k} 
    \overline{Y}^{k_1}_{\tau_1, s} \cdots \overline{Y}^{k_{n-1}}_{\tau_{n-1}, s} 
    R^{Y}_{st} Y^{k_{n+1}}_{st} \cdots Y^{k_j}_{st} 
    \boldsymbol{\eta}^{\tau_1}_{st} \cdots \boldsymbol{\eta}^{\tau_{n-1}}_{st} 
    \bigg|.
\end{align*}

For the trace we have by inequality \eqref{continuity_rough_integral}
\begin{align*}
    &|R_{st}^{Z} - R_{st}^{V}| \\  
    &= \bigg| \int_{s}^{t} b(X_r, \gamma_r)  -  b(Y_r, \nu_r) dr + \int_{s}^{t} \psi(\overline{X}, \gamma)_r d\boldsymbol{\zeta}_r - \psi(\overline{Y}, \nu)_r d\boldsymbol{\eta}_r  - \bigg( \sum_{|\beta| = 1}^{\floor{p} -1}  \psi^{\beta}(\overline{X}_s, \gamma_s) \boldsymbol{\zeta}^{\beta}_{st} - \psi^{\beta}(\overline{Y}_s, \nu_s) \boldsymbol{\eta}^{\beta}_{st} \bigg) \bigg| \\
    & \leq  \int_{s}^{t} |b(X_s, \gamma_s)  -  b(Y_s, \nu_s)| ds + \bigg | \int_{s}^{t} \psi(\overline{X}, \gamma)_r d\boldsymbol{\zeta}_r - \psi(\overline{Y}, \nu)_r d\boldsymbol{\eta}_r  -  \sum_{|\beta| = 1}^{\floor{p}}  \psi^{\beta}(\overline{X}, \gamma)_s \boldsymbol{\zeta}^{\beta}_{st} - \psi^{\beta}(\overline{Y}_s, \nu_s) \boldsymbol{\eta}^{\beta}_{st} \bigg| + \\ 
    & \quad \quad  + \bigg| \sum_{|\beta| = \floor{p}}  \psi^{\beta}(\overline{X}_s, \gamma_s) \boldsymbol{\zeta}^{\beta}_{st} - \psi^{\beta}(\overline{Y}_s, \nu_s) \boldsymbol{\eta}^{\beta}_{st} \bigg | \\
    & \lesssim_{b, p} \bigg( |X_s - Y_s| + \|X- Y\|_{p; [s,t]}  + \|\gamma - \nu\|_{\infty} \bigg)(t-s) + \sup_{u,v,z} \sum_{|\beta| = 1}^{\floor{p}}|R^{\psi(X, \gamma), \beta}_{uv}\boldsymbol{\zeta}^{\beta}_{vz} - R^{\psi(Y, \nu), \beta}_{uv}\boldsymbol{\eta}^{\beta}_{vz}|  + \\ 
    &\quad \quad +  \sum_{|\beta| = \floor{p}}  |\psi^{\beta}(\overline{Y}_s, \nu_s)|| \boldsymbol{\zeta}^{\beta}_{st} -\boldsymbol{\eta}^{\beta}_{st}| + |\psi^{\beta}(\overline{X}_s, \gamma_s) - \psi^{\beta}(\overline{Y}_s, \nu_s)| |\boldsymbol{\eta}^{\beta}_{st} |.
    \end{align*}
    Using the previous inequalities, this implies
\[  
    \begin{aligned}
     \left\|R^{Z} - R^{V}\right\|_{st} \lesssim_{b,p, L, \psi, M, \left\|X\right\|, \left\|Y\right\|} &(\left\|\gamma_s - \nu_s\right\| + \left\| \gamma - \nu \right\|_{\frac{p}{\floor{p}}; [s,t]} + \left\|\overline{X} - \overline{Y}\right\|_{p;[s,t]})(t-s + \left\|\boldsymbol{\zeta}\right\|_{p; [s,t]})   +  \left\|\boldsymbol{\zeta} - \boldsymbol{\eta}\right\|_{p; [s,t]}.
    \end{aligned}\\
\]
\end{proof}

With these stability estimates we are now ready to prove the existence and uniqueness result stated in Proposition \ref{existence_uniqueness_stability_control_process}

\begin{proof}[Proof of Proposition \ref{existence_uniqueness_stability_control_process}]
This proof follows the proof of Theorem 4.19 in \cite{friz2018differential}.\newline
We begin by defining the seminorms 
\begin{align*}
R^{X, k}_{st} &:= \max_{|\beta|, l \leq k} (R^{X, \beta}_{st} + |X_{st}|^{\floor{p} - l}), && k = 0, ...,  \floor{p} - 1, \\
\|R^{X, k}\|_{\frac{p}{\floor{p} - |\beta|}; [0, T]} &:= \max_{|\beta|\leq k} \|R^{X, \beta}\|_{\frac{p}{\floor{p} - |\beta|}; [0, T]} + \|X\|_{p; [0, T]},  && k = 0, ...,  \floor{p} - 1,
\end{align*}
and the closed set 
\begin{equation*}
B_t :=  \left\{\overline{X} \in \mathcal{D}_{\boldsymbol{\zeta}}(\mathbb{R}^e) : \overline{X}_0 = \left(x_0,  \lambda\left(\overline{X}_0, \gamma_0\right) \right), \|\overline{X}\|_{\mathcal{D}_{\boldsymbol{\zeta}};[0,t]} \leq  |X_0| + 1\, , \, \|R^{X,k}\|_{\frac{p}{\floor{p} - k}; [0,t]} < \delta_k \, , \, k = 0, ..., \floor{p}  - 1 \right\},
\end{equation*}
with $\delta_k > 0$ and the map 
\begin{equation*}
    \overline{\mathcal{M}}^{\gamma}_t : B_t \rightarrow B_t, \quad  \overline{\mathcal{M}}^{\gamma}_t(\overline{X}) := \bigg(x_0 + \int_{0}^{\cdot} b(X_r, \gamma_r) dr + \int_{0}^{\cdot}\lambda(\overline{X}, \gamma)_r d\boldsymbol{\zeta}_r, \lambda(\overline{X}, \gamma)_., ..., \lambda^{\floor{p} -1}(\overline{X}, \gamma)_.  \bigg).
\end{equation*}
The first step consists in showing that this mapping leaves $B_t$ invariant. From Lemma \ref{invariance_of_integration_map} we have
\begin{equation*}
\begin{aligned}
    \left\|R^{\mathcal{M}^{\gamma}, 0}\right\|_{p; [0,t]} 
    &= \left\|R^{\mathcal{M}^{\gamma}}\right\|_{p; [0,t]} + \left\|\mathcal{M}^{\gamma}\right\|_{p; [0,t]} \\
    &\lesssim_{\lambda, p, \left\|\overline{X}\right\|} 
    \Bigl(t  + \left\|\boldsymbol{\zeta}\right\|_{p; [0,t ]}\Bigr)
    \Bigl(1+ \left\|\gamma\right\|_{\frac{p}{\floor{p}}; [0,t]} + \left\|R^{X, \floor{p}-1}\right\|_{p; [0,t]}\Bigr) \\
    &\quad +  \left\|\mathcal{M}^{\gamma}\right\|_{p; [0,t]}.
\end{aligned}
\end{equation*}

and for $k \leq 1$
\begin{align*}
\left\|R^{\mathcal{M}^{\gamma}, k}\right\|_{\frac{p}{\floor{p}-k}; [0,t]} &= \max_{|\beta| \leq k} \left\|R^{\mathcal{M}^{\gamma}, \beta}\right\|_{\frac{p}{\floor{p}-|\beta|}; [0,t]} + \left\|\mathcal{M}^{\gamma}\right\|_{p; [0,t]}\\ 
& \lesssim_{\lambda, p, \left\|\overline{X}\right\|} \bigg(\left\|\gamma\right\|_{\frac{p}{\floor{p}}; [0,t ]} + \left\|R^{X, k-1}\right\|_{\frac{p}{\floor{p} - k + 1}; [0, t]} \bigg) + \left\|\mathcal{M}^{\gamma}\right\|_{p; [0,t]}\\
& \lesssim_{\lambda, p, \left\|\overline{X}\right\|}  \bigg(\left\|\gamma\right\|_{\frac{p}{\floor{p}}; [0,t ]} + \delta_{k-1} \bigg) + \left\|\mathcal{M}^{\gamma}\right\|_{p; [0,t]}. 
\end{align*}
Where the multiplicative constants appearing in the inequalities be chosen to be uniform across all values of k.\newline
Therefore if $t$ is chosen small enough so that \[
\left\|R^{\mathcal{M}^{\gamma}, 0}\right\|_{p; [0,t]} < \delta_0 \text{ and } \hbox{$\delta_k < C_{\lambda, p, \left\|\overline{X}\right\|} \bigg(\left\|\gamma\right\|_{\frac{p}{\floor{p}}; [0,t ]} + \delta_{k-1} \bigg) + \left\|\mathcal{M}^{\gamma}\right\|_{p; [0,t]}$},\]
for every $k \geq 1$, then $B_t$ is invariant under the map $\overline{\mathcal{M}}^{\gamma}_t$.\newline
For the contraction part we first introduce the class of norms
\begin{equation*}
    \left\|\overline{X}\right\|^{\kappa}_{p; [0,t]} := |\overline{X}_0| + \sum_{k = 0}^{\floor{p} - 1} \kappa_k\left\|R^{X, k}\right\|_{\frac{p}{\floor{p} - k}; [0, t]}, 
\end{equation*}
with $\kappa = (\kappa_0, ..., \kappa_{\floor{p} - 1})$ being a vector of positive entries.
For any two controlled paths $\overline{X}$, $\overline{Y} \in B_t$ we notice that 
\begin{align*}
    &\left\|\mathcal{M}^{\gamma}(\overline{X}) - \mathcal{M}^{\gamma}(\overline{Y})\right\|_{p; [0,t]} \\
    &= \left\|\int_0^t b(X_s, \gamma_s) - b(Y_s, \gamma_s)ds + \int_0^t \lambda(\overline{X}, \gamma)_s  - \lambda(\overline{Y}, \gamma)_s d\boldsymbol{\zeta}_s\right\|_{p; [0,t]}\\
    & \lesssim_{p, b, \lambda}\bigg( \sum_{|\beta| \leq \floor{p}-1}\left\|R^{X, \beta} - R^{Y, \beta}\right\|_{p; [0, t]}\bigg)(t + \left\|\boldsymbol{\zeta}\right\|_{p; [0, t]}).
\end{align*}
Then, from Lemma \ref{stability_of_integration_map} it follows that
\begin{align*}
    &\left\|R^{\mathcal{M}^{\gamma}(\overline{X}), 0} - R^{\mathcal{M}^{\gamma}(\overline{Y}), 0}\right\|_{\frac{p}{\floor{p}}; [0, t]}  \\
    &= \left\|R^{\mathcal{M}^{\gamma}(\overline{X})}  -R^{\mathcal{M}^{\gamma}(\overline{Y})}\right\|_{\frac{p}{\floor{p}}; [0, t]} + \left\|\mathcal{M}^{\gamma}(\overline{X}) - \mathcal{M}^{\gamma}(\overline{Y})\right\|_{p; [0, t]} \\
    &\lesssim_{p, b, \lambda, \left\|\overline{X}\right\|, \left\|\overline{Y}\right\|}\bigg( \sum_{|\beta| \leq \floor{p}-1}\left\|R^{X, \beta} - R^{Y, \beta}\right\|_{p; [0, t]}\bigg)\left(t + \left\|\boldsymbol{\zeta}\right\|_{p; [0, t]}\right),
\end{align*}
and when $k \geq 1$
\begin{align*}
    &\left\|R^{\mathcal{M}^{\gamma}(\overline{X}), k} - R^{\mathcal{M}^{\gamma}(\overline{Y}), k}\right\|_{\frac{p}{\floor{p}-k}; [0,t]} \\
    &= \max_{|\beta| \leq k} \left\|R^{\mathcal{M}^{\gamma}(\overline{X}), \beta} - R^{\mathcal{M}^{\gamma}(\overline{Y}), \beta}\right\|_{\frac{p}{\floor{p}-|\beta|}; [0, t]} + \left\|\mathcal{M}^{\gamma}(\overline{X}) - \mathcal{M}^{\gamma}(\overline{Y})\right\|_{p; [0, t]}\\
    &\lesssim_{p, L, \lambda, \left\|\overline{X}\right\|, \left\|\overline{Y}\right\| } \sum_{|\beta|< k} \left\|R^{X, \beta} -  R^{Y, \beta}\right\|_{\frac{p}{\floor{p} - |\delta|}; [0,t]}   + \left\|\mathcal{M}^{\gamma}(\overline{X}) - \mathcal{M}^{\gamma}(\overline{Y})\right\|_{p; [0, t]}\\
    & \lesssim_{p, L, b, \lambda, \left\|\overline{X}\right\|, \left\|\overline{Y}\right\|} \left\|R^{X, k-1} - R^{Y, k-1}\right\|_{\frac{p}{\floor{p} - k+1}; [0, t]} + \bigg( \sum_{|\beta| \leq \floor{p}-1}\left\|R^{X, \beta} - R^{Y, \beta}\right\|_{p; [0, t]}\bigg)\left(t + \left\|\boldsymbol{\zeta}\right\|_{p; [0, t]}\right).
\end{align*}
Which implies that
\begin{align*}
    &\|\mathcal{M}^{\gamma}(\overline{X}) - \mathcal{M}^{\gamma}(\overline{Y})\|^{\kappa}_{p; [0,t]} \\
    &= \kappa_0 C_{p, L, b,  \lambda} d(\overline{X}, \overline{Y})_{[0, t]}(t + \|\boldsymbol{\zeta}\|_{p; [0, t]}) \\
    &\quad \quad + C_{p, L, b, \lambda}\sum_{k=1}^{\floor{p}-1} \kappa_k \bigg(\|R^{X, k-1} - R^{Y, k-1}\|_{\frac{p}{\floor{p} - k+1}; [0,t]} + d(\overline{X}, \overline{Y})_{[0, t]}(t + \|\boldsymbol{\zeta}\|_{p; [0, t]}) \bigg)\\
    &\leq \kappa_0 C_{p, L, b,  \lambda} d(\overline{X}, \overline{Y})_{[0, t]}(t + \|\boldsymbol{\zeta}\|_{p; [0, t]})\\
    &\quad \quad + C_{p, L, b, \lambda} \sum_{k=1}^{\floor{p}-1} \kappa_k  d(\overline{X}, \overline{Y})_{[0, t]}(t + \|\boldsymbol{\zeta}\|_{p; [0, t]})  + C_{p, L, b, \lambda} \sum_{k=1}^{\floor{p}-1} \kappa_k \|R^{X, k-1} - R^{Y, k-1}\|_{\frac{p}{\floor{p} - k+1}; [0,t]}\\
    &\leq C_{\kappa, p, L, b,  \lambda}d(\overline{X}, \overline{Y})^\kappa_{[0, t]}(t+\|\boldsymbol{\zeta}\|_{p; [0, t]}) +  C_{p, L, b, \lambda} \sum_{k=1}^{\floor{p}-1} \kappa_k \|R^{X, k-1} - R^{Y, k-1}\|_{\frac{p}{\floor{p} - k+1}; [0, t]} \\
    &\leq C_{\kappa, p, L, b,  \lambda}d(\overline{X}, \overline{Y})^\kappa_{[0, t]}(t+\|\boldsymbol{\zeta}\|_{p; [0, t]}) +  C_{p, L, b, \lambda} \max_{1\leq k \leq \floor{p} - 1} \frac{\kappa_k}{\kappa_{k-1}} d(\overline{X}, \overline{Y})^\kappa_{[0, t]}.
\end{align*}
Therefore choosing first $\kappa$ in such a way that $C_{ p, L, b, \lambda}  \max\limits_{1\leq k \leq \floor{p} - 1} \frac{\kappa_k}{\kappa_{k-1}} < 1$ and then $t$ small enough that \[C_{\kappa, p, L, b,  \lambda}(t+\|\boldsymbol{\zeta}\|_{p; [0, t]}) < 1- C_{ p, L, b, \lambda}  \max\limits_{1\leq k \leq \floor{p} - 1} \frac{\kappa_k}{\kappa_{k-1}}\] allows to conclude that there exists a unique fixed point of the map $\mathcal{M}^{\gamma}$ over the interval $[0, t]$. Moreover, noticing that the  $t$ was chosen independently of $x_0$ and $\gamma_0$ a global solution for $[0, T]$ can be obtained by pasting together the local solutions. This concludes the contraction part of the argument.\newline
Lastly, using the results of Lemma \ref{stability_of_integration_map} again
\begin{align*}
&\begin{aligned}
    \left\|R^X - R^Y \right\|_{\frac{p}{\floor{p}}; [s,t]}
    &\leq  C_{b, p, L, M,  \psi, \left\|\overline{X}\right\|, \left\|\overline{Y}\right\|} 
    \Bigl(  \left\|\boldsymbol{\zeta} - \boldsymbol{\eta}\right\|_{p; [s,t]} \\
    &\quad + \Bigl(\left\|\gamma_s - \nu_s \right\| + \left\| \gamma - \nu \right\|_{\frac{p}{\floor{p}}; [s,t]} 
    + d(\overline{X}, \overline{Y})_{[s,t]}\Bigr) 
    \Bigl(t-s + \|\boldsymbol{\zeta}\|_{p; [s,t]}\Bigr) \Bigr),
\end{aligned}\\[1em]
&\begin{aligned}
    \left\|R^{\lambda(X, \gamma), \beta} - R^{\lambda(Y, \nu), \beta}\right\|_{p; [s, t]}
    &\leq C_{p, L, \lambda, M, \left\|\overline{X}\right\|, \left\|\overline{Y}\right\|} 
    \Bigg(
        \left\|\gamma_s - \nu_s\right\| + \left\| \gamma - \nu \right\|_{\frac{p}{\floor{p}}; [s,t]} \\
    &\quad +  \left\|\overline{X}_s - \overline{Y}_s\right\| 
      + \sum_{|\delta| \leq|\beta|-1} \left\|R^{X, \delta} -  R^{Y, \delta}\right\|_{\frac{p}{\floor{p} - |\delta|}; [s,t]} \\
    &\quad + \left\|\boldsymbol{\zeta} - \boldsymbol{\eta}\right\|_{p; [s,t]}
    \Bigg).
\end{aligned}
\end{align*}
From these we deduce that there exists a positive map $z \to C(z)$ with $C(z) \to 0^+$ as $z \to 0^+$ and for which
\begin{equation*}
    d(\overline{X},\overline{Y})_{ [s, t]} \leq C(t-s) \bigg(  |X_s - Y_s| + |\gamma_s - \nu_s| + \left\| \gamma - \nu \right\|_{\frac{p}{\floor{p}}; [s,t]} + d(\overline{X},\overline{Y})_{ [s, t]} + \left\|\boldsymbol{\zeta} - \boldsymbol{\eta}\right\|_{p; [s, t]} \bigg).
\end{equation*}
Choosing $t-s$ sufficiently small yields the desired local stability estimate
for the map. Pasting the estimates over finitely many subintervals then
concludes the proof.
\end{proof}

\section{Controlling the remainders}\label{controlling_remainders}

In this appendix we collect the auxiliary estimates that are too technical to
include in the main text but are needed for the proof of Lemma
\ref{lemma_3_6}. Combined with the estimates proved in the main body of the
paper, these results yield a sufficient condition on the running cost
functional ensuring that the control problem under consideration is
well-posed.

The estimates proved here complement the bounds obtained in the main text.
Their purpose is to close the argument by relating the controlled-path norm of
\(X\) to the controlled-path norm of the control \(\gamma\).

The complete argument proceeds in three steps. The first result is presented in Proposition \ref{fundamental inequalities}, where, for any \(\beta \ge1\) and any
\[
    \psi \in
    C_b^{\floor{p}+1}
    \bigl(
        \mathbb{R}^e \times \mathbb{R}^\ell,
        \mathcal{L}(\mathbb{R}^{d},\mathbb{R}^e)
    \bigr),
\]
we use the fact that \(X\) solves an RDE to estimate the relevant variation
norms of \(\psi(\overline X)^\beta\) in terms of $\gamma$ and the remainders of \(X\) up to
level \(|\beta|-1\) . We also show that the remainder of trace of \(\psi(\overline X)\) is
controlled by suitable combinations of the remainders of \(X\) and of the
control \(\gamma\).

Second, by progressively substituting these estimates, we derive bounds on the
remainders of \(X\) in terms only of the remainders of its trace and of
\(\gamma\); this is the content of Lemma
\ref{lemma_remainders_X_final}. Finally, combining this relation with Lemma
\ref{bound_allen_revisited}, which gives a complementary estimate of the trace
remainders in terms of the remainders and \(\gamma\), yields the
desired bound on the controlled-path norm of \(X\) in terms of \(\gamma\), as
stated in Lemma \ref{lemma_3_6}.

Throughout, the notation $A\lesssim_x B$ means that
$A\leq C_x B$ for some constant $C_x>0$ depending only on $x$, whose value may
change from line to line.
 
\begin{remark}\label{preliminary_bound_controlled_process}
A preliminary bound on increment of the Gubinelli derivatives of controlled path $X$ is given by 

\begin{equation*} 
        |\overline{X}_{\beta, st}| \leq C_{\lambda, p}\bigg( \left\|\boldsymbol{\zeta}\right\|_{p; [s, t]} + \left\|R^{X, \beta}\right\|_{\frac{p}{\floor{p} - |\beta|}; [s, t]} \bigg).
    \end{equation*}
\end{remark}

\begin{proposition} \label{fundamental inequalities}
Suppose that $X$ satisfies the RDE \eqref{control_process}, then the following estimates hold:
\begin{equation}
\begin{aligned}
&\left\|R^{\psi(X,\gamma), \beta}\right\|_{\frac{p}{\floor{p} - |\beta| +1}; [s,t]} 
\le C_{\lambda,\psi,p, L}\,(1+\|\gamma\|_{\frac{p}{\floor{p}}; [s,t]}) 
\Big( 1 + \sum_{j=1}^{\floor{p}-1} \|R^X\|_{\frac{p}{\floor{p}}; [s,t]}^j \Big) \quad  |\beta|=1,\\[2mm]
&\left\|R^{\psi(X,\gamma), \beta}\right\|_{\frac{p}{\floor{p} - |\beta| +1}} 
\le C_{\lambda,\psi,p, L}\Bigg(
\sum_{|\beta_i|=1}^{|\beta|-1} \|R^{X, \beta_i}\|_{\frac{p}{\floor{p}-|\beta_i|}; [s,t]} 
+ (1+\|\gamma\|_{\frac{p}{\floor{p}}; [s,t]}) \Big( 1 + \sum_{j=1}^{\floor{p}-1} \|R^X\|_{\frac{p}{\floor{p}}; [s,t]}^j \Big)
\Bigg) \quad |\beta|>1.
\end{aligned}
\end{equation}
\end{proposition}

\begin{proof}
If $|\beta| = 1$, then from Remark 4.15 in \cite{friz2018differential} and Remark \ref{preliminary_bound_controlled_process} it follows immediately that 
\begin{align*}
    \left\|R^{\psi(X,\gamma), \beta}\right\|_{\frac{p}{\floor{p} - |\beta| +1}; [s,t]} \lesssim_{\lambda,\psi,p, L}\left(1+\left\|\gamma\right\|_{\frac{p}{\floor{p}}; [s,t]}\right)\left( 1 + \sum\limits_{j=1}^{\floor{p} - 1 }\left\|R^{X}\right\|^j_{\frac{p}{\floor{p}}; [s,t]} \right).
\end{align*}
For the second estimate, using the definition of controlled paths, we can expand $\overline{X}^{k_{\cdot}}_{\beta_{\cdot}, t}$ to obtain 
\begin{align} \label{taylor_expansion_controlled_estimate}
    \psi(\overline{X}, \gamma)_{\beta, t} &= \sum_{{\substack{(\beta^1, ..., \beta^m) \in \overline{Sh}_1^{-1}(\beta^-) }}} \frac{\partial^{|k|} \psi(X_t, \gamma_t)^{\beta^.}}{\partial x^k} \overline{X}^{k_1}_{\beta^1, t} ... \overline{X}^{k_m}_{\beta^m, t} \notag \\
    &= \sum_{{\substack{(\beta^1, ..., \beta^m) \in \overline{Sh}_1^{-1}(\beta^-)  \\   |\epsilon_i| \leq \floor{p} - 1  - |\beta^i|}}} \frac{\partial^{|k|} \psi(X_t, \gamma_t)^{\beta^.}}{\partial x^k} \overline{X}^{k_1}_{(\epsilon_1, \beta^1), s} ... \overline{X}^{k_m}_{(\epsilon_m, \beta^m), s}  \boldsymbol{\zeta}^{\epsilon_1}_{st}  ... \boldsymbol{\zeta}^{\epsilon_m}_{st} + \tilde{R}^{\psi, \beta}_{st}.
\end{align}
$\tilde{R}^{\psi, \beta}_{st}$ contains at least a factor in $R^{\beta^1}, ..., R^{\beta^m}$ for $|\beta^i| < |\beta|$, therefore the following bound holds for any $0 \leq s \leq t \leq T$
\begin{equation*}
    |\tilde{R}^{\psi, \beta}_{st}| \lesssim_{p, \lambda, \psi, L}\sum_{ 1 \leq |\beta^i| <  |\beta|} |R^{X, \beta^i}_{st}|.
\end{equation*}
Now, using a Taylor expansion for $\frac{\partial^k \psi(X_t, \gamma_t)^{\beta^.}}{\partial x^k}$ around $s$: 
\begin{equation}\label{taylor_expansion_first_estimate}
\begin{aligned}
    &\frac{\partial^{|k|} \psi(X_t, \gamma_t)^{\beta^.}}{\partial x^k} \overline{X}^{k_1}_{(\epsilon_1, \beta^1), s} ... \overline{X}^{k_m}_{(\epsilon_m, \beta^m), s} \\
    &=  \frac{\partial^{|k|} \psi(X_s, \gamma_s)^{\beta^.}}{\partial x^k} \overline{X}^{k_1}_{(\epsilon_1, \beta^1), s} ... \overline{X}^{k_m}_{(\epsilon_m, \beta^m), s} 
      + \sum_{|k| \leq \floor{p} - 1-m} \frac{1}{k'!}\frac{\partial^{|k|+ |k'|)} \psi(X_s, \gamma_s)^{\beta^.}}{\partial x^{(k, k')}} \overline{X}^{k_1}_{(\epsilon_1, \beta^1), s} ... \overline{X}^{k_m}_{(\epsilon_m, \beta^m), s} X^{k'}_{st}  \\
    & \quad +\sum_{|k|\leq \floor{p} - m } \frac{|k|}{k!} \int_0^1 \frac{\partial^{|k|+ |k'|} \psi(X_s + rX_{st}, \gamma_s + r\gamma_{st})^{\beta^.}}{\partial x^{(k, k'_2, .., k'_{l})} \partial \gamma^{k'_1}} \overline{X}^{k_1}_{(\epsilon_1, \beta^1), s} ... \overline{X}^{k_m}_{(\epsilon_m, \beta^m), s} X^{k'_2
    }_{st}...X^{k'_{l}}_{st}\gamma^{k'_1}_{st} (1-r)^{|k'|-1} dr \\
    &\quad + \sum_{|k| = \floor{p} - m } \frac{|k|}{k!} \int_0^1 \frac{\partial^{|k| +|k'|} \psi(X_s + rX_{st}, \gamma_s + r\gamma_{st})^{\beta^.}}{\partial x^{(k, k')}} \overline{X}^{k_1}_{(\epsilon_1, \beta^1), s} ... \overline{X}^{k_m}_{(\epsilon_m, \beta^m), s} X^{k'}_{st}(1-r)^{|k'|-1} dr.
\end{aligned}
\end{equation}
Using the definition of controlled path, we can rewrite 
\begin{align*}
&\sum_{|k| \leq \floor{p} - 1-m} \frac{1}{k'!}\frac{\partial^{|k|+ |k'|} \psi(X_s, \gamma_s)^{\beta^.}}{\partial x^{(k, k')}} \overline{X}^{k_1}_{(\epsilon_1, \beta^1), s} ... \overline{X}^{k_m}_{(\epsilon_m, \beta^m), s} X^{k'}_{st} \\
&= \sum_{\substack{|k| \leq \floor{p} - 1-m \\1 \leq |\delta_1|, ..., |\delta_j| \leq \floor{p} - 1}} \frac{1}{k'!}\frac{\partial^{|k|+| k'|} \psi(X_s, \gamma_s)^{\beta^.}}{\partial x^{(k, k')}} \overline{X}^{k_1}_{(\epsilon_1, \beta^1), s} ... \overline{X}^{k_m}_{(\epsilon_m, \beta^m), s} \overline{X}^{k'_1}_{\delta_1, s}...\overline{X}^{k'_l}_{\delta_l, s} \boldsymbol{\zeta}^{\delta_1}_{st}\dots\boldsymbol{\zeta}^{\delta_l}_{st} + \tilde{\tilde{R}}^{\epsilon_1,\dots, \epsilon_m, k'_1, ..., k'_l}_{st}.
\end{align*}
Where $\tilde{\tilde{R}}^{\epsilon_1, ..., \epsilon_m, k'_1, ..., k'_m}_{st}$ depends on $\psi(X, \gamma)_s$, $X_{st}$ and at least a power of $R^{X}_{st}$, so that
\begin{align*}
    \bigg|\tilde{\tilde{R}}^{\epsilon_1, ..., \epsilon_m, k'_1, ..., k'_m}_{st}\bigg| \lesssim_{L, p, \lambda, \psi} \sum_{j=1}^{\floor{p}-1-m}|R^{X}_{st} |^j.
\end{align*}

Using Remark \eqref{preliminary_bound_controlled_process} and the definition of $\overline X$ as solution to the RDE \eqref{control_process}, we can obtain the following bound for the third term in the left side of \eqref{taylor_expansion_first_estimate}
\begin{align*}
    &\left|\sum_{|k|\leq \floor{p} - m } \frac{|k|}{k!} \int_0^1 \frac{\partial^{|k|+ |k'|} \psi(X_s + rX_{st}, \gamma_s + r\gamma_{st})^{\beta^.}}{\partial x^{(k, k'_2, .., k'_{l})} \partial \gamma^{k'_1}} \overline{X}^{k_1}_{(\epsilon_1, \beta^1), s} ... \overline{X}^{k_m}_{(\epsilon_m, \beta^m), s} X^{k'_2
    }_{st}...X^{k'_{l}}_{st}\gamma^{k'_1}_{st} (1-r)^{|k'|-1} dr \right|\\
    &\lesssim_{\lambda, p, \psi} \bigg( \sum_{j=1}^{\floor{p} - m -1 }\left\|\boldsymbol{\zeta}\right\|^j_{p; [s, t]} + \left\|R^{X}\right\|^j_{\frac{p}{\floor{p}}; [s, t]} \bigg) |\gamma_{st}|. 
\end{align*}
For the last term on the left hand side of \eqref{taylor_expansion_first_estimate},  recalling  again the definition of controlled rough path we get
\begin{align*}
    &\sum_{|k| = \floor{p} - m } \frac{|k|}{k!} \int_0^1 \frac{\partial^{(k, k')} \psi(X_s + rX_{st}, \gamma_s + r\gamma_{st})^{\beta^.}}{\partial x^{(k, k')}} \overline{X}^{k_1}_{(\epsilon_1, \beta^1), s} ... \overline{X}^{k_m}_{(\epsilon_m, \beta^m), s} X^{k'}_{st}(1-r)^{|k'|-1} dr \\ 
    &\lesssim_{\lambda, p, \psi} \bigg( \sum_{|\beta| = \floor{p}} \left\|\boldsymbol{\zeta}^{\beta}\right\|_{\frac{p}{\floor{p}}; [s, t]} +  \sum_{j = 1 }^{\floor{p} -m} \left\|R^{X}\right\|^j_{\frac{p}{\floor{p}}; [s, t]} \bigg).
\end{align*}
The remaining part of the proof, which consists in showing that what we identified as the remainder corresponds to $\psi(X, \gamma)_{\beta, t} - \sum\limits_{|\epsilon|=0}^{\floor{p} - 1 - |\beta|} \psi(X,\gamma)_{(\epsilon, \beta), s}\boldsymbol{\zeta}^{\epsilon}_{st}$
is identical to Remark 4.15 in \cite{friz2018differential}, therefore we omit it. 
\end{proof}

\begin{proof}[Proof of Lemma \ref{lemma_remainders_X_final}]
The second inequality can be shown by replacing $\psi$ with $\lambda$ and applying recursively the 
inequalities in Proposition $\ref{fundamental inequalities}$.\newline
For the trace we have
\begin{align*}
|R^{X}_{st}| &= \bigg| \int_{s}^{t} b(X_s, \gamma_s) ds + \int_{s}^{t} \lambda(\overline{X}, \gamma)_s d\boldsymbol{\zeta}_s  - \sum_{|\beta| = 1}^{\floor{p} -1}  \lambda(\overline{X}, \gamma)_{\beta, s} \boldsymbol{\zeta}^{\beta}_{st} \bigg|  \\
& \leq \bigg| \int_{s}^{t} b(X_s, \gamma_s) ds \bigg| + \bigg| \int_{s}^{t} \lambda(\overline{X}, \gamma)_s d\boldsymbol{\zeta}_s  - \sum_{|\beta| = 1}^{\floor{p}}  \lambda(\overline{X}, \gamma)_{\beta, s} \boldsymbol{\zeta}^{\beta}_{st} \bigg| + \bigg| \sum_{|\beta| = \floor{p}}  \lambda(\overline{X}, \gamma)_{\beta, s} \boldsymbol{\zeta}^{\beta}_{st} \bigg| \\
& \lesssim_{\lambda, b, p}  t-s + \sum_{|\beta| = 1}^{\floor{p}} \| R^{\lambda(X), \beta}\|_{\frac{p}{\floor{p} -|\beta|+1}; [s,t]} \|\boldsymbol{\zeta}\|_{\frac{p}{|\beta|}; [s,t]}  + \sum_{|\beta| = \floor{p}}\|\boldsymbol{\zeta}^{\beta}\|_{\frac{p}{\floor{p}}; [s,t]}. \
\end{align*}
This implies that 
\begin{align*}
\|R^X\|_{\frac{p}{\floor{p}}; [s,t]}  &\lesssim_{\lambda, b, p} (t-s + \|\boldsymbol{\zeta}\|_{p; [s,t]})\left(1 + \sum \limits_{|\beta|=1}^{\floor{p}} \| R^{\lambda(X), \beta}\|_{\frac{p}{\floor{p} -|\beta|+1}; [s,t]}\right)  \\
&\lesssim_{\lambda, b, p} (t-s + \|\boldsymbol{\zeta}\|_{p; [s,t]})(1+\|\gamma\|_{\frac{p}{\floor{p}}; [s,t]})\left( 1 + \sum\limits_{j=1}^{\floor{p} - 1 }\|R^{X}\|^j_{\frac{p}{\floor{p}}; [s, t]} \right) .
\end{align*}
Where in the third step we used the inequality in Proposition \ref{definition_rough_integral} and the previous Remark in the last step.\newline
\end{proof} 

\begin{proof}[Proof of Lemma \ref{bound_allen_revisited}]
We begin by defining the quantity
\[
\hat{R}^X_{st} = X_{st} - \sum_{|\beta| = 1}^{\floor{p}} \lambda(\overline{X}, \gamma)_{{\beta}, s} \boldsymbol{\zeta}^{\beta}_{st}
\]
from the sewing lemma (Lemma 2.6 in \cite{cass2022combinatorial}) and additivity of the increments of $X$, for every $s \leq u \leq t$ we have 
\[
\hat{R}^X_{st} - \hat{R}^X_{su} - \hat{R}^X_{ut} = \sum_{|\beta|=1}^{\floor{p}} R^{\lambda(X), \beta}_{su} \boldsymbol{\zeta}^{\beta}_{ut}.
\]
Iterating this one step identity over any partition $\{s=u_0 < \dots < u_M = t\}$ of $[s, t]$ yields
\begin{equation}\label{R_hat_eq}
\hat{R}_{st} = \sum_{i=1}^{M} \hat{R}_{u_{i-1},u_{i}} + \sum_{i=1}^{M-1} \sum_{|\beta|=1}^{\floor{p}} R^{\lambda(X), \beta}_{u_{i-1} u_{i}} \boldsymbol{\zeta}^{\beta}_{u_{i} t}. 
\end{equation}
Now, introduce the partition $\tilde{\mathcal{P}} = \{0=s_0 < \dots < s_m = T\}$ and label points in the union of the partitions $\mathcal{P}$ and $\tilde{\mathcal{P}}$ in following two distinct ways
\begin{align*}
    s_{j-1} = t_0^j < \dots < t_{n_j}^j = s_j, \quad \text{where } j =1, \dots, m\\
    t_{i-1} = s_0^i < \dots < s_{m_i}^i = t_i, \quad \text{where } i =1, \dots, n.
\end{align*}
Recalling that 
\[
R^X_{st} = \hat R^X_{st} + \sum_{|\beta| = \floor{p}} \lambda (\overline{X}, \gamma)_{\beta, s} \boldsymbol{\zeta}^\beta_{st},
\]
we use \eqref{R_hat_eq} to write for every $j = 1, \dots, m$
\begin{align*}
    \hat{R}^X_{s_{j-1}s_j} = \sum_{i=1}^{n_j} \hat{R}^X_{t^j_{i-1}t^j_{i}} + \sum_{i=1}^{n_j-1} \sum_{|\beta|=1}^{\floor{p}} R^{\lambda(X), \beta}_{t^j_{i-1} t^j_{i}} \boldsymbol{\zeta}^{\beta}_{t^j_{i} s_{j}},
\end{align*}
so that
\[
R^X_{s_{j-1},s_j} = \sum_{i=1}^{n_j} R^X_{t^j_it^j_{i+1}} + \sum_{i=1}^{n_j-1} \sum_{|\beta|=1}^{\floor{p}} R^{\lambda(X), \beta}_{t^j_{i-1} t^j_{i}} \boldsymbol{\zeta}^{\beta}_{t^j_{i} s_j} + \sum_{|\beta| = \floor{p}} \lambda (\overline{X}, \gamma)_{\beta, s_{j-1}} \boldsymbol{\zeta}^\beta_{s_{j-1}s_{j}} - \sum_{i=1}^{n_j} \sum_{|\beta| = \floor{p}} \lambda(\overline{X}, \gamma)_{\beta, t^j_{i-1}} \boldsymbol{\zeta}^\beta_{t^j_{i-1}t^j_{i}}.
\]
This identity can now be used to write for some $C_{p, \lambda} > 0$
\begin{align*}
    &\sum_{j=1}^{m} \left|R^X_{s_{j-1}s_j}\right|^{\frac{p}{\floor{p}}}\\
    &= \sum_{j=1}^{m} \left|\sum_{i=1}^{n_j} R^X_{t^j_{i-1},t^j_{i}} + \sum_{i=1}^{n_j-1} \sum_{|\beta|=1}^{\floor{p}} R^{\lambda(X), \beta}_{t^j_{i-1} t^j_{i}} \boldsymbol{\zeta}^{\beta}_{t^j_{i} s_j} - \sum_{|\beta| = \floor{p}} \lambda (\overline{X}, \gamma)_{\beta, s_{j-1}} \boldsymbol{\zeta}^\beta_{s_{j-1}s_{j}} + \sum_{i=1}^{n_j} \sum_{|\beta| = \floor{p}} \lambda(\overline{X}, \gamma)_{\beta, t^j_{i-1}} \boldsymbol{\zeta}^\beta_{t^j_{i-1}t^j_{i}}\right|^{\frac{p}{\floor{p}}}\\
    &\leq C_{p, \lambda}  n^{\frac{p}{\floor{p}} -1} \left( \sum_{j=1}^{m} \sum_{i=1}^{n_j}\left| R^X_{t^j_{i-1}t^j_{i}}\right|^{\frac{p}{\floor{p}}} + \sum_{|\beta| = \floor{p}}  \left\| \boldsymbol{\zeta}^\beta\right\|^{\frac{p}{\floor{p}}}_{\frac{p}{\floor{p}}; [0, T]} \right) + C_{p} \sum_{|\beta|=1}^{\floor{p}} \sum_{j=1}^{m} \left| \sum_{i=1}^{n_j-1}  \left|R^{\lambda(X), \beta}_{t^j_{i-1} t^j_{i}}\right| \|\boldsymbol{\zeta}^\beta\|_{\frac{p}{|\beta|};[s_{j-1}, s_j]} \right|^{\frac{p}{\floor{p}}}.
    \end{align*}
    Here we used H\"older's inequality in the last step.\newline
    From this point onwards, the constant $C_{p,\lambda}>0$ may change from line to line. We keep track only of the relevant dependencies, enlarging the subscript when further parameters enter the estimates.\newline
    To estimate the last term in the preceding inequality, we distinguish between the cases $|\beta|>1$ and $|\beta|=1$. 
If $|\beta|>1$, we apply Young's inequality for products with exponents
\(
\frac{\lfloor p \rfloor}{\lfloor p \rfloor-|\beta|+1}
\text{ and }
\frac{\lfloor p \rfloor}{|\beta|-1}.
\)
The case $|\beta|=1$ has to be treated separately, since the first exponent degenerates to $1$. In this case we simply use the bound
\[
\left|R^{\lambda(X),\beta}_{t^j_{i-1}t^j_i}\right|
\|\boldsymbol{\zeta}^{\beta}\|_{\frac{p}{|\beta|};[s_{j-1},s_j]}
\leq 
L\left|R^{\lambda(X),\beta}_{t^j_{i-1}t^j_i}\right|.
\]
Combining these two estimates gives
    \begin{align*}
    &\sum_{j=1}^{m} \left|R^X_{s_{j-1}s_j}\right|^{\frac{p}{\floor{p}}}\\
    &\leq C_{p, \lambda}  n^{\frac{p}{\floor{p}} -1} \left( \sum_{j=1}^{m} \sum_{i=1}^{n_j}\left| R^X_{t^j_{i-1}t^j_{i}}\right|^{\frac{p}{\floor{p}}} + \sum_{|\beta| = \floor{p}}  \left\| \boldsymbol{\zeta}^\beta\right\|^{\frac{p}{\floor{p}}}_{\frac{p}{\floor{p}}; [0, T]} \right)\\
    &\qquad + C_{p, L} \left(\sum_{|\beta|=1}^{\floor{p}} \sum_{j=1}^{m} \left| \sum_{i=1}^{n_j-1}  \left|R^{\lambda(X), \beta}_{t^j_{i-1} t^j_{i}}\right|\right| ^{\frac{p}{\floor{p} - |\beta| + 1}} +  \sum_{|\beta|=2}^{\floor{p}} \sum_{j=1}^{m} \left\| \boldsymbol{\zeta}\right\| ^{\frac{p}{|\beta|-1}}_{\frac{p}{|\beta|};[s_{j-1}, s_j]}\right)
    \end{align*}
    Using H\"older's inequality for the penultimate term allows now to find 
    \begin{align*}
    &\sum_{j=1}^{m} \left|R^X_{s_{j-1}s_j}\right|^{\frac{p}{\floor{p}}}\\
    &\leq C_{p, L, \lambda}  n^{\frac{p}{\lfloor p\rfloor}-1} \left( \sum_{i=1}^{n} \sum_{j=1}^{m_i}\left| R^X_{s^i_{j-1}s^i_{j}}\right|^{\frac{p}{\floor{p}}} +  \sum_{i=1}^{n}  \sum_{|\beta|=1}^{\floor{p}} \sum_{i=1}^{m_i-1}  \left| R^{\lambda(X), \beta}_{s^i_{j-1} s^i_{j}}\right|^{\frac{p}{\floor{p} - |\beta| +1}} +  \left\| \boldsymbol{\zeta}\right\|_{p; [0, T]} \right)\\
    &\leq C_{p, L, \lambda}  n^{\frac{p}{\lfloor p\rfloor}-1} \left( \sum_{i=1}^{n} \left\| R^X\right\|^{\frac{p}{\floor{p}}}_{{\frac{p}{\floor{p}}}; [{t_{i-1},t_{i}}]} +  \sum_{i=1}^{n}  \sum_{|\beta|=1}^{\floor{p}}   \left\| R^{\lambda(X), \beta}\right\|^{\frac{p}{\floor{p} - |\beta| +1}}_{{\frac{p}{\floor{p} - |\beta| +1}}; [t_{i-1}, t_i]} +  \left\| \boldsymbol{\zeta}\right\|_{p; [0, T]} \right).
\end{align*}
The inequality above is sufficient to conclude the proof after taking the supremum over $\tilde{\mathcal{P}}$.
\end{proof}

\section{Fundamentals of Fractional Differentiation and Integration}\label{fundamentals_fractional_integrations}

\begin{definition}
Let $r < T$ and $\alpha > 0$. For every $r \leq t \leq T$, the left-sided Riemann--Liouville fractional integral of order $\alpha$ with base point $r$ of a function $u \in L^1\big([r,T], \mathbb{R}^\ell\big)$ is defined by
\begin{equation*}
    I_{r^+}^{\alpha}u(t) := \frac{1}{\Gamma(\alpha)} \int_r^t \frac{u_s}{(t-s)^{1-\alpha}} \, ds .
\end{equation*}
Here, $\Gamma$ denotes the gamma function.
We will denote by $I_{r^+}^{\alpha}\big(L^1\big([r,T], \mathbb{R}^\ell\big)\big)$ the class of functions that can be represented in the form above.\newline
Similarly, the right-sided Riemann--Liouville fractional integral of order $\alpha$ with base point $T$ is given by
\begin{equation*}
    I_{T^-}^{\alpha}u(t) := \frac{1}{\Gamma(\alpha)} \int_t^T \frac{u_s}{(s-t)^{1-\alpha}} \, ds .
\end{equation*}

\end{definition}

\begin{definition}
The left-sided Riemann-Liouville derivative of order $\alpha \in (0, 1)$ with base point $r$ of a function $u \in I_{r^+}^{\alpha}\big(L^1\big([r,T], \mathbb{R}^\ell\big)\big)$ is given by
\begin{equation*}
    D_{r^+}^{\alpha}u(t) = \frac{1}{\Gamma(1- \alpha)} \frac{d}{dt} \int_r^t \frac{u_s}{(t-s)^\alpha} ds.
\end{equation*}
Here, $\Gamma$ denotes the gamma function.
\end{definition}
\noindent We define the space $AC^{\alpha}([0, T], \mathbb{R}^\ell)$ to be the class of functions  $\gamma$ that can be  expressed as
\begin{equation} \label{class_AC_alpha}
    \gamma_t = \gamma_0 + I_{0^+}^{\alpha}u (t), \quad u \in L^{\infty}([0,T], \mathbb{R}^\ell).
\end{equation}

\begin{definition}\label{definition_alpha_holder}
    Let $\alpha \in (0,1]$. A continuous path $\gamma: [0,T] \to \mathbb{R}^\ell$ is said to belong to $C^{\text{H\"ol}-\alpha}([0,T], \mathbb{R}^\ell)$, if  the following inequality holds:
\[
\|\gamma\|_{\alpha-\text{H\"ol}} := \sup_{0 \leq s < t \leq T} \frac{|\gamma_t - \gamma_s|}{|t-s|^{\alpha}} < \infty.
\]
\end{definition}
\noindent The following proposition presents fundamental properties of functions belonging to the class $AC^\alpha$, which will be utilized frequently in the subsequent sections.
\begin{proposition}\label{properties_AC}
\leavevmode
Let \(\alpha\in(0,1)\). For every \(k\in\mathbb N\), define
\[
AC_k^\alpha
:=
\left\{
\gamma\in AC^\alpha([0,T],\mathbb R^\ell):
\left\|D_{0^+}^\alpha(\gamma-\gamma_0)\right\|_\infty\leq k,
\ |\gamma_0|\leq k
\right\}.
\]
Then:
\begin{enumerate}
    \item
    \(
    AC^\alpha([0,T],\mathbb R^\ell)
    \subset
    \mathcal C^{\mathrm{H\ddot{o}l}-\alpha}
    ([0,T],\mathbb R^\ell).
    \)

    \item
    \(
    D_{0^+}^\alpha(\gamma-\gamma_0)(t)=u(t)
    \)
    for every \(\gamma\) defined as in \eqref{class_AC_alpha}.

    \item For every \(r>1/\alpha\), the set \(AC_k^\alpha\) is compact
    with respect to the topology induced by
    \[
    d_r(\gamma,\nu)
    :=
    |\gamma_0-\nu_0|
    +
    \|\gamma-\nu\|_{r\text{-var};[0,T]}.
    \]
    Consequently, \(AC^\alpha\) is \(\sigma\)-compact in this topology.
    In particular, if
    \(
    \alpha> \floor{p} / p,
    \)
    this conclusion applies with \(r=p/\floor{p}\).
\end{enumerate}
\end{proposition}

\begin{proof}
The first and second claims follow, respectively, from Theorem~3.1 and
Theorem~2.4 in \cite{samko1993fractional}.

We prove the third claim. Fix \(k\in\mathbb N\) and let
\((\gamma_n)_{n\in\mathbb N}\subset AC_k^\alpha\). Write
\[
u_n:=D_{0^+}^\alpha(\gamma_n-\gamma_{n,0}),
\]
so that
\[
\gamma_n(t)
=
\gamma_{n,0}
+
I_{0^+}^\alpha u_n(t),
\qquad
\|u_n\|_\infty\leq k,
\qquad
|\gamma_{n,0}|\leq k.
\]
The continuity properties of the fractional integral imply that there
exists a constant \(C_{\alpha,T}>0\), independent of \(n\), such that
\[
\|\gamma_n\|_\infty
+
\|\gamma_n\|_{\alpha\text{-Höl}}
\leq C_{\alpha,T}k.
\]
Thus, the sequence \((\gamma_n)_{n\in\mathbb N}\) is uniformly bounded
and equicontinuous. By the Arzel\`a--Ascoli theorem, after passing to a
subsequence, there exists
\(\gamma\in C([0,T],\mathbb R^\ell)\) such that
\[
\|\gamma_n-\gamma\|_\infty\to0.
\]

Since \((u_n)_{n\in\mathbb N}\) is bounded in
\(L^\infty([0,T],\mathbb R^\ell)\), after passing to a further
subsequence, there exists \(u\in L^\infty([0,T],\mathbb R^\ell)\)
such that
\[
u_n \overset{*}{\rightharpoonup} u
\quad\text{in } L^\infty([0,T],\mathbb{R}^\ell),
\qquad
\|u\|_\infty\leq k.
\]
For every \(t\in[0,T]\), the function
\[
s\to
\frac{\mathbbm 1_{[0,t]}(s)}
{\Gamma(\alpha)(t-s)^{1-\alpha}}
\]
belongs to \(L^1([0,T])\). Hence the weak-star convergence gives
\[
I_{0^+}^\alpha u_n(t)
\to
I_{0^+}^\alpha u(t).
\]
Since \(\gamma_{n,0}\to\gamma_0\), it follows that
\(
\gamma(t)=\gamma_0+I_{0^+}^\alpha u(t).
\)
Therefore,
\(
\gamma\in AC_k^\alpha.
\)

It remains to strengthen the uniform convergence to convergence in
\(r\)-variation. Fix \(r>1/\alpha\), and choose
\(
1/r<\beta<\alpha.
\)
The interpolation inequality for H\"older seminorms gives
\[
\|\gamma_n-\gamma\|_{\beta\text{-Höl}}
\leq
C_{\alpha,\beta}
\|\gamma_n-\gamma\|_\infty^{1-\beta/\alpha}
\|\gamma_n-\gamma\|_{\alpha\text{-Höl}}^{\beta/\alpha}.
\]
The \(\alpha\)-H\"older seminorms on the right-hand side are uniformly
bounded, while
\(\|\gamma_n-\gamma\|_\infty\to0\). Consequently,
\(
\|\gamma_n-\gamma\|_{\beta\text{-Höl}}
\to0.
\)
Since \(\beta r>1\), the embedding of \(\beta\)-H\"older paths into
paths of finite \(r\)-variation yields
\[
\|\gamma_n-\gamma\|_{r\text{-var};[0,T]}
\leq
T^\beta
\|\gamma_n-\gamma\|_{\beta\text{-H\"ol}}
\to0.
\]
Together with \(\gamma_{n,0}\to\gamma_0\), this proves that
\(AC_k^\alpha\) is compact with respect to \(d_r\).

Finally, every \(\gamma\in AC^\alpha\) belongs to \(AC_k^\alpha\) for
some \(k\in\mathbb N\). Hence
\[
AC^\alpha
=
\bigcup_{k=1}^\infty AC_k^\alpha,
\]
which proves the asserted \(\sigma\)-compactness. Taking
\(r=p/\floor{p}\) is admissible whenever
\(\alpha>\floor{p}/p\).
\end{proof}
In this work, we will use the operator $\gamma \rightarrow D_{0^+}^\alpha(\gamma - \gamma_0)$, known as the Caputo differential operator, which coincides with the Caputo derivative when $\gamma \in AC^1$. For further properties of these operators, the reader is referred to \cite{diethelm2010analysis}. Additional properties of the space $AC^\alpha$ and a detailed proof of the last  property of the previous proposition can be found in \cite{gomoyunov2020theory}.

\begin{lemma}
\label{lem:fractional_continuation_estimates}
Let $\alpha\in(0,1)$, $r\in[0,T]$, and
$\gamma\in AC^\alpha([0,T],\mathbb R^\ell)$. For $t\in[r,T]$, define
\[
    \nu_t
    :=
    \gamma_0
    +
    \frac{1}{\Gamma(\alpha)}
    \int_0^r
    \frac{D^\alpha_{0^+}(\gamma-\gamma_0)(s)}
         {(t-s)^{1-\alpha}}
    \,ds.
\]
Then $\nu_r=\gamma_r$, and, for every $t>r$,
\begin{equation}
\label{eq:fractional_continuation_representation}
    \nu_t
    =
    \int_0^\infty
    \gamma_{\left(r-(t-r)y\right)_+}
    \rho_\alpha(y)\,dy,
\end{equation}
where
\[
    \rho_\alpha(y)
    :=
    \frac{1}{\Gamma(\alpha)\Gamma(1-\alpha)}
    \frac{y^{-\alpha}}{1+y},
    \qquad y>0,
\]
and
\[
    \int_0^\infty \rho_\alpha(y)\,dy=1.
\]
In particular, for every $t,z\in[r,T]$,
\begin{equation}
\label{eq:fractional_continuation_oscillation}
    |\nu_t-\nu_z|
    \le
    \sup_{a,b\in[0,r]}|\gamma_a-\gamma_b|
    \le
    2\|\gamma-\gamma_0\|_{\infty;[0,r]}.
\end{equation}
Moreover, if $\gamma$ has finite $q$-variation on $[0,r]$ for some
$q\ge1$, then
\begin{equation}
\label{eq:fractional_continuation_qvar}
    \|\nu\|_{q\text{-var};[r,T]}
    \le
    \|\gamma\|_{q\text{-var};[0,r]}.
\end{equation}
\end{lemma}

\begin{proof}
Set
\[
    h:=\gamma-\gamma_0,
    \qquad
    \varphi:=D^\alpha_{0^+}h.
\]
Since $\gamma\in AC^\alpha$, we have
\(
    h=I^\alpha_{0^+}\varphi.
\)
In particular,
\[
    \nu_r
    =
    \gamma_0
    +
    I^\alpha_{0^+}\varphi(r)
    =
    \gamma_0+h(r)
    =
    \gamma_r.
\]

Now let $t>r$. By the Riesz identity for fractional integrals (identity 14.36 in \cite{samko1993fractional}),
\[
    \int_0^r
    (t-s)^{\alpha-1}\varphi(s)\,ds
    =
    \frac{(t-r)^\alpha}{\Gamma(1-\alpha)}
    \int_0^r
    \frac{h(v)}
         {(t-v)(r-v)^\alpha}\,dv.
\]
Consequently
\begin{equation}
\label{eq:nu_riesz_step}
    \nu_t-\gamma_0
    =
    \frac{(t-r)^\alpha}
         {\Gamma(\alpha)\Gamma(1-\alpha)}
    \int_0^r
    \frac{\gamma_v-\gamma_0}
         {(t-v)(r-v)^\alpha}\,dv.
\end{equation}

Set $\delta=t-r>0$ and make the change of variables $x=r-v$.
Since $t-v=\delta+x$, \eqref{eq:nu_riesz_step} becomes
\[
    \nu_t-\gamma_0
    =
    \frac{\delta^\alpha}
         {\Gamma(\alpha)\Gamma(1-\alpha)}
    \int_0^r
    \frac{\gamma_{r-x}-\gamma_0}
         {(\delta+x)x^\alpha}\,dx.
\]
Making the further substitution $x=\delta y$ gives
\[
    \nu_t-\gamma_0
    =
    \frac{1}{\Gamma(\alpha)\Gamma(1-\alpha)}
    \int_0^{r/\delta}
    \left(\gamma_{r-\delta y}-\gamma_0\right)
    \frac{y^{-\alpha}}{1+y}\,dy.
\]
Since $\gamma_0-\gamma_0=0$, the upper limit may be extended to
infinity by introducing
\[
    \theta_t(y):=\left(r-(t-r)y\right)_+.
\]
Hence
\begin{equation}
\label{eq:nu_before_probability}
    \nu_t-\gamma_0
    =
    \int_0^\infty
    \left(\gamma_{\theta_t(y)}-\gamma_0\right)
    \rho_\alpha(y)\,dy.
\end{equation}
On the other hand, by the beta-integral identity,
\[
    \int_0^\infty \frac{y^{-\alpha}}{1+y}\,dy
    =
    B(1-\alpha,\alpha)
    =
    \Gamma(1-\alpha)\Gamma(\alpha).
\]
Therefore
\[
    \int_0^\infty\rho_\alpha(y)\,dy=1,
\]
and \eqref{eq:nu_before_probability} yields
\eqref{eq:fractional_continuation_representation}.

It follows that, for $t,z\in[r,T]$,
\[
    \nu_t-\nu_z
    =
    \int_0^\infty
    \left(
        \gamma_{\theta_t(y)}
        -
        \gamma_{\theta_z(y)}
    \right)
    \rho_\alpha(y)\,dy.
\]
Thus
\[
    |\nu_t-\nu_z|
    \le
    \int_0^\infty
    \left|
        \gamma_{\theta_t(y)}
        -
        \gamma_{\theta_z(y)}
    \right|
    \rho_\alpha(y)\,dy
    \le
    \sup_{a,b\in[0,r]}|\gamma_a-\gamma_b|,
\]
which proves \eqref{eq:fractional_continuation_oscillation}.

Finally, suppose that $\gamma$ has finite $q$-variation and let
\[
    r\le t_0<t_1<\cdots<t_n\le T.
\]
By Minkowski's inequality,
\begin{align*}
    \left(
        \sum_{i=1}^n
        |\nu_{t_i}-\nu_{t_{i-1}}|^q
    \right)^{1/q}
    &\le
    \int_0^\infty
    \left(
        \sum_{i=1}^n
        \left|
            \gamma_{\theta_{t_i}(y)}
            -
            \gamma_{\theta_{t_{i-1}}(y)}
        \right|^q
    \right)^{1/q}
    \rho_\alpha(y)\,dy.
\end{align*}
For each fixed $y>0$, the map $t\mapsto\theta_t(y)$ is
non-increasing. Hence
\[
    \theta_{t_0}(y)
    \ge
    \theta_{t_1}(y)
    \ge
    \cdots
    \ge
    \theta_{t_n}(y),
\]
so, after reversing the order and removing possible repetitions, these
points form a partition of a subinterval of $[0,r]$. Therefore
\[
    \left(
        \sum_{i=1}^n
        \left|
            \gamma_{\theta_{t_i}(y)}
            -
            \gamma_{\theta_{t_{i-1}}(y)}
        \right|^q
    \right)^{1/q}
    \le
    \|\gamma\|_{q\text{-var};[0,r]}.
\]
Using again that $\rho_\alpha$ integrates to one gives
\[
    \left(
        \sum_{i=1}^n
        |\nu_{t_i}-\nu_{t_{i-1}}|^q
    \right)^{1/q}
    \le
    \|\gamma\|_{q\text{-var};[0,r]}.
\]
Taking the supremum over all partitions of $[r,T]$ proves
\eqref{eq:fractional_continuation_qvar}.
\end{proof}

\section{Auxiliary results}
In this appendix we collect several auxiliary results adapted from
Section~5 of \cite{gomoyunov2021viscosity}. They are included for the
reader's convenience and to make explicit the notation used in the present
setting. The proofs omitted below are obtained by arguments entirely
analogous to those in \cite{gomoyunov2021viscosity}, up to the
notational changes required here. We spell out the proof of Lemmas~\ref{lemma_5_7_g21} and \ref{lemma_5_8_g21} 
only to clarify minor changes in the arguments.
\begin{lemma}[Adapted from Lemma 5.1 in \cite{gomoyunov2021viscosity}]\label{lemma_5_4_g21}
    For every \(\epsilon>0\), the functional \(\varpi_\epsilon\) defined in
\eqref{auxiliary_functional} is non-negative and continuous.
Moreover, with \(a(\cdot\mid t,\gamma)\) defined as in
\eqref{auxiliary_functional_a}, the following identities hold:
    \begin{align*}
        \varpi_\epsilon(t, \gamma, t, \gamma) = 0, \quad \varpi_\epsilon(t, \gamma, \tau, \nu) = \varpi_\epsilon(\tau, \nu, t,  \gamma)
    \end{align*}
    and 
    \begin{align*}
        \varpi_\epsilon(\tilde{t}, a(\cdot| t, \gamma), \tilde \tau, a(\cdot| \tau, \nu)) = \varpi_\epsilon(t, \gamma, \tau, \nu)
    \end{align*}
    are valid for every $(t, \gamma), (\tau, \nu) \in [0, T] \times AC^\alpha([0, T], \mathbb{R}^\ell)$ and $\tilde t \in [t, T], \tilde \tau \in [\tau, T]$.
\end{lemma}

\begin{lemma}[Lemma 5.2 in \cite{gomoyunov2021viscosity}]\label{lemma_5_5_g21}
Let \(a\) be defined as in \eqref{auxiliary_functional_a}, and let
\(\varpi_\epsilon\) be defined as in \eqref{auxiliary_functional}.
Then there exists \(C>0\) such that
\begin{align}
&\left|a\bigl(T \mid t,\gamma\bigr)
      - a\bigl(T \mid \tau, \nu \bigr)\right|
 + \int_0^T
 \frac{
 \left|a\bigl(r \mid t, \gamma \bigr)
      - a\bigl(r \mid \tau, \nu\bigr)\right|
 }{(T-r)^{1-\alpha}}
 \,dr
\nonumber \\
&\qquad\leq
C \left(
\varpi_\epsilon\bigl(t,\gamma,\tau,\nu\bigr)
+ C_1 \epsilon^{\frac{c}{c-1}}
\right)^{\frac{1}{c}}
\label{eq:5.4}
\end{align}
for any $\epsilon > 0$, any
$(t,\gamma),(\tau,\nu) \in [0,T] \times AC^\alpha([0, T], \mathbb{R}^\ell)$ and where \(C_1\) is the constant appearing in \eqref{auxiliary_functional}.
\end{lemma}

\begin{lemma}[Adapted from Lemma 5.3 in \cite{gomoyunov2021viscosity}]\label{lemma_5_6_g21}
Let \(K\subset AC^\alpha([0,T],\mathbb R^\ell)\) be compact. For every
\(\varepsilon>0\), let
\[
    (t_\varepsilon,\gamma_\varepsilon),
    (\tau_\varepsilon,\nu_\varepsilon)
    \in [0,T]\times K
\]
be such that
\[
    \varpi_\varepsilon
    (t_\varepsilon,\gamma_\varepsilon,
    \tau_\varepsilon,\nu_\varepsilon)
    \to 0
    \qquad \text{as } \varepsilon\to0^+.
\]
Then
\[
\left\|
a(\cdot\mid t_\varepsilon,\gamma_\varepsilon)
-
a(\cdot\mid \tau_\varepsilon,\nu_\varepsilon)
\right\|_\infty
\to0 .
\]
Additionally suppose that, for some \(z>0\) and some $k > 0$,
\[
    t_\epsilon,\tau_\epsilon\in[0,T-z],
    \qquad
    \gamma^\epsilon,\nu^\epsilon\in AC_k^\alpha,
\]
for every sufficiently small \(\epsilon>0\). Assume further that
\(
    \varpi_\epsilon
    (t_\epsilon,\gamma^\epsilon,\tau_\epsilon,\nu^\epsilon)
    \leq C\epsilon^2.
\)
Then
\[
\sup_{r\in[0,T-z]}
\left|
a(r\mid t_\epsilon,\gamma^\epsilon)
-
a(r\mid\tau_\epsilon,\nu^\epsilon)
\right|
\leq
C_z\epsilon^{\frac{2\alpha}{1+\alpha c}}
=
C_z\epsilon^{\frac{2\alpha(2-\alpha)}{2+\alpha}}.
\]
If, in addition,
\(
    |t_\epsilon-\tau_\epsilon|
    \leq C\epsilon^{\frac{3}{2\alpha}},
\)
then
\(
    \left|
    \gamma^\epsilon_{t_\epsilon}
    -
    \nu^\epsilon_{\tau_\epsilon}
    \right|
    \leq
    C_z\epsilon^{\frac{2\alpha}{1+\alpha c}}.
\)
\end{lemma}

\begin{proof}
The first assertion follows from minor adaptations to Lemma 5.2 in \cite{gomoyunov2021viscosity} so we omit it. For the second claim, set
\[
    F_\epsilon(r)
    :=
    a(r\mid t_\epsilon,\gamma^\epsilon)
    -
    a(r\mid\tau_\epsilon,\nu^\epsilon),
    \qquad r\in[0,T].
\]
Since the map \(a(\cdot\mid s,\gamma)\) preserves \(AC_k^\alpha\) (as recalled after equation \eqref{auxiliary_functional}),
there exists a constant \(L_k>0\), independent of
\(\epsilon,t_\epsilon,\tau_\epsilon,\gamma^\epsilon,\nu^\epsilon\),
such that
\[
    |F_\epsilon(r)-F_\epsilon(u)|
    \leq L_k|r-u|^\alpha,
    \qquad r,u\in[0,T].
\]

By the definition of \(\varpi_\epsilon\),
\begin{align*}
&\int_0^T
\frac{|F_\epsilon(r)|^c}
     {(T-r)^{(1-\alpha-\beta)c}}
\,dr \leq
\varpi_\epsilon
(t_\epsilon,\gamma^\epsilon,\tau_\epsilon,\nu^\epsilon)
+
C_1\epsilon^{\frac{c}{c-1}}.
\end{align*}
Since
\[
    \frac{c}{c-1}=\frac{2}{\alpha}>2
\]
and \(0<\epsilon<1\), it follows that
\[
\int_0^T
\frac{|F_\epsilon(r)|^c}
     {(T-r)^{(1-\alpha-\beta)c}}
\,dr
\leq C\epsilon^2.
\]

Let
\[
    h_\epsilon
    :=
    \sup_{r\in[0,T-z]}|F_\epsilon(r)|
\]
and choose \(r_\epsilon\in[0,T-z]\) such that
\(
    |F_\epsilon(r_\epsilon)|=h_\epsilon.
\)
We may assume that \(h_\epsilon>0\), since otherwise there is
nothing to prove.

For sufficiently small \(\epsilon\), define
\[
    \delta_\epsilon
    :=
    \left(\frac{h_\epsilon}{2L_k}\right)^{1/\alpha}.
\]
Since the weighted integral above converges to zero as $\epsilon \to 0$, the uniform
Hölder estimate implies that \(h_\epsilon\to0\). Hence, after
reducing \(\epsilon\) if necessary, we have
\(
    \delta_\epsilon\leq z/2.
\)
For every
\(
    r\in[r_\epsilon,r_\epsilon+\delta_\epsilon],
\)
the Hölder estimate gives
\[
\begin{aligned}
    |F_\epsilon(r)|
    &\geq
    |F_\epsilon(r_\epsilon)|
    -
    |F_\epsilon(r)-F_\epsilon(r_\epsilon)| \geq
    h_\epsilon-L_k\delta_\epsilon^\alpha
    =
    \frac{h_\epsilon}{2}.
\end{aligned}
\]
Moreover,
\(
    r_\epsilon+\delta_\epsilon
    \leq T-\frac{z}{2}.
\)
The weight
\[
    r\to
    \frac{1}{(T-r)^{(1-\alpha-\beta)c}}
\]
is continuous and strictly positive on
\([0,T-z/2]\). Consequently, there exists \(m_z>0\) such that
\[
    \frac{1}{(T-r)^{(1-\alpha-\beta)c}}
    \geq m_z,
    \qquad r\in[0,T-z/2].
\]
It follows that
\begin{align*}
    C\epsilon^2
    &\geq
    \int_0^T
    \frac{|F_\epsilon(r)|^c}
         {(T-r)^{(1-\alpha-\beta)c}}
    \,dr\\
    &\geq
    m_z
    \int_{r_\epsilon}^{r_\epsilon+\delta_\epsilon}
    |F_\epsilon(r)|^c\,dr\\
    &\geq
    m_z\delta_\epsilon
    \left(\frac{h_\epsilon}{2}\right)^c\\
    &=
    \frac{m_z}{2^c(2L_k)^{1/\alpha}}
    h_\epsilon^{c+1/\alpha}.
\end{align*}
Therefore
\(
    h_\epsilon
    \leq
    C_z\epsilon^{\frac{2}{c+1/\alpha}}
    =
    C_z\epsilon^{\frac{2\alpha}{1+\alpha c}},
\)
which proves the second assertion.

Finally, recalling that
\[
    \gamma^\epsilon_{t_\epsilon}
    =
    a(t_\epsilon\mid t_\epsilon,\gamma^\epsilon),
    \qquad
    \nu^\epsilon_{\tau_\epsilon}
    =
    a(\tau_\epsilon\mid\tau_\epsilon,\nu^\epsilon),
\]
we obtain
\begin{align*}
\left|
\gamma^\epsilon_{t_\epsilon}
-
\nu^\epsilon_{\tau_\epsilon}
\right|
&\leq
\left|
a(t_\epsilon\mid t_\epsilon,\gamma^\epsilon)
-
a(t_\epsilon\mid\tau_\epsilon,\nu^\epsilon)
\right|\\
&\quad+
\left|
a(t_\epsilon\mid\tau_\epsilon,\nu^\epsilon)
-
a(\tau_\epsilon\mid\tau_\epsilon,\nu^\epsilon)
\right|\\
&\leq
C_z\epsilon^{\frac{2\alpha}{1+\alpha c}}
+
L_k|t_\epsilon-\tau_\epsilon|^\alpha\\
&\leq
C_z\epsilon^{\frac{2\alpha}{1+\alpha c}}
+
C\epsilon^{3/2}.
\end{align*}
Since
\[
    \frac{2\alpha}{1+\alpha c}<1<\frac32,
\]
the second term can be absorbed into the first, proving the result.
\end{proof}

\begin{lemma}[Adapted from Lemma 5.4 in \cite{gomoyunov2021viscosity}]\label{lemma_5_7_g21}
Let  $\varpi_{\epsilon}$, $C_1$  and $\beta$ be defined as in \eqref{auxiliary_functional}.
For any $\epsilon > 0$ and any fixed
$(\tau,\nu) \in [0,T] \times AC^\alpha([0, T], \mathbb{R}^\ell)$, the functional
$\mu_\epsilon^{(\tau,\nu)}(t, \gamma):= \varpi_{\epsilon}(t, \gamma, \tau, \nu)$ is ci-smooth of the order $\alpha$
and its ci-derivatives of the order $\alpha$ are given by
\begin{equation*}
\partial_t^\alpha
\mu_\epsilon^{(\tau,\nu)}(t,\gamma) = 0, \quad \partial_x^\alpha
\mu_\epsilon^{(\tau,\nu)}(t,\gamma) = 0 
\end{equation*}
and
\begin{align*}
&\nabla^\alpha
\mu_\epsilon^{(\tau,\nu)}(t,\gamma)
\nonumber \\
&\quad =
\frac{c}{\Gamma(\alpha)}
\left(
\frac{
a\bigl(T \mid t,\gamma\bigr) - a\bigl(T \mid \tau ,\nu\bigr)
}{
\left(
\epsilon^{\frac{2}{c-1}}
+
\left|a\bigl(T \mid t,\gamma\bigr)-a\bigl(T \mid \tau ,\nu\bigr)\right|^2
\right)^{1-\frac{c}{2}}
(T-t)^{1-\alpha}
}
\right.
\nonumber \\
&\qquad\qquad\left.
+
\int_t^T
\frac{
a\bigl(r \mid t,\gamma \bigr) - a\bigl(r \mid \tau,\nu \bigr)
}{
\left(
\epsilon^{\frac{2}{c-1}}
+
\left|a\bigl(r \mid t,\gamma\bigr)-a\bigl(r \mid \tau,\nu \bigr)\right|^2
\right)^{1-\frac{c}{2}}
(r-t)^{1-\alpha}
(T-r)^{(1-\alpha-\beta)c}
}
\,dr
\right)
\end{align*}
for all $(t,\gamma) \in [0,T) \times AC^\alpha([0, T], \mathbb{R}^\ell)$.
\end{lemma}
\begin{proof}
Although our definition of $\varpi_{\epsilon}$ differs from the one in \cite{gomoyunov2021viscosity}, in the proof of Lemma 5.4 the author shows that for every $\tilde t \in (t, T)$ and every path $\tilde \gamma \in AC^{\alpha}([0, T], \mathbb{R}^\ell)$ such that $\gamma(s) = \tilde \gamma (s)$ for $s \in [0, t]$, then there exists a constant $C>0$ for which
\begin{align*}
&\left|
\left(
\epsilon^{\frac{2}{c-1}}
+
\left|a\bigl(T \mid \tilde t ,\tilde \gamma \bigr)-a\bigl(T \mid \tau ,\nu \bigr)\right|^2
\right)^{\frac{c}{2}}
-
\left(
\epsilon^{\frac{2}{c-1}}
+
\left|a\bigl(T \mid t,\gamma\bigr)-a\bigl(T \mid \tau ,\nu \bigr) \right|^2
\right)^{\frac{c}{2}}
\right.
\nonumber \\
&\qquad\left.
-
\frac{
c\left\langle
a\bigl(T \mid t,\gamma \bigr)-a\bigl(T \mid \tau ,\nu \bigr), (I^{\alpha}_{0^+}(\tilde \gamma - \gamma_0)(\tilde t ) - I^{\alpha}_{0^+}(\gamma - \gamma_0)(t )) 
\right\rangle
}{
\Gamma(\alpha)
\left(
\epsilon^{\frac{2}{c-1}}
+
\left|a\bigl(T \mid t,\gamma\bigr)-a\bigl(T \mid \tau ,\nu \bigr)\right|^2
\right)^{1-\frac{c}{2}}
(T-t)^{1-\alpha}
}
\right|
\leq C(\tilde t - t)^2,
\end{align*}

and

\begin{align*}
&\left|
\int_{t}^T
\frac{
\left(
\epsilon^{\frac{2}{c-1}}
+
\left|a\bigl(r \mid \tilde t, \tilde \gamma\bigr)-a\bigl(r \mid \tau ,\nu\bigr)\right|^2
\right)^{\frac{c}{2}}
-
\left(
\epsilon^{\frac{2}{c-1}}
+
\left|a\bigl(r \mid t,\gamma\bigr)-a\bigl(r \mid \tau ,\nu\bigr) \right|^2
\right)^{\frac{c}{2}}
}{
(T-r)^{(1-\alpha-\beta)c}
}
\,dr
\right.
\nonumber \\
&\qquad\left.
-
\int_t^T
\frac{
c\left\langle
a\bigl(r \mid t,\gamma\bigr)-a\bigl(r \mid \tau ,\nu\bigr),(I^{\alpha}_{0^+}(\tilde \gamma - \gamma_0)(\tilde t ) - I^{\alpha}_{0^+}(\gamma - \gamma_0)(t )
\right\rangle
}{
\Gamma(\alpha)
\left(
\epsilon^{\frac{2}{c-1}}
+
\left|a\bigl(r \mid t,\gamma \bigr)-a\bigl(r \mid \tau ,\nu\bigr)\right|^2
\right)^{1-\frac{c}{2}}
(r-t)^{1-\alpha}
(T-r)^{(1-\alpha-\beta)c}
}
\,dr
\right|
\leq C(\tilde t - t)^{1+\alpha}.
\end{align*}
Comparing the above with Definition \ref{ci_smooth_definition} is enough to conclude our proof.
\end{proof}

\begin{lemma}[Adapted from Lemma 5.5 in \cite{gomoyunov2021viscosity}]\label{lemma_5_8_g21}
Let $\varpi_{\epsilon}$, $c$, $C_1$ and $\mu^{\cdot, \cdot}$ be defined as above, then there is a constant $C > 0$ such that
\begin{equation*}
\left|
\nabla^\alpha \mu_\epsilon^{(\tau,\nu)}(t,\gamma)
\right|
\leq
C
\left(
\varpi_\epsilon\bigl(t,\gamma,\tau,\nu\bigr)
+
C_1\epsilon^{\frac{c}{c-1}}
\right)^{\frac{c-1}{c}}
\end{equation*}
and
\begin{align*}
&\left|
\nabla^\alpha \mu_\epsilon^{(\tau,\nu)}(t,\gamma)
+
\nabla^\alpha \mu_\epsilon^{(t,\gamma)}(\tau,\nu)
\right|
\nonumber \\
&\qquad\leq
C
\left(
\epsilon^{\frac{2}{c-1}}
+
\left\|
a\bigl(\cdot \mid t,\gamma\bigr)
-
a\bigl(\cdot \mid \tau, \nu\bigr)
\right\|_\infty^2
\right)^{\frac{c-1}{2}}
|t-\tau|^\alpha
\end{align*}
for any $\epsilon > 0$ and any
$(t,\gamma),(\tau,\nu) \in [0,T) \times AC^\alpha([0, T], \mathbb{R}^\ell)$.
\end{lemma}
\begin{proof}
Although the definitions of $\varpi_{\epsilon}$ and
$\mu_{\epsilon}^{\cdot,\cdot}$ used here differ from those in
\cite{gomoyunov2021viscosity}, the explicit formula obtained above for
$\nabla^\alpha \mu_\epsilon^{(\tau,\nu)}$ coincides with the
corresponding expression in that work. Consequently, the desired estimates
follow directly from the proof of Lemma 5.5 in
\cite{gomoyunov2021viscosity}.
\end{proof}

\end{appendices}

\bibliographystyle{alpha}
\bibliography{biblio}

@article{cass2022combinatorial,
  title={A combinatorial approach to geometric rough paths and their controlled paths},
  author={Cass, Thomas and Driver, Bruce K and Litterer, Christian and Rossi Ferrucci, Emilio},
  journal={Journal of the London Mathematical Society},
  volume={106},
  number={2},
  pages={936--981},
  year={2022},
  publisher={Wiley Online Library}
}

@article{allan2020pathwise,
  author  = {Allan, Andrew L. and Cohen, Samuel N.},
  title   = {Pathwise stochastic control with applications to robust filtering},
  journal = {The Annals of Applied Probability},
  year    = {2020},
  volume  = {30},
  number  = {5},
  pages   = {2274--2310},
  doi     = {10.1214/19-AAP1558},
  url     = {https://doi.org/10.1214/19-AAP1558}
}

@article{friz2018differential,
  title={Differential equations driven by rough paths with jumps},
  author={Friz, Peter K and Zhang, Huilin},
  journal={Journal of Differential Equations},
  volume={264},
  number={10},
  pages={6226--6301},
  year={2018},
  publisher={Elsevier}
}

@book{samko1993fractional,
  title={Fractional Integrals and Derivatives},
  author={Samko, S. and Kilbas, A.A. and Marichev, O.},
  isbn={9782881248641},
  lccn={93026071},
  url={https://books.google.it/books?id=SO3FQgAACAAJ},
  year={1993},
  publisher={Taylor \& Francis}
}

@article{gomoyunov2020dynamic,
  title={Dynamic Programming Principle and Hamilton--Jacobi--Bellman Equations for Fractional-Order Systems},
  author={Gomoyunov, Mikhail I},
  journal={SIAM Journal on Control and Optimization},
  volume={58},
  number={6},
  pages={3185--3211},
  year={2020},
  publisher={SIAM}
}

@book{bardi1997optimal,
  author    = {Bardi, Martino and Capuzzo-Dolcetta, Italo},
  title     = {Optimal Control and Viscosity Solutions of Hamilton-Jacobi-Bellman Equations},
  publisher = {Birkh{\"a}user Boston},
  address   = {Boston, MA},
  year      = {1997},
  pages     = {xviii+574},
  doi       = {10.1007/978-0-8176-4755-1},
  url       = {https://doi.org/10.1007/978-0-8176-4755-1},
}

@article{bardi1997bellman,
  title={On the Bellman equation for some unbounded control problems},
  author={Bardi, Martino and Da Lio, Francesca},
  journal={Nonlinear Differential Equations and Applications NoDEA},
  volume={4},
  pages={491--510},
  year={1997},
  publisher={Springer}
}

@article{gomoyunov2021viscosity,
  title={On viscosity solutions of path-dependent Hamilton--Jacobi--Bellman--Isaacs equations for fractional-order systems},
  author={Gomoyunov, Mikhail I},
  journal={Journal of Differential Equations},
  volume={399},
  pages={335--362},
  year={2024},
  publisher={Elsevier}
}

@article{lukoyanov2007viscosity,
  title={On viscosity solution of functional Hamilton-Jacobi type equations for hereditary systems},
  author={Lukoyanov, N Yu},
  journal={Proceedings of the Steklov Institute of Mathematics},
  volume={259},
  pages={S190--S200},
  year={2007},
  publisher={Springer}
}

@article{lions1985differential,
  title={Differential games, optimal control and directional derivatives of viscosity solutions of Bellman’s and Isaacs’ equations},
  author={Lions, P-L and Souganidis, Panagiotis E},
  journal={SIAM journal on control and optimization},
  volume={23},
  number={4},
  pages={566--583},
  year={1985},
  publisher={SIAM}
}

@book{kim1999functional,
  author    = {Kim, A. V.},
  title     = {Functional Differential Equations: Application of i-Smooth Calculus},
  series    = {Mathematics and Its Applications},
  volume    = {479},
  publisher = {Springer},
  address   = {Dordrecht},
  year      = {1999},
  doi       = {10.1007/978-94-017-1630-7},
  url       = {https://doi.org/10.1007/978-94-017-1630-7},
  isbn      = {978-0-7923-5689-9}
}

@book{achdou2013hamilton,
  title     = {Hamilton-Jacobi Equations: Approximations, Numerical Analysis and Applications},
  author    = {Achdou, Yves and Barles, Guy and Ishii, Hitoshi and Litvinov, Grigory L.},
  series    = {Lecture Notes in Mathematics},
  volume    = {2074},
  publisher = {Springer},
  address   = {Berlin, Heidelberg},
  year      = {2013},
  doi       = {10.1007/978-3-642-36433-4},
  isbn      = {978-3-642-36433-4}
}

@article{gomoyunov2020theory,
  title={To the theory of differential inclusions with Caputo fractional derivatives},
  author={Gomoyunov, MI},
  journal={Differential Equations},
  volume={56},
  pages={1387--1401},
  year={2020},
  publisher={Springer}
}

@book{diethelm2010analysis,
  author    = {Diethelm, Kai},
  title     = {The Analysis of Fractional Differential Equations: An Application-Oriented Exposition Using Differential Operators of Caputo Type},
  series    = {Lecture Notes in Mathematics},
  volume    = {2004},
  publisher = {Springer},
  address   = {Berlin, Heidelberg},
  year      = {2010},
  doi       = {10.1007/978-3-642-14574-2},
  url       = {https://doi.org/10.1007/978-3-642-14574-2},
  isbn      = {978-3-642-14573-5}
}

@article{allan2019parameter,
  title={Parameter uncertainty in the Kalman--Bucy filter},
  author={Allan, Andrew L and Cohen, Samuel N},
  journal={SIAM Journal on Control and Optimization},
  volume={57},
  number={3},
  pages={1646--1671},
  year={2019},
  publisher={SIAM}
}

@inproceedings{caruana2011rough,
  title={A (rough) pathwise approach to a class of non-linear stochastic partial differential equations},
  author={Caruana, Michael and Friz, Peter K and Oberhauser, Harald},
  booktitle={Annales de l'Institut Henri Poincar{\'e} C, Analyse non lin{\'e}aire},
  volume={28},
  pages={27--46},
  year={2011},
  organization={Elsevier}
}

@book{friz2020course,
  author    = {Friz, Peter K. and Hairer, Martin},
  title     = {A Course on Rough Paths: With an Introduction to Regularity Structures},
  series    = {Universitext},
  edition   = {2},
  publisher = {Springer},
  address   = {Cham},
  year      = {2020},
  doi       = {10.1007/978-3-030-41556-3},
  url       = {https://doi.org/10.1007/978-3-030-41556-3},
  isbn      = {978-3-030-41555-6}
}

@article{ha1992deterministic,
  title={A deterministic approach to stochastic optimal control with application to anticipative control},
  author={HA Davis, Mark and Burstein, Gabriel},
  journal={Stochastics: An International Journal of Probability and Stochastic Processes},
  volume={40},
  number={3-4},
  pages={203--256},
  year={1992},
  publisher={Taylor \& Francis}
}

@inproceedings{ocone1989generalized,
  title={A generalized It{\^o}-Ventzell formula. Application to a class of anticipating stochastic differential equations},
  author={Ocone, Daniel and Pardoux, {\'E}tienne},
  booktitle={Annales de l'IHP Probabilit{\'e}s et statistiques},
  volume={25},
  pages={39--71},
  year={1989}
}

@article{lions1998fully,
  title={Fully nonlinear stochastic partial differential equations: non-smooth equations and applications},
  author={Lions, Pierre-Louis and Souganidis, Panagiotis E},
  journal={Comptes Rendus de l'Acad{\'e}mie des Sciences-Series I-Mathematics},
  volume={327},
  number={8},
  pages={735--741},
  year={1998},
  publisher={Elsevier}
}

@article{buckdahn2007pathwise,
  title={Pathwise stochastic control problems and stochastic HJB equations},
  author={Buckdahn, Rainer and Ma, Jin},
  journal={SIAM journal on control and optimization},
  volume={45},
  number={6},
  pages={2224--2256},
  year={2007},
  publisher={SIAM}
}

@article{diehl2017stochastic,
  title={Stochastic control with rough paths},
  author={Diehl, Joscha and Friz, Peter K and Gassiat, Paul},
  journal={Applied Mathematics \& Optimization},
  volume={75},
  pages={285--315},
  year={2017},
  publisher={Springer}
}

@article{rogers2002monte,
  title={Monte Carlo valuation of American options},
  author={Rogers, Leonard CG},
  journal={Mathematical Finance},
  volume={12},
  number={3},
  pages={271--286},
  year={2002},
  publisher={Wiley Online Library}
}

@article{rogers2007pathwise,
  title={Pathwise stochastic optimal control},
  author={Rogers, LCG},
  journal={SIAM Journal on Control and Optimization},
  volume={46},
  number={3},
  pages={1116--1132},
  year={2007},
  publisher={SIAM}
}

@inproceedings{wets1975relation,
  title={On the relation between stochastic and deterministic optimization},
  author={Wets, Roger J-B},
  booktitle={Control Theory, Numerical Methods and Computer Systems Modelling: International Symposium, Rocquencourt, June 17--21, 1974},
  pages={350--361},
  year={1975},
  organization={Springer}
}

@article{hairer2015geometric,
  author  = {Hairer, Martin and Kelly, David},
  title   = {Geometric versus non-geometric rough paths},
  journal = {Annales de l'I.H.P. Probabilit{\'e}s et statistiques},
  year    = {2015},
  volume  = {51},
  number  = {1},
  pages   = {207--251},
  doi     = {10.1214/13-AIHP564},
  url     = {https://www.numdam.org/articles/10.1214/13-AIHP564/}
}

@article{friz2010differential,
  author  = {Friz, Peter and Victoir, Nicolas},
  title   = {Differential equations driven by Gaussian signals},
  journal = {Annales de l'I.H.P. Probabilit{\'e}s et statistiques},
  year    = {2010},
  volume  = {46},
  number  = {2},
  pages   = {369--413},
  doi     = {10.1214/09-AIHP202},
  url     = {https://www.numdam.org/articles/10.1214/09-AIHP202/}
}

@article{bucci2018slow,
  title={Slow decay of impact in equity markets: insights from the ANcerno database},
  author={Bucci, Fr{\'e}d{\'e}ric and Benzaquen, Michael and Lillo, Fabrizio and Bouchaud, Jean-Philippe},
  journal={Market Microstructure and Liquidity},
  volume={4},
  number={03n04},
  pages={1950006},
  year={2018},
  publisher={World Scientific}
}

@book{webster2023handbook,
  title={Handbook of price impact modeling},
  author={Webster, Kevin T},
  year={2023},
  publisher={Chapman and Hall/CRC}
}

@article{bouchaud2003fluctuations,
author = {Jean-Philippe Bouchaud and Yuval Gefen and Marc Potters and Matthieu Wyart},
title = {Fluctuations and response in financial markets: the subtle nature of ‘random’ price changes},
journal = {Quantitative Finance},
volume = {4},
number = {2},
pages = {176--190},
year = {2004},
publisher = {Routledge},
doi = {10.1080/14697680400000022},
URL = {  
        https://www.tandfonline.com/doi/abs/10.1080/14697680400000022
},
eprint = { 
        https://www.tandfonline.com/doi/pdf/10.1080/14697680400000022
}
}

@article{abi2022optimal,
  author  = {Abi Jaber, Eduardo and Neuman, Eyal},
  title   = {Optimal Liquidation With Signals: The General Propagator Case},
  journal = {Mathematical Finance},
  year    = {2025},
  volume  = {35},
  number  = {4},
  pages   = {841--866},
  doi     = {10.1111/mafi.12465},
  url     = {https://doi.org/10.1111/mafi.12465}
}

@article{harang2021volterra,
  title={Volterra equations driven by rough signals},
  author={Harang, Fabian A and Tindel, Samy},
  journal={Stochastic Processes and their Applications},
  volume={142},
  pages={34--78},
  year={2021},
  publisher={Elsevier}
}

@article{gatheral2010no,
  title={No-dynamic-arbitrage and market impact},
  author={Gatheral, Jim},
  journal={Quantitative finance},
  volume={10},
  number={7},
  pages={749--759},
  year={2010},
  publisher={Taylor \& Francis}
}

@article{bel23,
author = {Carlo Bellingeri},
title = {{Quasi-geometric rough paths and rough change of variable formula}},
volume = {59},
journal = {Annales de l'Institut Henri Poincaré, Probabilités et Statistiques},
number = {3},
publisher = {Institut Henri Poincaré},
pages = {1398 -- 1433},
year = {2023},
doi = {10.1214/22-AIHP1297},
URL = {https://doi.org/10.1214/22-AIHP1297}
}

@article{gomoyunov2021path,
  author  = {Gomoyunov, Mikhail I. and Lukoyanov, Nikolai Yu. and Plaksin, Anton R.},
  title   = {Path-Dependent Hamilton--Jacobi Equations: The Minimax Solutions Revised},
  journal = {Applied Mathematics \& Optimization},
  volume  = {84},
  number  = {1},
  pages   = {1087--1117},
  year    = {2021},
  doi     = {10.1007/s00245-021-09794-4},
  url     = {https://doi.org/10.1007/s00245-021-09794-4}
}

@article{BAYRAKTAR20182096,
title = {Path-dependent Hamilton–Jacobi equations in infinite dimensions},
journal = {Journal of Functional Analysis},
volume = {275},
number = {8},
pages = {2096-2161},
year = {2018},
issn = {0022-1236},
doi = {https://doi.org/10.1016/j.jfa.2018.07.010},
url = {https://www.sciencedirect.com/science/article/pii/S002212361830274X},
author = {Erhan Bayraktar and Christian Keller},
}

@article{bayraktar2025viscosity,
  title={Viscosity and minimax solutions for path-dependent Hamilton-Jacobi equations in infinite dimensions and related differential games},
  author={Bayraktar, Erhan and Gomoyunov, Mikhail and Keller, Christian},
  journal={arXiv preprint arXiv:2509.16015},
  year={2025}
}

@article{brown2010information,
  title={Information relaxations and duality in stochastic dynamic programs},
  author={Brown, David B and Smith, James E and Sun, Peng},
  journal={Operations research},
  volume={58},
  number={4-part-1},
  pages={785--801},
  year={2010},
  publisher={INFORMS}
}

@article{DB09,
author = {Dumitru Baleanu and Ozlem Defterli and Om P. Agrawal},
title ={A Central Difference Numerical Scheme for Fractional Optimal Control Problems},

journal = {Journal of Vibration and Control},
volume = {15},
number = {4},
pages = {583-597},
year = {2009},
doi = {10.1177/1077546308088565},

URL = { 
    
        https://doi.org/10.1177/1077546308088565
    
    

},
}

@misc{friz2021rough,
  title         = {Rough stochastic differential equations},
  author        = {Friz, Peter K. and Hocquet, Antoine and L{\^e}, Khoa},
  year          = {2021},
  eprint        = {2106.10340},
  archivePrefix = {arXiv},
  primaryClass  = {math.PR},
  doi           = {10.48550/arXiv.2106.10340}
}

@misc{friz2024controlled,
  title         = {Controlled rough SDEs, pathwise stochastic control and dynamic programming principles},
  author        = {Friz, Peter K. and L{\^e}, Khoa and Zhang, Huilin},
  year          = {2024},
  eprint        = {2412.05698},
  archivePrefix = {arXiv},
  primaryClass  = {math.PR},
  doi           = {10.48550/arXiv.2412.05698}
}

@article{gomoyunov2019approximation,
  author  = {Gomoyunov, Mikhail},
  title   = {Approximation of Fractional Order Conflict-Controlled Systems},
  journal = {Progress in Fractional Differentiation and Applications},
  volume  = {5},
  number  = {2},
  pages   = {143--155},
  year    = {2019},
  doi     = {10.18576/pfda/050205}
}

@article{henrylabordere2016dual,
author = {Henry-Labord{\`e}re, Pierre and Litterer, Christian and Ren, Zhenjie},
title = {A Dual Algorithm for Stochastic Control Problems: Applications to Uncertain Volatility Models and CVA},
journal = {SIAM Journal on Financial Mathematics},
volume = {7},
number = {1},
pages = {159-182},
year = {2016},
doi = {10.1137/15M1019945},
URL = { 
        https://doi.org/10.1137/15M1019945
},
eprint = { 
        https://doi.org/10.1137/15M1019945
}
}

@misc{bank2026duality,
      title={Duality methods in stochastic optimal control}, 
      author={Peter Bank and Filippo de Feo},
      year={2026},
      eprint={2602.17823},
      archivePrefix={arXiv},
      primaryClass={math.OC},
      url={https://arxiv.org/abs/2602.17823}, 
}

@misc{lauriere2026dual,
      title={Dual Approaches to Stochastic Control via SPDEs and the Pathwise Hopf Formula}, 
      author={Mathieu Laurière and Jiefei Yang},
      year={2026},
      eprint={2604.08155},
      archivePrefix={arXiv},
      primaryClass={math.OC},
      url={https://arxiv.org/abs/2604.08155}, 
}

@article{carbotti2021note,
  author  = {Carbotti, Alessandro and Comi, Giovanni E.},
  title   = {A Note on Riemann--Liouville Fractional Sobolev Spaces},
  journal = {Communications on Pure and Applied Analysis},
  volume  = {20},
  number  = {1},
  pages   = {17--54},
  year    = {2021},
  doi     = {10.3934/cpaa.2020255}
}
\nocite{*}
\addcontentsline{toc}{chapter}{Bibliography}
\end{document}